\documentclass{amsart}
\usepackage{amssymb, latexsym}

\def\bi{\begin{itemize}}
\def\bs{\begin{split}}
\def\es{\end{split}}
\def\ba{\begin{align}}
\def\bas{\begin{align*}}
\def\ea{\end{align}}
\def\eas{\end{align*}}
\def\Im{{\operatorname{Im}}}

\def\Re{{\operatorname{Re}}}

\def\C{{\mathbb C}} 
\def\R{{{\mathbb R}}}
\def\Z{{{\mathbb Z}}}

\def\emph#1{{\it #1}}
\def\textbf#1{{\bf #1}}

\newcommand{\lqlr}{{L^q_t L^r_x}}
\newcommand{\ir}{{I \times \R^n}}
\newcommand{\izr}{{I_0\times \R^n}}
\newcommand{\nizr}{{\dot N^0(\izr)}}
\newcommand{\rr}{{\R \times \R^n}}
\newcommand{\lli}{{\lqlr (\ir)}}
\newcommand{\llr}{{\lqlr (\rr)}}
\newcommand{\sz}{{\dot S^0}}
\newcommand{\szi}{{\sz (\ir)}}
\newcommand{\szr}{{\sz (\rr)}}
\newcommand{\so}{{\dot S^1}}

\newcommand{\sor}{{\so (\rr)}}
\newcommand{\soizr}{{\so(\izr)}}

\newcommand{\sk}{{\dot S^k}}
\newcommand{\ski}{{\sk (\ir)}}

\newcommand{\hk}{{\dot H^k}}
\newcommand{\ho}{{\dot H^1}}

\newcommand{\nmt}{{\frac 4{n-2}}}

\newcommand{\dnnt}{{\frac {2n}{n-2}}}
\newcommand{\ijr}{{I_j\times \R^n}}
\newcommand{\hox}{{\dot H^1_x}}
\newcommand{\xoio}{{\dot X^1(I_1)}}

\newcommand{\sntnt}{{\dot S^0([t_0,T]\times\R^n)}}

\newcommand{\sftnt}{{\dot S^1([t_0,T]\times\R^n)}}

\newcommand{\nntnt}{{\dot N^0([t_0,T]\times\R^n)}}
\newcommand{\nnjn}{{\dot N^0(J_0\times\R^n)}}
\newcommand{\nfjn}{{\dot N^1(J_0\times\R^n)}}
\newcommand{\nnij}{{\dot N^0(I_j\times\R^n)}}
\newcommand{\sfjn}{{\dot S^1(J_0\times\R^n)}}
\newcommand{\snjn}{{\dot S^0(J_0\times\R^n)}}

\newcommand{\snij}{{\dot S^0(I_j\times\R^n)}}

\newcommand{\eps}{{\varepsilon}}
\newcommand{\uhi}{{u_{hi}}}
\newcommand{\ulo}{{u_{lo}}}

\newcommand{\pl}{{P_{lo}}}
\newcommand{\ph}{{P_{hi}}}

\newcommand{\propagateo}{{e^{i(t-t_1)\Delta}}}
\newcommand{\ntu}{{\lambda_1|u|^{p_1}u+\lambda_2|u|^{p_2}u}}
\newcommand{\ntv}{{\lambda_1|v|^{p_1}v+\lambda_2|v|^{p_2}v}}

\theoremstyle{plain}
\newtheorem{theorem}{Theorem}
\newtheorem{definition}[theorem]{Definition}
\newtheorem{remark}[theorem]{Remark}
\newtheorem{proposition}[theorem]{Proposition}
\newtheorem{lemma}[theorem]{Lemma}

\newtheorem{assumption}[theorem]{Assumption}

\numberwithin{equation}{section} \numberwithin{theorem}{section}

\begin{document}

\title[NLS with combined power-type nonlinearities]
{The nonlinear Schr\"odinger equation with combined power-type nonlinearities}
\author{Terence Tao}
\address{University of California, Los Angeles}
\author{Monica Visan}
\address{University of California, Los Angeles}
\author{Xiaoyi Zhang}
\thanks{The third author was supported in part by NSF grant No. 10426029.}
\address{Academy of Mathematics and System Sciences, Chinese Academy of Sciences}

\vspace{-0.3in}
\begin{abstract}
We undertake a comprehensive study of the nonlinear Schr\"odinger equation
$$
i u_t +\Delta u = \lambda_1|u|^{p_1} u+ \lambda_2 |u|^{p_2} u,
$$
where $u(t,x)$ is a complex-valued function in spacetime $\R_t\times\R^n_x$, $\lambda_1$ and $\lambda_2$
are nonzero real constants, and $0<p_1<p_2\le \frac 4{n-2}$. We address questions related to local and
global well-posedness, finite time blowup, and asymptotic behaviour.  Scattering is considered both in
the energy space $H^1(\R^n)$ and in the pseudoconformal
space $\Sigma:=\{f\in H^1(\R^n); \, xf\in L^2(\R^n)\}$.  Of particular interest is the case
when both nonlinearities are defocusing and correspond to the $L_x^2$-critical, respectively
$\dot H^1_x$-critical NLS, that is, $\lambda_1, \lambda_2>0$ and $p_1=\frac{4}{n}$, $p_2=\frac{4}{n-2}$.
The results at the endpoint $p_1 = \frac{4}{n}$ are conditional on a conjectured global existence and spacetime estimate
for the $L^2_x$-critical nonlinear Schr\"odinger equation.

As an off-shoot of our analysis, we also obtain a new, simpler proof of scattering in $H^1_x$ for solutions to the
nonlinear Schr\"odinger equation
$$
i u_t +\Delta u = |u|^{p} u,
$$
with $\frac{4}{n}<p<\frac{4}{n-2}$, which was first obtained by J. Ginibre and G. Velo, \cite{gv:scatter}.
\end{abstract}

\maketitle

\tableofcontents

\section{Introduction}

We study the initial value problem for the nonlinear Schr\"odinger equation with two power-type nonlinearities,
\begin{equation}\label{equation}
\begin{cases}
i u_t +\Delta u = \lambda_1|u|^{p_1} u+ \lambda_2 |u|^{p_2} u\\
u(0,x) = u_0(x),
\end{cases}
\end{equation}
where $u(t,x)$ is a complex-valued function in spacetime $\R_t\times \R^n_x$, $n\geq 3$, the initial data $u_0$
belongs to ${H}^1_x$ (or $\Sigma$), $\lambda_1$, $\lambda_2$ are nonzero real constants, and $0<p_1<p_2\le \frac 4{n-2}$.

This equation has Hamiltonian
\begin{equation}\label{energy}
E(u(t)):=\int_{\R^n} \Bigl[\tfrac{1}{2}|\nabla u (t,x)|^2
+\tfrac{\lambda_1}{p_1+2}|u(t,x)|^{p_1+2}+\tfrac{\lambda_2}{p_2+2}|u(t,x)|^{p_2+2}\Bigr]dx.
\end{equation}
As \eqref{energy} is preserved\footnote{To justify the energy conservation rigorously, one can approximate the data $u_0$ by smooth data, and also approximate the non-linearity by a smooth nonlinearity, to obtain a smooth approximate solution, obtain an energy conservation law for that solution, and then take limits, using the well-posedness and perturbation theory in Section \ref{local-sec}.  We omit the standard details.  Similarly for the mass conservation law and Morawetz type inequalities.}
 by the flow corresponding to \eqref{equation}, we shall refer to it as
the \emph{energy} and often write $E(u)$ for $E(u(t))$.

A second conserved quantity we will rely on is the mass $M(u(t)):=\|u(t)\|^2_{L^2_x(\R^n)}$. As the mass is conserved,
we will often write $M(u)$ for $M(u(t))$.

In this paper, we will systematically study the initial value problem \eqref{equation}.  We are interested in
local and global well-posedness, asymptotic behaviour (scattering), and finite time blowup.  More precisely, we will
prove that under certain assumptions on the parameters $\lambda_1, \lambda_2, p_1, p_2$, we have the phenomena
mentioned above.

One of the motivations for considering this problem is the failure of the equation to be scale invariant.
For $p>0$, there is a natural scaling associated to the nonlinear Schr\"odinger equation
\begin{align}\label{one nonlin eq}
iu_t +\Delta u= |u|^p  u,
\end{align}
which leaves the equation invariant.  More precisely, the map
\begin{align}\label{scaling}
u(t,x)\mapsto \lambda^{-\frac{2}{p}} u\Bigl(\frac{t}{\lambda^2},\frac{x}{\lambda}\Bigr)
\end{align}
maps a solution to \eqref{one nonlin eq} to another solution to \eqref{one nonlin eq}.  In the case when
$p=\frac{4}{n}$, the scaling \eqref{scaling} also leaves the mass invariant, which is why the nonlinearity
$|u|^{\frac{4}{n}}u$ is called $L_x^2$-critical.  When $p=\frac{4}{n-2}$, the scaling \eqref{scaling} leaves
the energy invariant and hence, the nonlinearity $|u|^{\nmt}u$ is called $\dot H^1_x$- or energy-critical.
As our combined nonlinearity obeys $p_1<p_2$, there is no scaling that leaves \eqref{equation} invariant.
On the other hand, one can use scaling and homogeneity to normalize both $\lambda_1$ and $\lambda_2$ to have magnitude one
without difficulty.

The classical techniques used to prove local and global well-posedness in $H^1_x$ for \eqref{one nonlin eq}
(i.e., Picard's fixed point theorem combined with a standard iterative argument) do not distinguish between
the various values of $p$ as long as the nonlinearity $|u|^p u$ is energy-subcritical, that is, $0<p<\frac{4}{n-2}$;
for details see, for example, \cite{kato, cazenave:book}.
However, the proof of global well-posedness for the energy-critical nonlinear Schr\"odinger equation,
\begin{equation}\label{critical equation}
\begin{cases}
(i\partial_t+\Delta )w= |w|^{\frac{4}{n-2}} w\\
w(0,x) = w_{0}(x)\in \dot H^1_x
\end{cases}
\end{equation}
relies heavily on the scale invariance for this equation; see \cite{ckstt:gwp, RV, Monica:thesis}.
Hence, adding an energy-subcritical perturbation to \eqref{critical equation}, which destroys the scale invariance,
is of particular interest. This particular problem was first pursued by the third author, \cite{xzhang}, who considered
the case $n=3$.  The perturbative approach used in \cite{xzhang} extends easily to dimensions $n=4,5,6$. However,
in higher dimensions ($n>6$) new difficulties arise, mainly related to the low power of the energy-critical
nonlinearity. For instance, the continuous dependence of solutions to \eqref{critical equation} upon the initial data
in energy-critical spaces is no longer Lipschitz. Until recently, it was not even known whether one has \emph{uniform}
continuity of the solution upon the initial data in energy-critical spaces. This issue was settled by the first two authors,
\cite{TV}, who established a local well-posedness and stability theory which is H\"older continuous
(of order $\tfrac{4}{n-2}$) in energy-critical spaces and that applies even for large initial data, provided
a certain spacetime norm is known to be bounded. Basing our analysis on the stability theory developed in \cite{TV},
specifically Theorem 1.4, we will treat all dimensions $n\geq 3$ in a unified manner.

The local theory for \eqref{equation} is considered in Section~3.  Standard techniques involving Picard's fixed point
theorem can be used to construct local-in-time solutions to \eqref{equation}; in the case when an
energy-critical nonlinearity is present, that is, $p_2=\frac{4}{n-2}$, the time of existence for these local solutions
depends on the profile of the data, rather than on its $H^1_x$-norm.  After reviewing these classical statements,
we will develop a stability theory for the $L_x^2$-critical nonlinear Schr\"odinger equation and record the stability
result for the energy-critical NLS obtained by the first two authors, \cite{TV}.

Our first main result addresses the question of global well-posedness for \eqref{equation} in the energy space:

\begin{theorem}[Global well-posedness]\label{global wellposedness}
Let $u_0\in H_x^1$. Then, there exists a unique global solution to \eqref{equation} in each of the following cases:
\begin{enumerate}
\item $0<p_1<p_2<\frac 4n$ and $\lambda_1,\lambda_2\in\R$;
\item $0<p_1<p_2\leq \frac 4{n-2}$ and $\lambda_1\in\R$, $\lambda_2>0$.
\end{enumerate}
Moreover, for all compact intervals $I$, the global solution satisfies the following spacetime bound\footnote{In this paper we use $C$ to denote various large finite constants, which depend on the dimension $n$, the exponents $p_1,p_2$, the coefficients $\lambda_1, \lambda_2$, and any other quantities indicated by the parentheses (in this case, $|I|$ and $\|u_0\|_{H^1_x}$).  The exact value of $C$ will vary from line to line.}:
\begin{align}\label{bound wellposedness}
\|u\|_{S^1(I\times\R^n)}\le C\bigl(|I|, \|u_0\|_{H^1_x}\bigr).
\end{align}
\end{theorem}

We prove this theorem in Section~4.
The global existence of solutions to \eqref{equation} under the hypotheses of Theorem~\ref{global wellposedness}
is obtained as a consequence of three factors: the conservation of mass, an \emph{a priori} estimate on the kinetic
energy, and a `good' local well-posedness statement, by which we mean that the time of existence for local solutions
to \eqref{equation} in the two cases described in Theorem~\ref{global wellposedness} depends only on the
$H^1_x$-norm of the initial data.  This `good' local well-posedness statement coincides with the standard local
well-posedness statement when $0<p_1<p_2<\nmt$.  However, when $p_2=\nmt$ further analysis is needed as the standard
local well-posedness statement asserts that the time of existence for local solutions depends instead on the profile
of the initial data.  In order to upgrade the standard statement to the `good' statement we will make use of the
stability result in \cite{TV}.

In Section~5, we consider the asymptotic behaviour of these global solutions.  We will be able to obtain unconditional results in
the regime $\frac{4}{n} < p_1 < p_2 \leq \frac{4}{n-2}$.  It is natural to also seek the endpoint $p_1 = \frac{4}{n}$ for these results, but there is a difficulty because the defocusing $L_x^2$-critical NLS,
\begin{equation}\label{l2critical}
\begin{cases}
iv_t+\Delta v = |v|^{\frac{4}{n}} v\\
v(0,x) = v_{0}(x)
\end{cases}
\end{equation}
is not currently known to have a good scattering theory (except when the mass is small).
However, we will be able to obtain conditional results  in the $p_1 = \frac{4}{n}$ case assuming that a good theory for
\eqref{l2critical} exists.  More precisely, we will need the following

\begin{assumption}\label{assumption}
Let $v_0\in H^1_x$.  Then, there exists a unique global solution $v$ to \eqref{l2critical} and moreover,
$$
\|v\|_{L^{\frac{2(n+2)}n}_{t,x}(\R\times\R^n)}\le C(\|v_0\|_{L_x^2}).
$$
\end{assumption}

We can now state our second main result.

\begin{theorem}[Energy space scattering]\label{scattering}
Let $u_0\in H_x^1$, $\frac 4n\le p_1<p_2\le\frac 4{n-2}$, and let $u$ be the unique solution to \eqref{equation}.
If $p_1 = \frac{4}{n}$, then we also assume Assumption \ref{assumption}.
Then, there exists unique $u_{\pm}\in H^1_x$ such that
$$
\|u(t)-e^{it\Delta}u_{\pm}\|_{H^1_x}\to 0 \ \  \mbox{ as } \ \  t\to \pm \infty
$$
in each of the following two cases:
\begin{enumerate}
\item $\lambda_1$, $\lambda_2>0$;
\item $\lambda_1<0$, $\lambda_2>0$, and we have the small mass condition $M \leq c(\|\nabla u_0\|_2)$ for some suitably small
quantity $c(\|\nabla u_0\|_2) > 0$ depending only on $\|\nabla u_0\|_2$.
\end{enumerate}
\end{theorem}

\begin{remark} Note that in each of the two cases described in Theorem~\ref{scattering}, the unique solution
to \eqref{equation} is global by Theorem~\ref{global wellposedness}.
\end{remark}

We prove this theorem in Section ~5.
The scattering result for case (1) of the theorem is obtained in three
stages:

First, we develop an \emph{a priori} interaction Morawetz estimate; see subsection~5.1. This estimate is
particularly useful when both nonlinearities are defocusing (that is, both $\lambda_1$ and $\lambda_2$
are positive) and, as such, is an expression of dispersion (quantifying how mass is interacting with itself).
As a consequence of the interaction Morawetz inequality, we obtain $|\nabla|^{-\frac{n-3}{4}}u \in L_{t,x}^4$.
Interpolating between this estimate and the estimate on the kinetic energy (which is obviously bounded
when both nonlinearities are defocusing by the conservation of energy), we obtain control over the solution in the
$L_t^{n+1}L_x^{\frac{2(n+1)}{n-1}}$-norm.

The second step is to upgrade this bound to a global Strichartz bound using the stability results for the
$L_x^2$-critical and the energy-critical NLS; see subsection~5.2 through 5.5.  When both nonlinearities are defocusing
and $\frac{4}{n}<p_1<p_2=\nmt$, we view \eqref{equation} as a perturbation to the energy-critical NLS
\eqref{critical equation}, which is globally wellposed (see \cite{ckstt:gwp, RV, Monica:thesis}) and moreover,
the global solution satisfies
$$
\|w\|_{L_{t,x}^{\frac{2(n+2)}{n-2}}(\R\times\R^n)}\leq C(\|w_0\|_{\dot H^1_x}).
$$
Whenever the two nonlinearities are defocusing and $\frac{4}{n}=p_1<p_2<\nmt$, we view \eqref{equation}
as a perturbation to the pure power equation \eqref{l2critical} (normalizing $\lambda_1=1$).

Of particular interest (and difficulty) is the case when the nonlinearities are defocusing and $p_1=\frac{4}{n}$, $p_2=\nmt$.
In this case, the low frequencies of the solution are well approximated by the $L_x^2$-critical problem, while
the high frequencies are well approximated by the energy-critical problem.  The medium frequencies will eventualy be controlled by a Morawetz estimate.  Thus, in this case, the global Strichartz
bounds we derive are again conditional upon a satisfactory theory for the $L_x^2$-critical
NLS, that is, we need Assumption~\ref{assumption}.

In the intermediate case $\frac 4n<p_1<p_2<\frac 4{n-2}$, we reprove the classical scattering (in $H_x^1$) result
for solutions to \eqref{one nonlin eq} with $\frac 4n<p<\nmt$ due to J.~Ginibre and G.~Velo, \cite{gv:scatter}.
The proof we present in subsection~5.3 (Ode to Morawetz) is a new and simpler one based on the interaction Morawetz
estimate.

The last step required to obtain scattering under the assumptions described in Case 1) of Theorem~\ref{scattering}
is to prove that finite global Strichartz bounds imply scattering; see subsection~5.8.

In order to obtain finite global Strichartz norms that imply scattering in the second case described in
Theorem~\ref{scattering}, we make use of the small mass assumption (as a substitute for the interaction Morawetz
estimate) and of the stability result for the energy-critical NLS.  See subsections 5.6 and 5.7.

In the remaining two sections, we present our results on finite time blowup and global well-posedness and scattering
for \eqref{equation} with initial data $u_0\in\Sigma$, where $\Sigma$ is the space of all functions on $\R^n$ whose norm
$$ \| f\|_{\Sigma} := \|f\|_{H^1_x} + \| xf\|_{L^2_x}
$$
is finite (as usual, we identify functions which agree almost everywhere).

The finite time blowup result is a consequence of an argument of Glassey's, \cite{glassey}; in Section~6 we prove the
following

\begin{theorem}[Blowup]\label{blowup}
Let $u_0\in\Sigma$, $\lambda_2<0$, and $\frac 4n<p_2\leq\frac 4{n-2}$.  Let $y_0:=\Im \int_{\R^n}r\bar u_0\partial_r u_0 dx$
denote the weighted mass current and assume $y_0>0$. Then blowup occurs in each of the following three cases:\\
1) $\lambda_1>0$, $0<p_1<p_2$, and $E(u_0)<0$;\\
2) $\lambda_1<0$, $\frac 4n<p_1<p_2$, and $E(u_0)<0$;\\
3) $\lambda_1<0$, $0<p_1\le \frac 4n$, and $E(u_0)+CM(u_0)<0$ for some suitably large constant $C$ (depending as usual on
$n$, $p_1$, $p_2$, $\lambda_1$, $\lambda_2$).
\\
More precisely, in any of the above cases there exists $0<T_*\le C\frac{\|xu_0\|_2^2}{y_0}$ such that
$$
\lim_{t\to T_*}\|\nabla u(t)\|_{L_x^2}=\infty.
$$
\end{theorem}

\begin{remark}  When comparing the conditions in Case 1) and Case 3) of Theorem~\ref{blowup}, one might be puzzled
at first by the fact that we need stronger assumptions to prove blowup in the case of a focusing nonlinearity than in the case of a
defocusing nonlinearity.  However, one should notice that the condition
$$
E(t)=\int_{\R^n} \Bigl[\tfrac 12|\nabla u(t,x)|^2+\tfrac {\lambda_1}{p_1+2}|u(t,x)|^{p_1+2}+\tfrac{\lambda_2}{p_2+2}|u(t,x)|^{p_2+2}\Bigr]dx<0
$$
is easier to satisfy when $\lambda_1<0$ and $\lambda_2<0$.  Specifically, even when the kinetic energy of $u$
is small, which, in particular, implies that there is no blowup for $\|\nabla u(t)\|_2$, the energy $E$ can still be
made negative just by requiring that the mass be sufficiently large.  Hence, in order to obtain blowup of the kinetic
energy in this case, it is necessary to add a size restriction on the mass of the initial data.
\end{remark}

\begin{remark} In Theorem \ref{blowup}, we do not treat the endpoint $p_2=\frac 4n$.
For the focusing $L_x^2$-critical nonlinear Schr\"odinger equation, it is known that the blowup criterion is intimately
related to the properties of the unique spherically-symmetric, positive ground state of the elliptic equation
$$
-\Delta Q + \lambda_2 |Q|^{\frac 4n}Q = -Q.
$$
For results on this problem and a more detailed list of references see \cite{merle-raphael1, merle-raphael2}.
\end{remark}

In Section~7 we prove scattering in $\Sigma$ for solutions to \eqref{equation} with defocusing nonlinearities and
initial data $u_0\in\Sigma$.  More precisely, we obtain the following

\begin{theorem}[Pseudoconformal space scattering]\label{sigmascattering}
Let $u_0\in\Sigma$, $\lambda_1$ and $\lambda_2$ be positive numbers, and $\alpha(n)<p_1<p_2\le \frac 4{n-2}$
where $\alpha(n)$ is the Strauss exponent
$\alpha(n):=\frac{2-n+\sqrt{n^2+12n+4}}{2n}$.  Let $u$ be the unique global solution to \eqref{equation}.
Then, there exist unique scattering states $u_{\pm}\in \Sigma$ such that
$$
\|e^{-it\Delta}u(t)-u_{\pm}\|_{\Sigma}\to 0 \ \ \mbox{ as } \ \ t\to \pm\infty.
$$
\end{theorem}

We summarize our results in Table~\ref{Table1}.

\begin{table}[h]
\begin{center}
\def\up{\vrule width 0mm height 2.8ex depth 1.4ex}

\begin{tabular}{|c|c|c|c|c|l|}
\hline $\lambda_2$   & $\lambda_1$         & $p_1,p_2$                               & GWP &Scattering&Provided\\
\hline $\lambda_2>0$ & $\lambda_1\in\R$    & $0<p_1<p_2\le\frac{4}{n-2}$             &\checkmark & ? & --- \up \\
\hline $\lambda_2>0$ & $\lambda_1>0$       & $\frac{4}{n}\le p_1<p_2\le\frac{4}{n-2}$         &\checkmark & in $H_x^1$    & --- \up\\
\hline $\lambda_2>0$ & $\lambda_1\in\R$    & $\frac 4n\le p_1<p_2\le\frac{4}{n-2}$   &\checkmark & in $H_x^1$ & $M(u_0)\ll 1$\up \\
\hline $\lambda_2>0$ & $\lambda_1>0$       & $\alpha(n)<p_1<p_2\le\frac{4}{n-2}$     &\checkmark & in $\Sigma$ & $u_0\in\Sigma$\up\\
\hline $\lambda_2<0$ & $\lambda_1\in\R$    & $0<p_1<p_2<\frac{4}{n}$                 &\checkmark & ? & --- \up\\
\hline $\lambda_2<0$ & $\lambda_1>0$       & $0<p_1<p_2$, $\frac4n<p_2\le\frac 4{n-2}$  &$\times$ &$\times$ &$y_0>0$, $E(u_0)<0$ \up\\
\hline $\lambda_2<0$ & $\lambda_1<0$       & $\frac{4}{n}<p_1<p_2\le\frac{4}{n-2}$   &$\times$ & $\times$    & $y_0>0$, $E(u_0)<0$\up \\
\hline $\lambda_2<0$ & $\lambda_1<0$       & $0<p_1\le\frac 4n<p_2\le\frac{4}{n-2}$  &$\times$ &$\times$ & $y_0>0$, $E(u_0)+CM(u_0)<0$\up \\
 \hline
\end{tabular}
\end{center}
\caption{Summary of Results.  In all cases the initial data is assumed to lie in $H^1_x$.  The positive scattering results when $p_1=\frac{4}{n}$ are conditional on Assumption \ref{assumption}.}\label{Table1}
\end{table}

%
%
%
%

\section{Preliminaries}
We will often use the notation $X \lesssim Y$ whenever there exists
some constant $C$ so that $X \leq CY$; as before, $C$ can depend on $n$, $p_1, p_2$, $\lambda_1$, $\lambda_2$. Similarly, we will use $X
\sim Y$ if $X \lesssim Y \lesssim X$.  We use $X \ll Y$ if $X \leq
cY$ for some small constant $c$. The derivative operator $\nabla$
refers to the space variable only. We will occasionally use
subscripts to denote spatial derivatives and will use the summation
convention over repeated indices.

We use $L_x^r(\R^n)$ to denote the Banach space of functions $f:\R^n\to \C$ whose norm
$$
\|f\|_r:=\Bigl(\int_{\R^n} |f(x)|^r dx\Bigr)^{\frac{1}{r}}
$$
is finite, with the usual modifications when $r=\infty$, and identifying functions which agree almost everywhere.
For any non-negative integer $k$, we denote by
$W^{k,r}(\R^n)$ the Sobolev space defined as the closure of test functions in the norm
$$
\|f\|_{W^{k,r}(\R^n)}:=\sum_{|\alpha|\leq k}\Bigl\|\frac{\partial ^\alpha}{\partial x^\alpha}f\Bigr\|_r.
$$

We use $\lqlr$ to denote the spacetime norm
$$
\|u\|_{\llr} :=\Bigl(\int_{\R}\Bigl(\int_{\R^n} |u(t,x)|^r dx \Bigr)^{q/r} dt\Bigr)^{1/q},
$$
with the usual modifications when $q$ or $r$ is infinity, or when the domain $\R \times \R^n$ is
replaced by some smaller spacetime region.  When $q=r$ we abbreviate $\lqlr$ by $L^q_{t,x}$.

We define the Fourier transform on $\R^n$ to be
$$
\hat f(\xi) := \int_{\R^n} e^{-2 \pi i x \cdot \xi} f(x) dx.
$$

We will make use of the fractional differentiation operators $|\nabla|^s$ defined by
$$
\widehat{|\nabla|^sf}(\xi) := |\xi|^s \hat f (\xi).
$$
These define the homogeneous Sobolev norms
$$
\|f\|_{\dot H^s_x} := \| |\nabla|^s f \|_{L^2_x}.
$$

Let $e^{it\Delta}$ be the free Schr\"odinger propagator.  In physical space this is given by the formula
$$
e^{it\Delta}f(x) = \frac{1}{(4 \pi i t)^{n/2}} \int_{\R^n} e^{i|x-y|^2/4t} f(y) dy
$$
for $t\neq 0$ (using a suitable branch cut to define $(4\pi it)^{n/2}$),
while in frequency space one can write this as
\begin{equation}\label{fourier rep}
\widehat{e^{it\Delta}f}(\xi) = e^{-4 \pi^2 i t |\xi|^2}\hat f(\xi).
\end{equation}
In particular, the propagator obeys the \emph{dispersive inequality}
\begin{equation}\label{dispersive ineq}
\|e^{it\Delta}f\|_{L^\infty_x} \lesssim
|t|^{-\frac{n}{2}}\|f\|_{L^1_x}
\end{equation}
for all times $t\neq 0$.

We also recall \emph{Duhamel's formula}
\begin{align}\label{duhamel}
u(t) = e^{i(t-t_0)\Delta}u(t_0) - i \int_{t_0}^t e^{i(t-s)\Delta}(iu_t + \Delta u)(s) ds.
\end{align}

We say that a pair of exponents $(q,r)$ is Schr\"odinger-\emph{admissible} if $\tfrac{2}{q} +
\tfrac{n}{r} = \frac{n}{2}$ and $2 \leq q,r \leq \infty$. If $I \times \R^n$ is a spacetime slab, we
define the $\szi$ \emph{Strichartz norm} by
\begin{equation}\label{s0}
\|u\|_{\szi} := \sup \|u \|_{\lli}
\end{equation}
where the $\sup$ is taken over all admissible pairs $(q,r)$.  We
define the $\dot S^1 (\ir)$ \emph{Strichartz norm} to be
$$
\|u\|_{\dot S^1 (\ir)} := \| \nabla u \|_{\szi}.
$$
We also use $\dot N^0(\ir)$ to denote the dual space of $\dot
S^0(\ir)$ and
$$
\dot N^1(\ir):=\{u;\ \nabla u\in \dot N^0(\ir)\}.
$$

By definition and Sobolev embedding, we obtain

\begin{lemma}\label{lemma strichartz norms}
For any $\dot{S^1}$ function $u$ on $\ir$, we have
\begin{align*}\label{strichartz norms}
&\|\nabla u\|_{L^\infty_t L^2_x} + \|\nabla u \|_{L^{\frac{2(n+2)}{n-2}}_t L^{\frac{2n(n+2)}{n^2+4}}_x}
+
\|\nabla u\|_{L^{\frac{2(n+2)}{n}}_{t,x}} + \|\nabla u\|_{L^2_t L^{\frac{2n}{n-2}}_x} \\
&\quad + \|u\|_{L^\infty_t L^{\frac{2n}{n-2}}_x} +
\|u\|_{L^{\frac{2(n+2)}{n-2}}_{t,x}} +\|u\|_{L^{\frac{2(n+2)}{n}}_t
L^{\frac{2n(n+2)}{n^2-2n-4}}_x} \lesssim \|u\|_{\so}
\end{align*}
where all spacetime norms are on $\ir$.
\end{lemma}

Let us also record the following standard Strichartz estimates that we will invoke throughout the paper
(for a proof see \cite{tao:keel}):

\begin{lemma}\label{lemma linear strichartz}
Let $I$ be a compact time interval, $k=0,1$, and let $u : \ir \rightarrow \C$ be an $\dot{S^k}$ solution
to the forced Schr\"odinger equation
\begin{equation*}
i u_t + \Delta u = F
\end{equation*}
for a function $F$.  Then we have
\begin{equation}
\|u\|_{\ski} \lesssim \|u(t_0)\|_{\hk (\R^n)} +  \|F\|_{\dot N^k(\ir)},
\end{equation}
 for any time $t_0 \in I$.
\end{lemma}

We will also need some Littlewood-Paley theory.  Specifically, let $\varphi(\xi)$ be a smooth bump supported in the ball
$|\xi| \leq 2$ and equalling one on the ball $|\xi| \leq 1$.  For each dyadic number $N \in 2^\Z$ we define the
Littlewood-Paley operators
\begin{align*}
\widehat{P_{\leq N}f}(\xi) &:=  \varphi(\xi/N)\hat f (\xi),\\
\widehat{P_{> N}f}(\xi) &:=  [1-\varphi(\xi/N)]\hat f (\xi),\\
\widehat{P_N f}(\xi) &:=  [\varphi(\xi/N) - \varphi (2 \xi /N)] \hat f (\xi).
\end{align*}
Similarly we can define $P_{<N}$, $P_{\geq N}$, and $P_{M < \cdot \leq N} := P_{\leq N} - P_{\leq M}$, whenever $M$ and
$N$ are dyadic numbers.  We will frequently write $f_{\leq N}$ for $P_{\leq N} f$ and similarly for the other operators.
We recall some standard Bernstein type inequalities:

\begin{lemma}\label{lemma fourier projection}
For any $1\le p\le q\le\infty$ and $s>0$, we have
\begin{align*}
\|P_{\geq N} f\|_{L^p_x} &\lesssim N^{-s} \| |\nabla|^s P_{\geq N} f \|_{L^p_x}\\
\| |\nabla|^s  P_{\leq N} f\|_{L^p_x} &\lesssim N^{s} \| P_{\leq N} f\|_{L^p_x}\\
\| |\nabla|^{\pm s} P_N f\|_{L^p_x} &\sim N^{\pm s} \| P_N f \|_{L^p_x}\\
\|P_{\leq N} f\|_{L^q_x} &\lesssim N^{\frac{n}{p}-\frac{n}{q}} \|P_{\leq N} f\|_{L^p_x}\\
\|P_N f\|_{L^q_x} &\lesssim N^{\frac{n}{p}-\frac{n}{q}} \| P_N f\|_{L^p_x}.
\end{align*}
\end{lemma}

For our analysis and the sake of the exposition, it is convenient to introduce a number of function spaces.
We will need the following Strichartz spaces defined on a slab $\ir$ as the closure of the test functions
under the appropriate norms:
\begin{align*}
\|u\|_{V(I)}&:= \|u\|_{L_{t,x}^{\frac{2(n+2)}{n}}(\ir)}\\
\|u\|_{W(I)}&:= \|u\|_{L_{t,x}^{\frac{2(n+2)}{n-2}}(\ir)}\\
\|u\|_{Z(I)}&:= \|u\|_{L^{n+1}_tL_x^{\frac{2(n+1)}{n-1}}(\ir)}.
\end{align*}

\begin{definition}\label{def xi spaces}
Let $\ir$ be an arbitrary spacetime slab.  For $0<p_1<p_2<\frac 4{n-2}$, we define the space
$$
\dot X^0(I):=\cap_{i=1,2}L_t^{\gamma_i}L_x^{\rho_i}(I\times\R^n),
$$
and
$$
\dot X^1(I):=\{u;\, \nabla u\in\dot X^0(I)\},\ \ \ X^1(I):=\dot X^0(I)\cap\dot X^1(I),
$$
where $\gamma_i:=\frac {4(p_i+2)}{p_i(n-2)}$, $\rho_i:=\frac{n(p_i+2)}{n+p_i}$.  It is not hard to check that
$(\gamma_i,\rho_i)$ is a Schr\"odinger admissible pair and thus $\dot S^0\subset \dot X^0$.

In the case $0<p_1<p_2=\nmt$, we define the spaces
$$
\dot X^0(I):=L_t^{\gamma_1}L_x^{\rho_1}(\ir)\cap L_t^{\frac{2(n+2)}{n-2}}L_x^{\frac{2n(n+2)}{n^2+4}}(\ir)\cap V(I)
$$
and
$$
\dot X^1(I):=\{u;\, \nabla u\in\dot X^0(I)\},\ \ \ X^1(I):=\dot X^0(I)\cap\dot X^1(I).
$$
\end{definition}

Define $\rho_i^*:=\frac {n(p_i+2)}{n-2}$ and let $\gamma_i'$, $\rho_i'$ be the dual the exponents of $\gamma_i$,
respectively $\rho_i$ introduced in Definition~\ref{def xi spaces}.  It is easy to verify that the following identities
hold:
\begin{align}
\frac 1{\gamma_i'}&=1-\frac {p_i(n-2)}4+\frac {p_i+1}{\gamma} \label{gamma}\\
\frac 1{\rho'_i}&=\frac 1{\rho_i}+\frac p{\rho_i^*}\label{rho}\\
\frac 1{\rho_i^*}&=\frac 1{\rho_i}-\frac 1n. \label{rho*}
\end{align}

Using \eqref{gamma} through \eqref{rho*} as well as H\"older and Sobolev embedding, we derive the first estimates
for our nonlinearity.

\begin{lemma}\label{control by x}
Let $I$ be a compact time interval, $0<p_1<p_2\leq \nmt$, $\lambda_1$ and $\lambda_2$ be nonzero real numbers,
and $k=0,1$.  Then,
$$
\bigl\|\ntu\bigr\|_{\dot N^k(I\times\R^n)}
\lesssim \sum_{i=1,2} |I|^{1-\frac{p_i(n-2)}4}\|u\|_{\dot X^1(I)}^{p_i}\|u\|_{\dot X^k(I)}
$$
and
\begin{align*}
\bigl\|\bigl(\ntu\bigr)-&\bigl(\ntv\bigr)\bigr\|_{\dot N^0(I\times\R^n)}\\
&\lesssim\sum_{i=1,2}|I|^{1-\frac{p_i(n-2)}4}\bigl(\|u\|_{\dot X^1(I)}^{p_i}+\|v\|_{\dot X^1(I)}^{p_i}\bigr)\|u-v\|_{\dot X^0(I)}.
\end{align*}
\end{lemma}

When the length of the time interval $I$ is infinite, the estimates in Lemma~\ref{control by x} are useless.
In this case we will use instead the following

\begin{lemma}\label{lemma bound by w01}
Let $\ir$ be an arbitrary spacetime slab, $\frac 4n\le p\le \nmt$, and $k=0,1$.  Then
\begin{equation}\label{bound by w01}
\||u|^pu\|_{\dot N^k(\ir)}
\lesssim \|u\|_{V(I)}^{2-\frac{p(n-2)}2}\|u\|_{W(I)}^{\frac{np}2-2}\||\nabla|^k u\|_{V(I)}.
\end{equation}
\end{lemma}
\begin{proof} By the boundedness of the Riesz potentials on any $L_x^p$, $1<p<\infty$, H\"older, and interpolation, we estimate
\begin{align*}
\||u|^pu\|_{\dot N^k(\ir)}
&\lesssim \bigl\||\nabla|^k(|u|^pu)\bigr\|_{L_{t,x}^{\frac{2(n+2)}{n+4}}(\ir)}\\
&\lesssim \||u|^p\|_{L_{t,x}^{\frac{n+2}2}(\ir)}\||\nabla|^k u\|_{L_{t,x}^{\frac{2(n+2)}n}(\ir)}\\
&\lesssim \|u\|_{L_{t,x}^{\frac{(n+2)p}2}(\ir)}^p\||\nabla|^k u\|_{L_{t,x}^{\frac{2(n+2)}n}(\ir)}\\
&\lesssim\|u\|_{L_{t,x}^{\frac{2(n+2)}n}(\ir)}^{2-\frac{p(n-2)}2}\|u\|_{L_{t,x}^{\frac{2(n+2)}{n-2}}(\ir)}^{\frac{np}2-2}\||\nabla|^k u\|_{L_{t,x}^{\frac{2(n+2)}n}(\ir)},
\end{align*}
which proves \eqref{bound by w01}.
\end{proof}

When deriving global spacetime bounds that imply scattering, we would like to involve the $Z$-norm which corresponds
to the control given by the \emph{a priori} interaction Morawetz estimate.
For $k=0,1$ and $\frac 4n<p<\frac 4{n-2}$, applying H\"older's inequality we estimate
\begin{align*}
\bigl\||\nabla|^k(|u|^pu)&\bigr\|_{L_t^2L_x^{\frac{2n}{n+2}}(\ir)}\\
&\lesssim \||\nabla|^k u\|_{L_t^2L_x^{\frac{2n}{n-2}}(\ir)} \|u\|^p_{L_t^\infty L_x^{\frac{np}{2}}(\ir)}\\
&\lesssim \||\nabla|^k u\|_{L_t^2L_x^{\frac{2n}{n-2}}(\ir)}\|u\|^{2-\frac{p(n-2)}{2}}_{L_t^\infty L_x^2(\ir)}
   \|u\|^{\frac{np-4}{2}}_{L_t^\infty L_x^{\frac{2n}{n-2}}(\ir)}.
\end{align*}
In order to use the \emph{a priori} $L_t^{n+1}L_x^{\tfrac{2(n+1)}{n-1}}$ control (given to us by the interaction
Morawetz estimate in the case when both nonlinearities are defocusing), we would like to replace the
$L_t^\infty L_x^{\frac{np}{2}}$-norm by a norm which belongs to the open triangle determined by the mass
($L_t^\infty L_x^2$), the potential energy ($L_t^\infty L_x^{\frac{2n}{n-2}}$), and the
$L_t^{n+1}L_x^{\tfrac{2(n+1)}{n-1}}$-norm. This can be achieved by increasing the time exponent in the
$L_t^2W_x^{k,\frac{2n}{n-2}}$-norm by a tiny amount, while maintaining the scaling (by which we mean replacing the pair
$(2,\tfrac{2n}{n-2})$ by another Schr\"odinger-admissible pair). We obtain the following

\begin{lemma}\label{m interpolation lemma}
Let $k=0,1$ and $\frac 4n<p<\frac 4{n-2}$.  Then, there exists $\theta>0$ large enough such that on every slab $\ir$
we have
\begin{align}
\bigl\||\nabla|^k( |u|^pu)\bigr\|_{L^2_tL_x^{\frac{2n}{n+2}}}
&\lesssim \||\nabla|^k u\|_{L_t^{2+\frac 1\theta}L_x^{\frac{2n(2\theta+1)}{n(2\theta+1)-4\theta}}}
   \|u\|^{\frac{n+1}{2(2\theta+1)}}_{L^{n+1}_t L_x^{\frac{2(n+1)}{n-1}}}
   \|u\|_{L^{\infty}_tL_x^2}^{\alpha(\theta)} \|u\|^{\beta(\theta)}_{L^{\infty}_tL_x^{\frac{2n}{n-2}}}\notag\\
&\lesssim \|u\|_{\dot S^k(\ir)} \|u\|_{Z(I)}^{\frac{n+1}{2(2\theta+1)}}\|u\|_{L^{\infty}_tL_x^2}^{\alpha(\theta)}
   \|u\|^{\beta(\theta)}_{L^{\infty}_tL_x^{\frac{2n}{n-2}}},\label{minter}
\end{align}
where
$$
\alpha(\theta):=p\bigl(1-\tfrac n2\bigr)+\tfrac {8\theta+1}{2(2\theta+1)} \quad \text{and} \quad
\beta(\theta):=\tfrac n2\bigl(p-\tfrac{n+8\theta+2}{n(2\theta+1)}\bigr).
$$
\end{lemma}
\begin{proof}
First note that the pair $(2+\frac 1\theta,\ \frac
{2n(2\theta+1)}{n(2\theta+1)-4\theta})$ is Schr\"odinger-admissible.

Once $\alpha(\theta)$ and $\beta(\theta)$ are positive, \rm
(\ref{minter}) is a direct consequence of H\"older's inequality, as
the reader can easily check. It is not hard to see that
$\theta\mapsto\alpha(\theta)$ and $\theta\mapsto\beta(\theta)$ are
increasing functions and moreover,
$$
\alpha(\theta)\to p(1-\tfrac n2)+2 \quad \text{and} \quad
\beta(\theta)\to \tfrac n2(p-\tfrac 4n) \quad \text{as} \quad
\theta\to\infty.
$$
As $\frac 4n<p<\frac 4{n-2}$, the two limits are positive. Thus, for
$\theta$ sufficiently large we obtain
$$
\alpha(\theta)>0 \quad \text{and} \quad \beta(\theta)>0.
$$
This concludes the proof of the lemma.
\end{proof}

When $p=\nmt$, we can still control the $\dot N^0$-norm of $|u|^pu$ in terms of the $Z$-norm.
The idea is simple: Note that by H\"older, we have
\begin{align}\label{model}
\||u|^{\nmt}u\|_{\dot N^0(\ir)}
&\lesssim \||u|^{\nmt} u\|_{L^2_tL_x^{\frac{2n}{n+2}}(\ir)}\notag\\
&\lesssim \|u\|_{L^2_tL_x^{\frac{2n}{n-2}}(\ir)}\|u\|_{L_t^\infty L_x^{\frac{2n}{n-2}}(\ir)}^{\frac{4}{n-2}}.
\end{align}
In order to get a small fractional power of $\|u\|_{Z(I)}$ on the right-hand side of \eqref{model}, we need to perturb
the above estimate a little bit.  More precisely, we will replace the norm $L^2_tL_x^{\frac{2n}{n-2}}$ by
$L^{2+\eps}_tL_x^{\frac{2n}{n-2-\eps}}$ for a small constant $\eps>0$.  The latter norm interpolates between the
$\dot S^0$-norm $L^{2+\eps}_tL_x^{\frac{2n(2+\eps)}{n(2+\eps)-4}}$ and the $\dot S^1$-norm
$L^{2+\eps}_tL_x^{\frac{2n(2+\eps)}{n(2+\eps)-2(4+\eps)}}$, provided $\eps$ is sufficiently small.  Thus,
\begin{align}\label{variation in s1}
\|u\|_{L^{2+\eps}_tL_x^{\frac{2n}{n-2-\eps}}(\ir)}\lesssim \|u\|_{S^1(\ir)},
\end{align}
provided $\eps$ is chosen sufficiently small.  Keeping the $L^2_tL_x^{\frac{2n}{n+2}}$-norm on the left-hand side of
\eqref{model}, this change forces us to replace the norm $L_t^{\infty}L_x^{\frac{2n}{n-2}}$ (which appears
on the right-hand side of \eqref{model}) with a norm which lies in the open triangle determined by the potential energy
($L_t^{\infty}L_x^{\frac{2n}{n-2}}$), the mass ($L_t^\infty L_x^2$), and the $Z$-norm.  Therefore, we have the following

\begin{lemma}\label{h1 morawetz control}
Let $\ir$ be a spacetime slab.  Then, there exists a small constant $0<\theta<1$ such that
\begin{align}\label{est in Z}
\||u|^{\nmt}u\|_{\dot N^0(\ir)}
\lesssim \|u\|_{Z(I)}^{\theta}\|u\|_{S^1(I\times\R^n)}^{\frac{n+2}{n-2}-\theta}.
\end{align}
\end{lemma}

\begin{proof}
We will in fact prove the following estimate
\begin{align}\label{est for crit in Z}
\||u|^{\nmt}u\|_{L_t^2L_x^{\frac{2n}{n+2}}}
\lesssim \|u\|_{L^{2+\eps}_tL_x^{\frac{2n}{n-2-\eps}}} \|u\|_{Z(I)}^{\frac{(n+1)\eps}{2(2+\eps)}}
  \|u\|_{L_t^\infty L_x^2}^{a(\eps)} \|u\|_{L_t^\infty L_x^{\frac{2n}{n-2}}}^{b(\eps)},
\end{align}
which holds for $\eps>0$ sufficiently small.  Here, all spacetime norms are on the slab $\ir$ and
$$
a(\eps):=\frac{(1+\eps)\eps}{2(2+\eps)} \ \ \text{and} \ \ b(\eps):=\frac{4}{n-2}-\frac{(n+2+\eps)\eps}{2(2+\eps)}.
$$

In order to prove \eqref{est for crit in Z}, we just need to check that for $\eps>0$ sufficiently small, $a(\eps)$ and
$b(\eps)$ are positive, since then the estimate is a simple consequence of H\"older's inequality.

It is easy to see that as functions of $\eps$, $a$ is increasing and $a(0)=0$, while $b$ is decreasing and $b(0)=\nmt$.
Thus, taking $\eps>0$ sufficiently small, we have $a(\eps)>0$ and $b(\eps)>0$, which yields \eqref{est for crit in Z}
for the reasons discussed above.

Taking $\theta:=\frac{(n+1)\eps}{2(2+\eps)}$ and using \eqref{variation in s1}, we obtain \eqref{est in Z}.

\end{proof}

\begin{remark}\label{variation}
An easy consequence of \eqref{est for crit in Z} are estimates for nonlinearities of the form $|u|^{\nmt}v$.
More precisely, on every spacetime slab $\ir$ we have
$$
\||u|^{\nmt}v\|_{\dot N^0(\ir)}
\lesssim \|u\|_{Z(I)}^{\theta}\|u\|_{S^1(\ir)}^{\nmt-\theta}\|v\|_{S^1(\ir)}
$$
and
$$
\||u|^{\nmt}v\|_{\dot N^0(\ir)}
\lesssim \|v\|_{Z(I)}^{\theta}\|u\|_{S^1(\ir)}^{\nmt}\|v\|_{S^1(\ir)}^{1-\theta}.
$$
\end{remark}

%
%
%
%

\section{Local theory}\label{local-sec}

In this section we present the local theory for the initial value problem \eqref{equation}.
We start by recording the local well-posedness statements.  As the material is classical, we prefer to omit
the proofs and instead send the reader to the detailed expositions in \cite{cwI, cazenave:book, kato, katounique}.

\begin{proposition}[Local well-posedness for \eqref{equation} with $\dot H^1_x$-subcritical nonlinearities]\label{subcritical lwp}\leavevmode\\
Let $u_0\in H_x^1$, $\lambda_1$ and $\lambda_2$ be nonzero real constants, and $0<p_1<p_2<\nmt$.
Then, there exists $T=T(\|u_0\|_{H^1_x})$ such that \eqref{equation} with the parameters given above admits
a unique strong $H_x^1$-solution $u$ on $[-T,T]$.  Let $(-T_{min}, T_{max})$ be the maximal time interval
on which the solution $u$ is well-defined.  Then, $u\in S^1 (\ir)$ for every compact time interval
$I\subset (-T_{min}, T_{max})$ and the following properties hold:\\
$\bullet$ If $T_{max}<\infty$, then
$$
\lim_{t\to T_{max}}\|u(t)\|_{\dot H_x^1}=\infty;
$$
similarly, if $T_{min}<\infty$, then
$$
\lim_{t\to -T_{min}}\|u(t)\|_{\dot H_x^1}=\infty.
$$
$\bullet$ The solution $u$ depends continuously on the initial data $u_0$ in the following sense:
There exists $T=T(\|u_0\|_{H^1_x})$ such that if $u_0^{(m)}\to u_0$ in $H^1_x$ and if $u^{(m)}$ is the
solution to \eqref{equation} with initial data $u_0^{(m)}$, then $u^{(m)}$ is defined on $[-T,T]$ for $m$
sufficiently large and $u^{(m)}\to u$ in $S^1([-T,T]\times\R^n)$.
\end{proposition}

\begin{proposition}[Local well-posedness for \eqref{equation} with a $\dot H^1_x$-critical nonlinearity]\label{critical lwp}\leavevmode\\
Let $u_0\in H_x^1$, $\lambda_1$ and $\lambda_2$ be nonzero real constants, and $0<p_1<p_2=\nmt$.
Then, for every $T>0$, there exists $\eta=\eta(T)$ such that if
$$
\|e^{it\Delta} u_0\|_{X^1([-T,T])}\leq \eta,
$$
then \eqref{equation} with the parameters given above admits a unique strong $H_x^1$-solution $u$ defined on
$[-T,T]$.  Let $(-T_{min}, T_{max})$ be the maximal time interval on which $u$ is well-defined.  Then,
$u\in S^1 (\ir)$ for every compact time interval $I\subset (-T_{min}, T_{max})$ and the following properties hold:\\
$\bullet$ If $T_{max}<\infty$, then
$$
\text{either} \ \ \lim_{t\to T_{max}}\|u(t)\|_{\dot H_x^1}=\infty \ \ \text{or}\ \
\|u\|_{\dot S^1((0,T_{max})\times\R^n)}=\infty.
$$
Similarly, if $T_{min}<\infty$, then
$$
\text{either}\ \ \lim_{t\to -T_{min}}\|u(t)\|_{\dot H_x^1}=\infty \ \ \text{or}\ \
\|u\|_{\dot S^1((-T_{min},0)\times\R^n)}=\infty.
$$
$\bullet$ The solution $u$ depends continuously on the initial data $u_0$ in the following sense:
The functions $T_{min}$ and $T_{max}$ are lower semicontinuous from $H^1_x$ to $(0, \infty]$.  Moreover,
if $u_0^{(m)}\to u_0$ in $H^1_x$ and if $u^{(m)}$ is the maximal solution to \eqref{equation} with initial data
$u_0^{(m)}$, then $u^{(m)}\to u$ in $L^q_t H_x^1([-S,T]\times\R^n)$ for every $q<\infty$ and every interval
$[-S,T]\subset(-T_{min}, T_{max})$.
\end{proposition}

We record also the following companion to Proposition~\ref{critical lwp}.

\begin{lemma}[Blowup criterion]\label{x1 implies existence}
Let $u_0\in H^1_x$ and let $u$ be the unique strong $H_x^1$-solution to \eqref{equation} with $p_2=\nmt$ on the slab
$[0,T_0]\times\R^n$ such that
$$
\|u\|_{\dot X^1([0,T_0])}<\infty.
$$
Then, there exists $\delta=\delta(u_0)>0$ such that the solution $u$ extends to a strong $H^1_x$-solution
on the slab $[0, T_0+\delta]\times\R^n$.
\end{lemma}

The proof is standard (see, for example, \cite{cazenave:book}).  In the contrapositive, this lemma asserts that
if a solution cannot be continued strongly beyond a time $T_0$, then the $\dot X^1$-norm (and all other $\dot S^1$-norms) must blow up at that time.

Next, we will establish a stability result for the $L_x^2$-critical NLS, by which we mean the following property:
Given an \emph{approximate} solution
\begin{equation*}
\begin{cases}
i \tilde v_t +\Delta \tilde v = |\tilde v|^{\frac{4}{n}}\tilde v + e\\
\tilde v(0,x)= \tilde v_0(x) \in L^2(\R^n)
\end{cases}
\end{equation*}
to \eqref{l2critical}, with $e$ small in a suitable space and $\tilde v_0 - v_0$ small in $L^2_x$,
there exists a \emph{genuine} solution $v$ to \eqref{l2critical} which stays very close to $\tilde v$
in $L_x^2$-critical norms.

\begin{lemma}[Short-time perturbations]\label{lemma l2 perturbation-0}
Let $I$ be a compact interval and let $\tilde v$ be an approximate solution to \eqref{l2critical} in the sense that
$$
(i\partial_t+\Delta)\tilde v=|\tilde v|^{\frac 4n}\tilde v+e,
$$
for some function $e$. Assume that
\begin{align}\label{small mass-0}
\|\tilde v\|_{L_t^{\infty}L_x^2(I\times\R^n)}&\le M
\end{align}
for some positive constant $M$. Let $t_0\in I$ and let $v(t_0)$ close to $\tilde v(t_0)$ in the sense that
\begin{align}\label{close l2-0}
\|v(t_0)-\tilde v(t_0)\|_{L_x^2}\le M'
\end{align}
for some $M'>0$.  Assume also the smallness conditions
\begin{align}
\|\tilde v\|_{V(I)}&\le \eps_0 \label{small norm-0}\\
\bigr\|e^{i(t-t_0)\Delta}\bigr(v(t_0)-\tilde v(t_0)\bigl)\bigl\|_{V(I)} &\leq \eps \label{closer l2-0}\\
\|e\|_{\dot N^0(I\times\R^n)}&\le \eps,\label{error small-0}
\end{align}
for some $0<\eps\le \eps_0$ where $\eps_0=\eps_0(M,M')>0$ is a small constant.  Then, there exists a solution $v$
to \eqref{l2critical} on $\ir$ with initial data $v(t_0)$ at time $t=t_0$ satisfying
\begin{align}
\|v- \tilde v\|_{V(I)} &\lesssim \eps \label{concl1-0}\\
\|v- \tilde v\|_{\dot S^0(I\times\R^n)} &\lesssim M' \label{concl2-0}\\
\|v\|_{\dot S^0(I\times\R^n)} &\lesssim M+ M'\label{concl3-0}\\
\|(i\partial_t+\Delta)(v-\tilde v)+e\|_{\dot N^0(\ir)} &\lesssim \eps. \label{concl4-0}
\end{align}
\end{lemma}

\begin{remark}\label{strich1} Note that by Strichartz,
$$
\bigr\|e^{i(t-t_0)\Delta}\bigr(v(t_0)-\tilde v(t_0)\bigl)\bigl\|_{V(I)}
\lesssim \|v(t_0)-\tilde v(t_0)\|_{L_x^2},
$$
so the hypothesis \eqref{closer l2-0} is redundant if $M'=O(\eps)$.
\end{remark}

\begin{proof}
By time symmetry, we may assume $t_0=\inf I$.  Let $z:=v-\tilde v$.  Then $z$ satisfies the following initial
value problem
\begin{equation*}
\begin{cases}
iz_t+\Delta z=|\tilde v+z|^{\frac 4n}(\tilde v+z)-|\tilde v|^{\frac 4n}\tilde v -e\\
z(t_0)=v(t_0)-\tilde v(t_0).
\end{cases}
\end{equation*}

For $t\in I$ define
$$
S(t):=\|(i\partial_t+\Delta)z+e\|_{\dot N^0([t_0,t]\times\R^n)}.
$$
By \eqref{small norm-0}, we have
\begin{align}\label{S(t)}
S(t)
&\lesssim \|(i\partial_t+\Delta)z+e\|_{L_{t,x}^{\frac{2(n+2)}{n+4}}([t_0,t]\times\R^n)}\notag\\
&\lesssim \|z\|_{V([t_0,t])}^{1+\frac4n}+\|z\|_{V([t_0,t])}\|\tilde v\|_{V([t_0,t])}^{\frac 4n}\notag\\
&\lesssim \|z\|_{V([t_0,t])}^{1+\frac 4n}+\eps_0^{\frac{4}{n}}\|z\|_{V([t_0,t])}.
\end{align}

On the other hand, by Strichartz, \eqref{closer l2-0}, and \eqref{error small-0}, we get
\begin{align}\label{z}
\|z\|_{V([t_0,t])}
&\lesssim\|e^{i(t-t_0)\Delta}z(t_0)\|_{V([t_0,t])}+S(t) +\|e\|_{\dot N^0([t_0,t]\times\R^n)}\notag\\
&\lesssim S(t)+\eps.
\end{align}

Combining \eqref{S(t)} and \eqref{z}, we obtain
$$
S(t)\lesssim (S(t)+\eps)^{1+\frac 4n}+\eps_0^{\frac4n}(S(t)+\eps)+\eps.
$$
A standard continuity argument then shows that if $\eps_0=\eps_0(M,M')$ is taken sufficiently small, then
$$
S(t)\le \eps \ \ \text{for any}\ \ t\in I,
$$
which implies \eqref{concl4-0}.  Using \eqref{concl4-0} and \eqref{z}, one easily derives \eqref{concl1-0}.
Moreover, by Strichartz, \eqref{close l2-0}, \eqref{error small-0}, and \eqref{concl4-0},
\begin{align*}
\|z\|_{\dot S^0(\ir)}
&\lesssim \|z(t_0)\|_{L_x^2}+\|(i\partial_t+\Delta )z+e\|_{\dot N^0(I\times\R^n)}+\|e\|_{\dot N^0([t_0,t]\times\R^n)}\\
&\lesssim M'+\eps,
\end{align*}
which establishes \eqref{concl2-0}.

To prove \eqref{concl3-0}, we use Strichartz, \eqref{small mass-0}, \eqref{close l2-0}, \eqref{small norm-0},
\eqref{error small-0}, and \eqref{concl4-0}:
\begin{align*}
\|v\|_{\dot S^0(\ir)}
&\lesssim \|v(t_0)\|_{L_x^2}+\|(i\partial_t+\Delta )v\|_{\dot N^0(I\times\R^n)}\\
&\lesssim \|\tilde v(t_0)\|_{L_x^2}+ \|v(t_0)-\tilde v(t_0)\|_{L_x^2}+ \|(i\partial_t+\Delta )(v-\tilde v)+e\|_{\dot N^0(I\times\R^n)}\\
&\quad + \|(i\partial_t+\Delta )\tilde v\|_{\dot N^0(I\times\R^n)}+\|e\|_{\dot N^0(\ir)}\\
&\lesssim M + M'+ \eps +\|(i\partial_t+\Delta )\tilde v\|_{L_{t,x}^{\frac{2(n+2)}{n+4}}(I\times\R^n)}\\
&\lesssim M + M'+ \|\tilde v\|^{1+\frac{4}{n}}_{V(I)}\\
&\lesssim M + M'+ \eps_0^{1+\frac{4}{n}}\\
&\lesssim M + M'.
\end{align*}
\end{proof}

Building upon the previous lemma, we have the following

\begin{lemma}[$L_x^2$-critical stability result]\label{l2 stability}
Let $I$ be a compact interval and let $\tilde v$ be an approximate solution to \eqref{l2critical} in the sense that
$$
(i\partial_t+\Delta)\tilde v=|\tilde v|^{\frac 4n}\tilde v+e,
$$
for some function $e$. Assume that
\begin{align}
\|\tilde v\|_{L_t^{\infty}L_x^2(I\times\R^n)}&\le M\label{small mass-1}\\
\|\tilde v\|_{V(I)}&\leq L \label{finite norm-1}
\end{align}
for some positive constants $M$ and $L$. Let $t_0\in I$ and let $v(t_0)$ close to $\tilde v(t_0)$ in the sense that
\begin{align}\label{close-1}
\|v(t_0)-\tilde v(t_0)\|_{L_x^2}\le M'
\end{align}
for some $M'>0$.  Moreover, assume the smallness conditions
\begin{align}
\bigr\|e^{i(t-t_0)\Delta}\bigr(v(t_0)-\tilde v(t_0)\bigl)\bigl\|_{V(I)} &\leq \eps \label{closer l2-1}\\
\|e\|_{\dot N^0(I\times\R^n)}&\le \eps,\label{error small-1}
\end{align}
for some $0<\eps\le \eps_1$ where $\eps_1=\eps_1(M,M',L)>0$ is a small constant.  Then, there exists a solution $v$
to \eqref{l2critical} on $\ir$ with initial data $v(t_0)$ at time $t=t_0$ satisfying
\begin{align}
\|v-\tilde v\|_{V(I)} &\le \eps C(M,M',L)\label{concl1-1}\\
\|v-\tilde v\|_{\dot S^0(I\times\R^n)} &\le  C(M,M',L)M' \label{concl2-1}\\
\|v\|_{\dot S^0(I\times\R^n)} &\le  C(M,M',L). \label{concl3-1}
\end{align}
\end{lemma}

\begin{remark} By Strichartz, the hypothesis \eqref{closer l2-1} is redundant if $M'=O(\eps)$; see Remark~\ref{strich1}.
Assumption \ref{assumption} is not explicitly used in the proof of this lemma, although in practice one needs an assumption like this if one wishes to obtain the hypothesis \eqref{finite norm-1}.
\end{remark}

\begin{proof}
Subdivide $I$ into $N\sim (1+\frac{L}{\eps_0})^{\frac{2(n+2)}{n}}$ subintervals $I_j=[t_j, t_{j+1}]$ such that
$$
\|\tilde v\|_{V(I_j)}\sim \eps_0,
$$
where $\eps_0=\eps_0(M,2M')$ is as in Lemma~\ref{lemma l2 perturbation-0}.  We need to replace $M'$ by $2M'$ as
the mass of the difference $v-\tilde v$ might grow slightly in time.

Choosing $\eps_1$ sufficiently small depending on $N$, $M$, and $M'$, we can apply Lemma~\ref{lemma l2 perturbation-0}
to obtain for each $j$ and all $0<\eps<\eps_1$
\begin{align*}
\|v- \tilde v\|_{V(I_j)} &\lesssim C(j)\eps \\
\|v- \tilde v\|_{\dot S^0(I_j\times\R^n)} &\lesssim C(j)M'\\
\|v\|_{\dot S^0(I_j\times\R^n)} &\lesssim C(j)(M+ M')\\
\|(i\partial_t+\Delta)(v-\tilde v)+e\|_{\dot N^0(I_j\times\R^n)} &\lesssim C(j)\eps,
\end{align*}
provided we can prove that \eqref{close-1} and \eqref{closer l2-1} hold with $t_0$ replaced by $t_j$.
In order to verify this, we use an inductive argument.  By Strichartz, \eqref{close-1}, \eqref{error small-1},
and the inductive hypothesis,
\begin{align*}
\|v(t_j)-\tilde v(t_j)\|_{L_x^2}
&\lesssim \|v(t_0)-\tilde v(t_0)\|_{L_x^2}+\|(i\partial_t+\Delta)(v-\tilde v)+e\|_{\dot N^0([t_0,t_j]\times\R^n)}\\
&\quad + \|e\|_{\dot N^0([t_0,t_j]\times\R^n)}\\
&\lesssim M' + \sum_{k=0}^{j} C(k)\eps + \eps.
\end{align*}
Similarly, by Strichartz, \eqref{closer l2-1}, \eqref{error small-1}, and the inductive hypothesis,
\begin{align*}
\bigr\|e^{i(t-t_j)\Delta}\bigr(v(t_j)-\tilde v(t_j)\bigl)\bigl\|_{V(I)}
&\lesssim \, \bigr\|e^{i(t-t_0)\Delta}\bigr(v(t_0)-\tilde v(t_0)\bigl)\bigl\|_{V(I)}
  +\|e\|_{\dot N^0([t_0,t_j]\times\R^n)}\\
&\quad + \|(i\partial_t+\Delta)(v-\tilde v)+e\|_{\dot N^0([t_0,t_j]\times\R^n)}\\
&\lesssim \eps + \sum_{k=0}^{j} C(k)\eps.
\end{align*}
Here, $C(k)$ depends on $k$, $M$, $M'$, and $\eps_0$.  Choosing $\eps_1$ sufficiently small depending on
$N$, $M$, and $M'$, we can continue the inductive argument.
\end{proof}

The corresponding stability result for the $\dot H^1_x$-critical NLS in dimensions $3\leq n\leq 6$ is
derived by similar techniques as the ones presented above.  However, the higher dimensional case, $n>6$,
is more delicate as derivatives of the nonlinearity are merely H\"older continuous of order
$\frac{4}{n-2}$ rather than Lipschitz. A stability theory for the $\dot H^1_x$-critical NLS in higher
dimensions was established by the first two authors, \cite{TV}.  They made use of exotic Strichartz
estimates and fractional chain rule type estimates in order to avoid taking a full derivative, but still
remain energy-critical with respect to the scaling\footnote{A very similar technique was employed by K. Nakanishi,
\cite{nakanishi}, for the energy-critical non-linear Klein-Gordon equation in high dimensions.}. We
record their result below.

\begin{lemma}[$\dot H^1_x$-critical stability result]\label{h1 stability}
Let $I$ be a compact time interval and let $\tilde w$ be an approximate solution to \eqref{critical equation}
on $I\times\R^n$ in the sense that
$$
(i\partial_t+\Delta)\tilde w=|\tilde w|^{\frac{4}{n-2}}\tilde w+e
$$
for some function $e$. Assume that
\begin{align}
\|\tilde w\|_{W(I)}&\leq L \label{finite S norm} \\
\|\tilde w\|_{L_t^{\infty}\dot{H}^1_x(\ir)}&\leq E_0 \label{finite energy}
\end{align}
for some constants $L, E_0>0$. Let $t_0\in I$ and let $w(t_0)$ close to $\tilde w(t_0)$ in the sense that
\begin{align}\label{close}
\|w(t_0)-\tilde w(t_0)\|_{\dot{H}^1_x}\leq E'
\end{align}
for some $E'>0$. Assume also the smallness conditions
\begin{align}
\Bigl(\sum_N \|P_N \nabla e^{i(t-t_0)\Delta}\bigl(w(t_0)-\tilde w(t_0)\bigr)\|^2_{L_t^{\frac{2(n+2)}{n-2}}L_x^{\frac{2n(n+2)}{n^2+4}}(\ir)}\Bigr)^{1/2} &\leq \eps \label{closer} \\
\|\nabla e\|_{\dot{N}^0(\ir)}&\leq \eps \label{error small}
\end{align}
for some $0<\eps \leq \eps_2$, where $\eps_2=\eps_2(E_0, E', L)$ is a small constant. Then, there exists a solution
$w$ to \eqref{critical equation} on $\ir$ with the specified initial data $w(t_0)$ at time $t=t_0$ satisfying
\begin{align}
\|\nabla(w-\tilde w)\|_{L_t^{\frac{2(n+2)}{n-2}}L_x^{\frac{2n(n+2)}{n^2+4}}(\ir)}&\leq C(E_0 ,E', L)\bigl(\eps+\eps^{\frac{7}{(n-2)^2}}\bigr) \label{close in L^p}\\
\|w-\tilde w\|_{\dot{S}^1(\ir)}&\leq C(E_0 ,E', L)\bigl(E'+\eps+ \eps^{\frac{7}{(n-2)^2}}\bigr) \label{close in S^1}\\
\|w\|_{\dot{S}^1(\ir)}&\leq C(E_0, E', L). \label{u in S^1}
\end{align}
Here, $C(E_0,E',L) > 0$ is a non-decreasing function of $E_0,E',L$, and the dimension~$n$.
\end{lemma}

\begin{remark}\label{redundant}
By Strichartz and Plancherel, on the slab $\ir$ we have
\begin{align*}
\Bigl(\sum_N \|P_N \nabla e^{i(t-t_0)\Delta}\bigl(w(t_0)-\tilde w(t_0)&\bigr)\|^2_{L_t^{\frac{2(n+2)}{n-2}}L_x^{\frac{2n(n+2)}{n^2+4}}(\ir)}\Bigr)^{1/2}\\
&\lesssim \Bigl(\sum_N \|P_N \nabla(w(t_0)-\tilde w(t_0)\bigr)\|^2_{L_t^\infty L_x^2}\Bigr)^{1/2}\\
&\lesssim \|\nabla (w(t_0)-\tilde w(t_0)\bigr)\|_{L_t^\infty L_x^2} \\
&\lesssim E',
\end{align*}
so the hypothesis \eqref{closer} is redundant if $E'=O(\eps)$.
\end{remark}

We end this section with two well-posedness results concerning the $L_x^2$-critical and the $\dot H^1_x$-critical
nonlinear Schr\"odinger equations. More precisely, we show that control of the solution in a specific norm
($\|\cdot\|_V$ for the $L_x^2$-critical NLS and $\|\cdot\|_W$ for the $\dot H^1_x$-critical NLS),
yields control of the solution in all the $S^1$-norms .

\begin{lemma}\label{v implies s1}
Let $k=0,1$, $I$ be a compact time interval, and let $v$ be the unique solution to \eqref{l2critical} on $\ir$
obeying the bound
\begin{align}\label{bound-l2}
\|v\|_{V(I)}\le L.
\end{align}
Then, if $t_0\in I$ and $v(t_0)\in H^k_x$, we have
\begin{align}
\|v\|_{\dot S^k(I\times\R^n)}& \le C(L)\|v(t_0)\|_{\dot{H}^k_x}.\label{concl bound-l2}
\end{align}
\end{lemma}

\begin{proof}

Subdivide the interval $I$ into $N\sim (1+\frac{L}{\eta})^{\frac{2(n+2)}{n}}$ subintervals $I_j=[t_j,t_{j+1}]$ such that
\begin{align*}
\|v\|_{V(I_j)}\le\eta,
\end{align*}
where $\eta$ is a small positive constant to be chosen momentarily. By Strichartz, on each $I_j$ we obtain
\begin{align*}
\|v\|_{\dot S^k(I_j\times\R^n)}
&\lesssim \|v(t_j)\|_{\dot H^k_x} + \||\nabla|^k\bigr(|v|^{\frac 4n}v\bigl)\|_{L_{t,x}^{\frac{2(n+2)}{n+4}}(I_j\times\R^n)}\\
&\lesssim \|v(t_j)\|_{\dot H^k_x} + \|v\|_{V(I_j)}^{\frac 4n}\|v\|_{\dot S^k(I_j\times\R^n)}\\
&\lesssim \|v(t_j)\|_{\dot H^k_x} + \eta^{\frac 4n}\|v\|_{\dot S^k(I_j\times\R^n)}.
\end{align*}
Choosing $\eta$ sufficiently small, we obtain
\begin{align*}
\|v\|_{\dot S^k(I_j\times\R^n)}\lesssim \|v(t_j)\|_{\dot H^k_x}.
\end{align*}
Adding these estimates over all the subintervals $I_j$, we obtain \eqref{concl bound-l2}.
\end{proof}

\begin{lemma}\label{w implies s1}
Let $k=0,1$, $I$ be a compact time interval, and let $w$ be the unique solution to \eqref{critical equation} on $\ir$
obeying the bound
\begin{align}\label{bound-h1}
\|w\|_{W(I)}\le L.
\end{align}
Then, if $t_0\in I$ and $w(t_0)\in H^k_x$, we have
\begin{align}
\|w\|_{\dot S^k(I\times\R^n)}& \le C(L)\|w(t_0)\|_{\dot{H}^k_x}.\label{concl bound-h1}
\end{align}
\end{lemma}

\begin{proof}
The proof is similar to that of Lemma~\ref{v implies s1}. Subdivide the interval $I$ into
$N\sim (1+\frac{L}{\eta})^{\frac{2(n+2)}{n-2}}$ subintervals $I_j=[t_j,t_{j+1}]$ such that
\begin{align*}
\|w\|_{W(I_j)}\le\eta,
\end{align*}
where $\eta$ is a small positive constant to be chosen later. By Strichartz, on each $I_j$ we obtain
\begin{align*}
\|w\|_{\dot S^k(I_j\times\R^n)}
&\lesssim \|w(t_j)\|_{\dot H^k_x} + \||\nabla|^k\bigr(|w|^{\frac{4}{n-2}}w\bigl)\|_{L_{t,x}^{\frac{2(n+2)}{n+4}}(I_j\times\R^n)}\\
&\lesssim \|w(t_j)\|_{\dot H^k_x} + \|w\|_{W(I_j)}^{\frac{4}{n-2}}\|w\|_{\dot S^k(I_j\times\R^n)}\\
&\lesssim \|w(t_j)\|_{\dot H^k_x} + \eta^{\frac{4}{n-2}}\|w\|_{\dot S^k(I_j\times\R^n)}.
\end{align*}
Choosing $\eta$ sufficiently small, we obtain
\begin{align*}
\|w\|_{\dot S^k(I_j\times\R^n)}\lesssim \|w(t_j)\|_{\dot H^k_x}.
\end{align*}
The conclusion \eqref{concl bound-h1} follows by adding these estimates over all subintervals $I_j$.
\end{proof}

%
%
%
%

\section{Global well-posedness}
Our goal in this section is to prove Theorem~\ref{global wellposedness}.  We shall abbreviate the energy $E(u)$ as $E$, and
the mass $M(u)$ as $M$.

There are two ingredients to proving the existence of global solutions to \eqref{equation} in the cases described in
Theorem~\ref{global wellposedness}.  One of them is a `good' local well-posedness statement, by which we mean that
the time of existence of an $H^1_x$-solution depends only on the $H^1_x$-norm of its initial data.  The second
ingredient is an \emph{a priori} bound on the kinetic energy of the solution, i.e., its $\dot H^1_x$-norm.
These two ingredients together with the conservation of mass are sufficient to yield the existence of global solutions
via the standard iterative argument.

Before we continue with our proof, we should make a few remarks:

The existence of global $L_x^2$-solutions for \eqref{equation} when both nonlinearities are $L_x^2$-subcritical,
i.e., $0<p_1<p_2<\frac{4}{n}$, follows from the local theory for these equations and the conservation of
mass.  Indeed, the time of existence of local solutions to \eqref{equation} in this case depends only
on the $L_x^2$-norm of the initial data and global well-posedness in $L^2_x$ follows from the conservation of mass
via the standard iterative argument. For details see \cite{kato, katounique, cazenave:book}.  However,
we are interested in the existence of global $H^1_x$-solutions so, in order to iterate, we also need to control
the increment of the kinetic energy in time.

Moreover, while in the case when both nonlinearities are energy-subcritical the time of existence of $H^1_x$-solutions
depends on the $H^1_x$-norm of the initial data, in the presence of an energy-critical nonlinearity, i.e.,
$p_2=\frac{4}{n-2}$, the local theory asserts that the time of existence for $H^1_x$-solutions depends instead
on the profile of the initial data.  In order to prove a `good' local well-posedness statement in the latter case,
we will treat the energy-subcritical nonlinearity $|u|^{p_1}u$ as a perturbation to the energy-critical NLS, which
is globally wellposed (see \cite{ckstt:gwp, RV, Monica:thesis}).


\subsection{Kinetic energy control}
In this subsection we prove an \emph{a priori} bound on the kinetic energy of the solution, which is
uniform over the time of existence and which depends only on the energy and the mass of the initial data.
More precisely, we prove that for all times $t$ for which the solution is defined, we have
\begin{align}\label{kinetic control}
\|u(t)\|_{\dot H^1_x}\leq C(E,M).
\end{align}

As the energy
$$
E(u(t))=\int_{\R^n}\bigl[\tfrac12|\nabla u(t,x)|^2+\tfrac{\lambda_1}{p_1+2}|u(t,x)|^{p_1+2}
        + \tfrac{\lambda_2}{p_2+2}|u(t,x)|^{p_2+2}\bigr]dx
$$
is conserved, we immediately see that when both $\lambda_1$ and $\lambda_2$ are positive, we obtain
$$
\|\nabla u(t)\|_2^2\lesssim E,
$$
uniformly in time.

Whenever $\lambda_1<0$ and $\lambda_2>0$, we remark the inequality
$$
-\frac{|\lambda_1|}{p_1+2}|u(t,x)|^{p_1+2}+\frac{|\lambda_2|}{p_2+2}|u(t,x)|^{p_2+2}\ge
-C(\lambda_1,\lambda_2)|u(t,x)|^2,
$$
which immediately yields
$$
\|\nabla u(t)\|_2^2\lesssim E+M,
$$
uniformly over the time of existence.

When both $\lambda_1$ and $\lambda_2$ are negative, the hypotheses of Theorem~\ref{global wellposedness} also force
$0<p_1<p_2<\frac 4n$.  By interpolation and Sobolev embedding, for all times $t$ we obtain
\begin{align*}
\|u(t)\|_{p_i+2}
&\lesssim \|u(t)\|_2^{1-\frac {p_i n}{2(p_i+2)}}\|u(t)\|_{\frac {2n}{n-2}}^{\frac {np_i}{2(p_i+2)}}\\
&\lesssim M^{\frac{1}{2}-\frac {p_i n}{4(p_i+2)}}\|\nabla u(t)\|_2^{\frac {np_i}{2(p_i+2)}},
\end{align*}
where $i=1,2$.
Thus,
\begin{align}\label{perturb bdd}
\|u(t)\|_{p_i+2}^{p_i+2}\lesssim M^{1-\frac {(n-2)p_i}4}\|\nabla u(t)\|_2^{\frac {np_i}2}.
\end{align}
Next, we make use of Young's inequality,
\begin{equation} \label{young}
ab\lesssim \eps a^q+\eps^{-\frac {q'}q}b^{q'},
\end{equation}
valid for any $a,{}b,{}\eps>0$, with $1<q<\infty$ and $q'$ the dual exponent to $q$.
Taking $a=\|\nabla u(t)\|_2^{\frac{np_i}{2}}$, $b=1$, and $q=\frac{4}{np_i}$, we obtain
\begin{align*}
\|u(t)\|_{p_i+2}^{p_i+2}
&\lesssim M^{1-\frac {(n-2)p_i}4}\|\nabla u(t)\|_2^{\frac {p_in}2}\\
&\lesssim M^{1-\frac {(n-2)p_i}4}(\eps \|\nabla u(t)\|_2^2+\eps^{-\frac {np_i}{4-p_in}}).
\end{align*}
Choosing $\eps$ sufficiently small, more precisely $\eps=c^2M^{\frac {(n-2)p_i}2-2}$ for some positive constant $c\ll 1$,
we get
\begin{align*}
\|u(t)\|_{p_i+2}^{p_i+2} \leq c\|\nabla u(t)\|_2^2+C(M)
\end{align*}
Thus, by the conservation of energy,
$$
\|\nabla u(t)\|_2\le C(E,M)
$$
uniformly in $t$.

\subsection{`Good' local well-posedness}
In this subsection we prove a `good' local well-posedness statement for \eqref{equation} in the presence of an
energy-critical nonlinearity, i.e., $p_2=\frac{4}{n-2}$.  More precisely, we will find $T=T(\|u_0\|_{H^1_x})$ such that
in this case, \eqref{equation} admits a unique strong solution $u\in S^1([-T,T]\times\R^n)$ and moreover,
\begin{align}\label{good lwp bound}
\|u\|_{S^1([-T,T]\times\R^n)}\leq C(E,M).
\end{align}
In the case when both nonlinearities are energy-subcritical, this is a consequence of Proposition~\ref{subcritical lwp}.
The bound \eqref{bound wellposedness} follows easily from \eqref{good lwp bound} by subdividing the interval $I$
into subintervals of length $T$, deriving the corresponding $S^1$-bounds on each of these subintervals,
and finally adding these bounds.

To simplify notation we assume without loss of generality that $|\lambda_1|=|\lambda_2|=1$.
Moreover, by the local theory (see Section~2, specifically Proposition~\ref{critical lwp} and
Lemma~\ref{x1 implies existence}) it suffices to prove \emph{a priori} $\dot X^1$-bounds on $u$ on a time interval
whose size depends only on the $H^1_x$-norm of the initial data.  That is, we may assume that there exists a
strong solution $u$ to \eqref{equation} with $p_2=\frac{4}{n-2}$ on the slab $[-T,T]\times\R^n$ and show that
$u$ has finite $\dot X^1$-bounds on this slab as long as $T=T(\|u_0\|_{H_x^1})$ is sufficiently small.

In establishing this local well-posedness result, our approach is entirely perturbative.  More precisely,
we view the first nonlinearity $|u|^{p_1}u$ as a perturbation to the energy-critical NLS, which is globally wellposed,
\cite{ckstt:gwp, RV, Monica:thesis}.

Let therefore $w$ be the unique strong global solution to the energy-critical equation \eqref{critical equation} with initial data $w_0=u_0$ at time $t=0$. By the main results in
\cite{ckstt:gwp, RV, Monica:thesis}, we know that such a $w$ exists and moreover,
\begin{align}
\|w\|_{\sor}&\leq C(\|u_0\|_{\ho_x}).\label{v in s1}
\end{align}
Furthermore, by Lemma~\ref{w implies s1}, we also have
\begin{align*}
\|w\|_{\szr} &\leq C(\|u_0\|_{\ho_x})\|u_0\|_{L^2_x}\le C(E,M).
\end{align*}

By time reversal symmetry it suffices to solve the problem forward in time.  By \eqref{v in s1}, we can subdivide $\R_+$
into $J=J(E,\eta)$ subintervals $I_j=[t_j,t_{j+1}]$ such that
\begin{align}\label{ass v}
\|w\|_{\dot{X}^1(I_j)}\sim \eta
\end{align}
for some small $\eta$ to be specified later.

We are only interested in those subintervals $I_j$ that have a nonempty intersection with $[0,T]$. We may assume
(renumbering, if necessary) that there exists $J'<J$ such that for any $0\leq j\leq J'-1$,
$[0,T]\cap I_j \neq \varnothing$.  Thus, we can write
$$[0,T]=\cup_{j=0}^{J'-1}([0,T]\cap I_j).$$

The nonlinear evolution $w$ being small on the interval $I_j$ implies that the free evolution $e^{i(t-t_j)\Delta}w(t_j)$
is small on $I_j$. Indeed, this follows from Strichartz, Sobolev embedding, and \eqref{ass v}:
\begin{align*}
\|e^{i(t-t_j)\Delta} w(t_j)\|_{\dot X^1(I_j)}
&\le \|w\|_{\dot X^1(I_j)}+ \|\nabla (|w|^{\frac{4}{n-2}}w)\|_{L_{t,x}^{\frac{2(n+2)}{n+4}}(\ijr)}\\
&\le \|w\|_{\dot X^1(I_j)}+C\|\nabla w\|_{L^{\frac{2(n+2)}{n}}_{t,x}(\ijr)}\|w\|^{\frac 4{n-2}}_{L^{\frac {2(n+2)}{n-2}}_{t,x}(\ijr)}\\
&\le \eta + C\|w\|_{\dot X^1(I_j)}^{\frac {n+2}{n-2}}\\
&\le \eta + C\eta^{\frac{n+2}{n-2}},
\end{align*}
where $C$ is an absolute constant that depends on the Strichartz constant. Thus, taking $\eta$ sufficiently small,
for any $0\leq j\leq J'-1$, we obtain
\begin{align}\label{free evol v}
\|e^{i(t-t_j)\Delta} w(t_j)\|_{\dot X^1(I_j)}\leq 2\eta.
\end{align}

Next, we use \eqref{ass v} and \eqref{free evol v} to derive estimates on $u$. On the interval $I_0$, recalling that
$u(0)=w(0)=u_0$, we use Lemma~\ref{control by x} to estimate
\begin{align*}
\|u\|_{\dot X^1(I_0)}
&\le \|e^{it\Delta} u_0\|_{\dot X^1(I_0)}+C|I_0|^{1-\frac {p_1(n-2)}4} \|u\|_{\dot X^1(I_0)}^{p_1+1}+ C\|u\|_{\dot X^1(I_0)}^{\frac {n+2}{n-2}}\\
&\le 2\eta + C T^{1-\frac {(n-2)p_1}4}\|u\|_{\dot X^1(I_0)}^{p_1+1}
+C\|u\|_{\dot X^1(I_0)}^{\frac {n+2}{n-2}}.
\end{align*}
Assuming $\eta$ and $T$ are sufficiently small, a standard continuity argument then yields
\begin{align}\label{bdd on u}
\|u\|_{\dot X^1(I_0)}\le 4\eta.
\end{align}
Thus, \eqref{finite S norm} holds on $I:=I_0$ for $L:= 4C \eta$.
Moreover, in the previous subsection we proved that \eqref{finite energy} holds with $E_0:=C(E, M)$.
Also, as \eqref{close} holds with $E':=0$, we are in the position to apply the stability result Lemma~\ref{h1 stability}
provided the error, which in this case is the first nonlinearity, is sufficiently small.
As by H\"older and \eqref{bdd on u},
\begin{align}\label{error}
\|\nabla e\|_{\nizr}
&\lesssim T^{1-\frac{(n-2)p_1}4}\|u\|_{\dot{X}^1(I_0)}^{p_1+1} \lesssim  T^{1-\frac {(n-2)p_1}4} \eta^{p_1+1},
\end{align}
we see that by choosing $T$ sufficiently small (depending only on the energy and the mass of the initial data),
we get
$$
\|\nabla e\|_{\nizr}< \eps,
$$
where $\eps=\eps(E, M)$ is a small constant to be chosen later. Thus, taking $\eps$ sufficiently small,
the hypotheses of Lemma~\ref{h1 stability} are satisfied, which implies that the conclusion holds. In particular,
\begin{align}\label{bdd I0}
\|u-w\|_{\soizr}\le C(E,M)\eps^c
\end{align}
for a small positive constant $c$ that depends only on the dimension $n$.

By Strichartz, \eqref{bdd I0} implies
\begin{align}
\|u(t_1)-w(t_1)\|_{\ho_x}&\le C(E,M)\eps^c,\label{diff 1-1}\\
\|\propagateo(u(t_1)-w(t_1))\|_{\dot X^1(I_1)}&\le C(E,M)\eps^c.\label{diff 1-2}
\end{align}
Using \eqref{free evol v}, \eqref{diff 1-1}, \eqref{diff 1-2}, and Strichartz, we estimate
\begin{align*}
\|u\|_{\xoio}
&\le \|\propagateo u(t_1)\|_{\xoio}+ C|I_1|^{1-\frac {p_1(n-2)}4}\|u\|_{\xoio}^{p_1+1}+ C\|u\|_{\xoio}^{\frac {n+2}{n-2}}\\
&\le \|\propagateo w(t_1)\|_{\xoio}+ \|\propagateo(u(t_1)-w(t_1))\|_{\xoio} \\
&\quad + CT^{1-\frac {p_1(n-2)}4}\|u\|_{\xoio}^{p_1+1} + C\|u\|_{\xoio}^{\frac {n+2}{n-2}}\\
&\le 2\eta+C(E,M)\eps^c +C T^{1-\frac{(n-2)p_1}4}\|u\|_{\xoio}^{p_1+1}+C\|u\|_{\xoio}^{\frac {n+2}{n-2}}.
\end{align*}
A standard continuity argument then yields
$$
\|u\|_{\xoio}\le 4\eta,
$$
provided $\eps$ is chosen sufficiently small depending on $E$ and $M$, which amounts to taking $T$ sufficiently small
depending on $E$ and $M$.  Thus \eqref{error} holds with $I_0$ replaced by $I_1$ and we are again in the position
to apply Lemma~\ref{h1 stability} on $I:=I_1$ to obtain
$$
\|u-w\|_{\dot S^1(I_1)}\le C(E,M)\eps^{c^2}.
$$
By induction, for every $0\leq j\leq J'-1$ we obtain
\begin{align}\label{bdd u j}
\|u\|_{\dot X^1(I_j)}\le 4\eta,
\end{align}
provided $\eps$ (and hence $T$) is sufficiently small depending on the energy and the mass of the initial data.
Adding \eqref{bdd u j} over all $0\leq j\leq J'-1$ and recalling that $J'<J=J(E,\eta)$, we get
\begin{align}\label{x1 bound}
\|u\|_{\dot X^1([0,T])}\lesssim 4J'\eta\le C(E).
\end{align}

Next, we show that \eqref{x1 bound} implies $S^1$-control over the solution $u$ on the slab $[0,T]\times\R^n$.  This type of argument will appear
repeatedly in Section~5.  However each time, the hypotheses will be slightly different; this is why we choose not to encapsulate
it into a lemma.

By Strichartz, Lemma~\ref{control by x}, \eqref{kinetic control}, \eqref{x1 bound},
and recalling that $T=T(E,M)$, we obtain
\begin{align}
\|u\|_{\dot S^1([0,T]\times\R^n)}
&\lesssim \|u_0\|_{\dot H^1_x} + T^{1-\frac{p_1(n-2)}{4}}\|u\|_{\dot X^1([0,T])}^{1+p_1}
   +\|u\|_{\dot X^1([0,T])}^{\frac{n+2}{n-2}}\notag\\
&\leq C(E,M).\label{good lwp u s1}
\end{align}
Similarly,
\begin{align}
\|u\|_{\dot S^0([0,T]\times\R^n)}
&\lesssim \|u_0\|_{L^2_x} + T^{1-\frac{p_1(n-2)}{4}}\|u\|_{\dot X^1([0,T])}^{p_1}\|u\|_{\dot X^0([0,T])}\notag\\
&\quad +\|u\|_{\dot X^1([0,T])}^{\frac{4}{n-2}}\|u\|_{\dot X^0([0,T])}\notag\\
&\lesssim M^{\frac 12} + C(E,M)\|u\|_{\dot X^1([0,T])}^{p_1}\|u\|_{\dot S^0([0,T]\times\R^n)}\label{u s0 to bound}\\
&\quad +\|u\|_{\dot X^1([0,T])}^{\frac{4}{n-2}}\|u\|_{\dot S^0([0,T]\times\R^n)}.\notag
\end{align}
Subdividing $[0,T]$ into $N=N(E,M,\delta)$ subintervals $J_k$ such that
$$
\|u\|_{\dot X^1(J_k)}\sim \delta
$$
for some small constant $\delta>0$, the computations that lead to \eqref{u s0 to bound} now yield
\begin{align*}
\|u\|_{\dot S^0(J_k\times\R^n)}
\lesssim M^{\frac 12} + C(E,M)\delta^{p_1}\|u\|_{\dot S^0(J_k\times\R^n)}+\delta^{\frac{4}{n-2}}\|u\|_{\dot S^0(J_k\times\R^n)}.
\end{align*}
Choosing $\delta$ sufficiently small depending on $E$ and $M$, we obtain
\begin{align*}
\|u\|_{\dot S^0(J_k\times\R^n)}\leq C(E,M)
\end{align*}
on every subinterval $J_k$.  Adding these bounds over all subintervals $J_k$, we get
\begin{align}\label{good lwp u s0}
\|u\|_{\dot S^0([0,T]\times\R^n)}\leq C(E,M).
\end{align}

Collecting \eqref{good lwp u s1} and \eqref{good lwp u s0}, we obtain
$$
\|u\|_{S^1([0,T])}\le C(E,M).
$$
This concludes the proof of Theorem~\ref{global wellposedness}.

%
%
%
%

\section{Scattering results} In this section we prove Theorem~\ref{scattering}.
As before we shall abbreviate the energy $E(u)$ as $E$, and
the mass $M(u)$ as $M$.
The key ingredient is a good
spacetime bound; scattering then follows by standard techniques (see subsection 4.8).

In the case when both nonlinearities are defocusing, i.e., $\lambda_1$, $\lambda_2>0$, an \emph{a priori}
spacetime estimate for the solution is provided by the interaction Morawetz inequality, which we develop
in subsection 5.1.  In subsections 5.2 through 5.5, we upgrade this bound to a spacetime bound that implies
scattering, thus covering Case 1) of Theorem~\ref{scattering}.  Case 2) is treated in subsections 5.6 and 5.7.
In subsection 5.8 we construct the scattering states and prove the scattering result.

\subsection{The interaction Morawetz inequality}
The goal of this subsection is to prove

\begin{proposition}[Morawetz control]\label{Morawetz control proposition}
Let $I$ be a compact interval, $\lambda_1$ and $\lambda_2$ positive real numbers, and $u$ a solution to
\eqref{equation} on the slab $I\times\R^n$. Then
\begin{align}\label{Morawetz control}
\|u\|_{L_t^{n+1}L_x^{\frac{2(n+1)}{n-1}}(I\times\R^n)}\lesssim \|u\|_{L_t^\infty H^1_x(\ir)}.
\end{align}
\end{proposition}


We will derive Proposition~\ref{Morawetz control proposition} from the following:

\begin{proposition}[The interaction Morawetz estimate]\label{prop Morawetz}
Let $I$ be a compact time interval and $u$ a solution to \eqref{equation} on the slab $\ir$.  Then,
we have the following \emph{a priori} estimate:
\begin{align*}
&-(n-1)\int_I\int_{\R^{n}}\int_{\R^{n}} \Delta\bigr(\tfrac{1}{|x-y|}\bigl)|u(t,y)|^{2} |u(t,x)|^{2}dx\,dy\,dt \\
&+\sum_{i=1,2}2(n-1)\frac{\lambda_ip_i}{p_i+2}\int_I\int_{\R^{n}}\int_{\R^{n}}\frac{|u(t,y)|^{2}|u(t,x)|^{p_i+2}}
{|x-y|}dx\,dy\,dt \\
&\qquad \qquad \leq 4\|u\|_{L_t^\infty H^1_x(\ir)}^4.
\end{align*}
\end{proposition}

\begin{proof}

The calculations in this section are difficult to justify without additional assumptions (particularly,
smoothness) of the solution. This obstacle can be dealt with in the standard manner: mollify the initial
data and the nonlinearity to make the interim calculations valid and observe that the mollifications can
be removed at the end using the local well-posedness theory. For clarity, we omit these details and keep
all computations on a formal level.

We start by recalling the standard Morawetz action centered at a point.  Let $a$ be a spatial function and $u$ satisfy
\begin{align}\label{phi equation}
iu_{t}+\Delta u=\mathcal{N}
\end{align}
on $I\times \R^{n}$.  For $t\in I$, we define the Morawetz action centered at zero to be
$$
M_{a}^0(t)=2\int_{\R^{n}}a_{j}(x)\Im(\overline{u(t,x)}u_{j}(t,x))dx.
$$
A computation establishes the following
\begin{lemma}
$$
\partial_{t}M_{a}^0=\int_{\R^n} (-\Delta\Delta a)|u|^{2}
  +4\int_{\R^n} a_{jk}\Re(\overline{u_{j}}u_{k})+2\int_{\R^n} a_{j}\{\mathcal{N},u\}_{p}^{j},
$$
where we define the momentum bracket to be $\{f,g\}_p:=\Re(f\nabla \bar{g}-g\nabla\bar{f})$ and repeated
indices are implicitly summed.
\end{lemma}

Note that in the particular case when the nonlinearity is
$\mathcal{N}:=\sum_{i=1,2}\lambda_i |u|^{p_i}u$, we have
$\{\mathcal{N},u\}_p=-\sum_{i=1,2}\frac{\lambda_ip_i}{p_i+2}\nabla(|u|^{p_i+2})$.

Now let $a(x):=|x|$. For this choice of the function $a$, one should interpret $M_a^0$ as a spatial
average of the radial component of the $L^2_x$-mass current. Easy computations show that in dimension
$n\geq 3$ we have the following identities:
\begin{align*}
a_{j}(x)=&\frac{x_{j}}{|x|} \\
a_{jk}(x)=&\frac{\delta_{jk}}{|x|}-\frac{x_{j}x_{k}}{|x|^{3}} \\
\Delta a(x)=&\frac{n-1}{|x|} \\
-\Delta \Delta a(x)=&-(n-1)\Delta\Bigr(\frac{1}{|x|}\Bigl)
\end{align*}
and hence
\begin{align*}
\partial_{t}M_{a}^0
&=-(n-1)\int_{\R^n}\Delta\Bigr(\frac{1}{|x|}\Bigl) |u(x)|^{2}dx
   +4\int_{\R^n} \Bigr(\frac{\delta_{jk}}{|x|}-\frac{x_{j}x_{k}}{|x|^{3}}\Bigl) \Re(\overline{u_{j}}u_{k})(x)dx\\
&\phantom{=-(n-1)\int_{\R^n}\Delta\Bigr(\frac{1}{|x|}\Bigl) |u(x)|^{2}dx,} +2\int_{\R^n} \frac{x_{j}}{|x|}\{\mathcal{N},u\}_{p}^{j}(x)dx \\
&=-(n-1)\int_{\R^n}\Delta\Bigr(\frac{1}{|x|}\Bigl) |u(x)|^{2}dx
   +4\int_{\R^n} \frac{1}{|x|} |\nabla_{0}u(x)|^{2}dx\\
&\phantom{=-(n-1)\int_{\R^n}\Delta\Bigr(\frac{1}{|x|}\Bigl) |u(x)|^{2}dx,}
   +2\int_{\R^n} \frac{x}{|x|} \{\mathcal{N},u\}_{p}(x)dx,
\end{align*}
where we use $\nabla_{0}$ to denote the complement of the radial portion of the gradient, that is,
$\nabla_{0}:=\nabla-\frac{x}{|x|}(\frac{x}{|x|}\cdot\nabla)$.

We may center the above argument at any other point $y\in \R^{n}$. Choosing $a(x):=|x-y|$, we define the
Morawetz action centered at $y$ to be
$$
M_{a}^y(t)=2\int_{\R^{n}}\frac{x-y}{|x-y|}\Im(\overline{u(t,x)}\nabla u(t,x))dx.
$$
The same considerations now yield
\begin{align*}
\partial_{t}M_a^{y}
&=-(n-1)\int_{\R^n}\Delta\Bigr(\frac{1}{|x-y|}\Bigl) |u(x)|^{2}dx
   +4\int_{\R^n} \frac{1}{|x-y|} |\nabla_{y}u(x)|^{2}dx\\
&\phantom{\qquad=-(n-1)\int_{\R^n}\Delta\Bigr(\frac{1}{|x|}\Bigl) |u(x)|^{2}dx} +2\int_{\R^n} \frac{x-y}{|x-y|}
\{\mathcal{N},u\}_{p}(x)dx.
\end{align*}

We are now ready to define the interaction Morawetz potential, which is a way of quantifying how mass is
interacting with (moving away from) itself:
\begin{align*}
M^{interact}(t)
&:=\int_{\R^n} |u(t,y)|^{2}M_a^{y}(t)dy \\
&=2\Im \int_{\R^n} \int_{\R^n}
     |u(t,y)|^{2}\frac{x-y}{|x-y|}\nabla u(t,x)\overline{u(t,x)}dx\,dy.
\end{align*}
One obtains immediately the easy estimate
$$
|M^{interact}(t)| \leq 2 \|u(t)\|_{L_{x}^{2}}^{3} \|u(t)\|_{\dot{H}_{x}^{1}}.
$$

Calculating the time derivative of the interaction Morawetz potential, we get the following virial-type
identity:
\begin{align}
\partial_{t}M^{interact}
=&-(n-1)\int_{\R^n} \int_{\R^n} \Delta\Bigr(\frac{1}{|x-y|}\Bigl)|u(y)|^{2}|u(x)|^{2}dx\,dy  \label{vi1} \\
 &+4\int_{\R^n} \int_{\R^n} \frac{|u(y)|^{2}|\nabla_{y}u(x)|^{2}}{|x-y|}dx\,dy  \label{vi2} \\
 &+2\int_{\R^n} \int_{\R^n} |u(y)|^{2}\frac{x-y}{|x-y|} \{\mathcal{N},u\}_{p}(x)dx\,dy  \label{vi3} \\
 &+2 \int_{\R^n} \partial_{y_{k}} \Im(u\bar{u}_k)(y)M_a^{y}dy  \label{vi4} \\
 &+4\Im \int_{\R^n} \int_{\R^n} \{\mathcal{N},u\}_{m}(y)\frac{x-y}{|x-y|}\nabla u(x)\overline{u(x)}dx\,dy, \label{vi5}
\end{align}
where the mass bracket is defined to be $\{f,g\}_m:=\Im(f\bar{g})$. Note that for the nonlinearity of interest
$\mathcal{N}:=\sum_{i=1,2}\lambda_i |u|^{p_i}u$, we have $\{\mathcal{N},u\}_m=0$.\\

As far as the terms in the above identity are concerned, we will establish
\begin{lemma}\label{lemmaviterms}
\eqref{vi4} $\geq$ --\eqref{vi2}.
\end{lemma}

Thus, integrating with respect to time over a compact interval $I$, we get

\begin{proposition}[General interaction Morawetz inequality] \label{intmorineq}
\begin{align*}
&-(n-1)\int_I\int_{\R^n} \int_{\R^n} \Delta\Bigr(\frac{1}{|x-y|}\Bigl)|u(y)|^{2}|u(x)|^{2}dx\,dy\,dt \\
&\phantom{(n-1)(n-3)\int_I} +2\int_I \int_{\R^n} \int_{\R^n} |u(t,y)|^{2}\frac{x-y}{|x-y|}\{\mathcal{N},u\}_{p}(t,x)dx\,dy\,dt \\
&\phantom{(n-1)(n)} \leq 4\|u\|_{L_{t}^{\infty}L_{x}^{2}(I\times\R^{n})}^{3}
\|\nabla u\|_{L_{t}^{\infty}L_x^2(I\times\R^{n})} \\
&\phantom{(n-1)(n-3)\int_{I_{0}}}+4\int_I\int_{\R^n}\int_{\R^n}|\{\mathcal{N},u\}_{m}(t,y)u(t,x)\nabla u(t,x)|dx\,dy\,dt.
\end{align*}
\end{proposition}

Note that in the particular case $\mathcal{N}=\sum_{i=1,2}\lambda_i |u|^{p_i}u$, after performing an integration by
parts in the momentum bracket term, the inequality becomes
\begin{align*}
&-(n-1)\int_I\int_{\R^n} \int_{\R^n} \Delta\Bigr(\frac{1}{|x-y|}\Bigl)|u(y)|^{2}|u(x)|^{2}dx\,dy\,dt\\
&+2(n-1)\sum_{i=1,2}\frac{\lambda_ip_i}{p_i+2}\int_I\int_{\R^{n}}\int_{\R^{n}}\frac{|u(t,y)|^{2}|u(t,x)|^{p_i+2}}{|x-y|}dx\,dy\,dt\\
&\quad\quad\quad \leq 4\|u\|_{L_{t}^{\infty}L_{x}^{2}(I\times\R^{n})}^{3} \|\nabla u\|_{L_{t}^{\infty}L_x^2(I\times\R^{n})},
\end{align*}
which proves Proposition \ref{prop Morawetz}.

We turn now to the proof of Lemma \ref{lemmaviterms}. We write
$$
\eqref{vi4}=4\int_{\R^n} \int_{\R^n} \partial_{y_{k}}
\Im(u(y)\overline{u_{k}(y)})\frac{x_j-y_j}{|x-y|}\Im(\overline{u(x)}u_j(x))dx\,dy,
$$
where we sum over repeated indices. We integrate by parts moving $\partial_{y_k}$ to the unit vector
$\frac{x-y}{|x-y|}$. Using the identity
$$
\partial_{y_k}\Bigl(\frac{x_j-y_j}{|x-y|}\Bigr)=-\frac{\delta_{kj}}{|x-y|}+\frac{(x_k-y_k)(x_j-y_j)}{|x-y|^3}
$$
and the notation $p(x)=2\Im(\overline{u(x)}\nabla u(x))$ for the momentum density, we rewrite
\eqref{vi4} as
$$
-\int_{\R^n} \int_{\R^n} \Bigl[ p(y)p(x)-\Bigl( p(y)\frac{x-y}{|x-y|} \Bigr) \Bigl(
p(x)\frac{x-y}{|x-y|} \Bigr) \Bigr] \frac{dx\,dy}{|x-y|}.
$$
Note that the quantity between the square brackets represents the inner product between the projections
of the momentum densities $p(x)$ and $p(y)$ onto the orthogonal complement of $(x-y)$. But
\begin{align*}
|\pi_{(x-y)^\perp}p(y)|
   &=\Bigl|p(y)-\frac{x-y}{|x-y|}\Bigl(\frac{x-y}{|x-y|}p(y)\Bigr)\Bigr|
    =2|\Im(\overline{u(y)}\nabla_x u(y))| \\
   &\leq 2|u(y)||\nabla_x u(y))|.
\end{align*}
As the same estimate holds when we switch $y$ and $x$, we get
\begin{align*}
\eqref{vi4}
    &\geq -4\int_{\R^n} \int_{\R^n} |u(y)||\nabla_x u(y))| |u(x)| |\nabla_y u(x)|\frac{dx\,dy}{|x-y|} \\
    &\geq -2 \int_{\R^n} \int_{\R^n} \frac{|u(y)|^{2}|\nabla_{y}u(x)|^{2}}{|x-y|}dx\,dy
       -2 \int_{\R^n} \int_{\R^n} \frac{|u(x)|^{2}|\nabla_{x}u(y)|^{2}}{|x-y|}dx\,dy \\
    &\geq -\eqref{vi2}.
\end{align*}
which proves Lemma \ref{lemmaviterms}.
\end{proof}

We return now to the proof of Proposition~\ref{Morawetz control proposition}.  We are assuming that both
nonlinearities are defocusing, i.e., $\lambda_1>0$ and $\lambda_2>0$.
An immediate corollary of Proposition \ref{prop Morawetz} in this case is the following estimate:
\begin{align}\label{first control}
-\int_I\int_{\R^n} \int_{\R^n} \Delta\Bigr(\frac{1}{|x-y|}\Bigl)|u(t,y)|^{2}|u(t,x)|^{2}dx\,dy\,dt
\lesssim \|u\|_{L_t^\infty H^1_x(\ir)}^4.
\end{align}

In dimension $n=3$, we have $-\Delta(\frac{1}{|x|})=4\pi\delta$, so \eqref{first control} yields
$$
\int_I\int_{\R^3}|u(t,x)|^4dx\,dt \lesssim \|u\|_{L_t^\infty H^1_x(I\times\R^3)}^4,
$$
which proves Proposition \ref{Morawetz control proposition} in this case.

In dimension $n\geq 4$, an easy computation shows $-\Delta(\frac{1}{|x|})=\frac{n-3}{|x|^3}$.
As convolution with $\frac{1}{|x|^3}$ is (apart from a constant factor) essentially the fractional integral operator
$|\nabla|^{-(n-3)}$, from \eqref{first control} we derive
\begin{align}\label{u^2 control}
\|\nabla^{-\frac{n-3}{2}}|u|^2\|_{L_{t,x}^2(I\times\R^n)}\lesssim \|u\|^2_{L_t^\infty H^1_x(I\times\R^n)}.
\end{align}
A consequence of \eqref{u^2 control} is
\begin{align}\label{u control}
\|\nabla^{-\frac{n-3}{4}}u\|_{L_{t,x}^4(I\times\R^n)}\lesssim
\|u_0\|_{L_t^\infty H^1_x(I\times\R^n)},
\end{align}
as can be seen by taking $f=u$ in the following

\begin{lemma}
\begin{align}\label{u^2-->u}
\||\nabla|^{-\frac{n-3}{4}}f\|_4\lesssim \||\nabla|^{-\frac{n-3}{2}}|f|^2\|_2^{1/2}.
\end{align}
\end{lemma}

\begin{proof}
As $|\nabla|^{-\frac{n-3}{4}}$ and $|\nabla|^{-\frac{n-3}{2}}$ correspond to convolutions with positive
kernels, it suffices to prove \eqref{u^2-->u} for a positive Schwartz function $f$. For such an $f$, we
will show the pointwise inequality
\begin{align}\label{square function}
S(|\nabla|^{-\frac{n-3}{4}}f)(x)\lesssim [(|\nabla|^{-\frac{n-3}{2}}|f|^2)(x)]^{1/2},
\end{align}
where $S$ denotes the Littlewood-Paley square function $Sf := (\sum_N |P_N f|^2)^{1/2}$. Clearly \eqref{square function} implies
\eqref{u^2-->u}:
$$
\||\nabla|^{-\frac{n-3}{4}}f\|_4\lesssim \|S(|\nabla|^{-\frac{n-3}{4}}f)\|_4
    \lesssim\|(|\nabla|^{-\frac{n-3}{2}}|f|^2)^{1/2}\|_4\lesssim\||\nabla|^{-\frac{n-3}{2}}|f|^2\|_2^{1/2}.
$$

In order to prove \eqref{square function} we will estimate each of the dyadic pieces,
$$
P_N(|\nabla|^{-\frac{n-3}{4}}f)(x)=\int e^{2\pi i x\xi}\hat{f}(\xi)|\xi|^{-\frac{n-3}{4}}m(\xi/N)d\xi,
$$
where $m(\xi):=\phi(\xi)-\phi(2\xi)$ in the notation introduced in Section~2. As
$|\xi|^{-\frac{n-3}{4}}m(\xi/N)\sim N^{-\frac{n-3}{4}}\tilde{m}(\xi/N)$ for $\tilde{m}$ a multiplier
with the same properties as $m$, we have
\begin{align*}
P_N(|\nabla|^{-\frac{n-3}{4}}f)(x)
 &\sim f* \bigl( N^{-\frac{n-3}{4}} [\tilde{m}(\xi/N)]\hbox{\Large$\check{\ }$}(x) \bigr)
    = N^{\frac{3(n+1)}{4}} f*\check{\tilde{m}}(Nx)  \\
 &= N^{\frac{3(n+1)}{4}} \int f(x-y) \check{\tilde{m}}(Ny) dy.
\end{align*}
Hence,
$$
|P_N(|\nabla|^{-\frac{n-3}{4}}f)(x)|\lesssim N^{\frac{3(n+1)}{4}} \int_{|y|\leq N^{-1}} f(x-y)dy
$$
and a simple application of Cauchy-Schwartz yields
\begin{align*}
S(|\nabla|^{-\frac{n-3}{4}}f)(x)
  &=\Bigl(\sum_{N}|P_N(|\nabla|^{-\frac{n-3}{4}}f)(x)|^2\Bigr)^{1/2} \\
  &\lesssim \Bigl(\sum_{N} N^{\frac{3(n+1)}{2}} \bigl|\int_{|y|\leq N^{-1}} f(x-y)dy\bigr|^2\Bigr)^{1/2} \\
  &\lesssim \Bigl(\sum_{N} N^{\frac{3(n+1)}{2}} N^{-n} \int_{|y|\leq N^{-1}} |f(x-y)|^2 dy \Bigr)^{1/2} \\
  &\lesssim \Bigl(\sum_{N} N^{\frac{n+3}{2}}\int_{|y|\leq N^{-1}} |f(x-y)|^2 dy \Bigr)^{1/2}.
\end{align*}
As
$$
\sum_{N} N^{\frac{n+3}{2}}\chi_{\{|y|\leq N^{-1}\}}(y)\lesssim \sum_{|y|\leq N^{-1}} N^{\frac{n+3}{2}}
   \lesssim |y|^{-\frac{n+3}{2}},
$$
we get
$$
S(|\nabla|^{-\frac{n-3}{4}}f)(x)\lesssim \Bigl( \int \frac{|f(x-y)|^2}{|y|^{\frac{n+3}{2}}}
dy\Bigr)^{1/2}
    \sim [(|\nabla|^{-\frac{n-3}{2}}|f|^2)(x)]^{1/2},
$$
and the claim follows.
\end{proof}

Proposition~\ref{Morawetz control proposition} in the case $n\geq 4$ follows by interpolation between
\eqref{u control} and the bound on the kinetic energy,
$$
\|\nabla u\|_{L_t^\infty L_x^2}\lesssim E^{\frac{1}{2}},
$$
which is an immediate consequence of the conservation of energy when both nonlinearities are defocusing.

\subsection{Global bounds in the case $\frac 4n=p_1<p_2<\frac 4{n-2}$ and $\lambda_1, \lambda_2>0$}

In this subsection, we upgrade the spacetime bound given by Proposition~\ref{Morawetz control proposition} to a
spacetime bound that implies scattering for solutions to \eqref{equation} with $\frac 4n=p_1<p_2<\frac 4{n-2}$ and
$\lambda_1$, $\lambda_2$ positive.  To simplify notation, we assume, without loss of generality, that
$\lambda_1=\lambda_2=1$.

Our approach in proving the desired spacetime estimate is perturbative.  More precisely, we view the second
nonlinearity as a perturbation to the $L_x^2$-critical NLS.  The result we obtain is thus conditional upon a satisfactory
global theory for the $L_x^2$-critical problem; specifically, it holds under the Assumption~\ref{assumption}.

By Proposition \ref{Morawetz control proposition} and the conservation of energy and mass, we get
$$
\|u\|_{Z(\R)}\lesssim \|u\|_{L_t^{\infty}H_x^1(\R\times\R^n)}\le C(E,M).
$$
Let $\eps>0$ be a small constant to be chosen later.  Split $\R$ into $J=J(E,M,\eps)$ subintervals $I_j$,
$0\le j\le J-1$, such that
$$
\|u\|_{Z(I_j)}\sim \eps.
$$
We will show that on each slab $I_j\times\R^n$, $u$ obeys good Strichartz bounds.

Our analysis in this subsection will be carried out in the following space:
On the slab $I\times\R^n$, we define
$$
\dot Y^0(I):=L^{2+\frac 1{\theta}}_tL_x^{\frac{2n(2\theta+1)}{n(2\theta+1)-4\theta}}(I\times\R^n)\cap V(I),
$$
where $\theta>0$ is chosen sufficiently large so that \eqref{minter} holds.  This allows us to control the second
nonlinearity in terms of the $\dot Y^0$-norm, the $H^1_x$-norm, and the $Z$-norm:
\begin{align}
\||u|^{p_2}u\|_{\nnij}
&\lesssim \|u\|_{L^{2+\frac 1{\theta}}_tL_x^{\frac{2n(2\theta+1)}{n(2\theta+1)-4\theta}}(I_j\times \R^n)}
  \|u\|_{L_t^{\infty}H_x^1(I_j\times\R^n)}^{\alpha(\theta)+\beta(\theta)}\|u\|_{Z(I_j)}^c\nonumber\\
&\le C(E,M)\eps^c \|u\|_{\dot Y^0(I_j)},\label{small up2}
\end{align}
for all $0\leq j\leq J-1$ and a constant $c:=\frac{n+1}{2(2\theta+1)}$.

In what follows, we fix an interval $I_{j_0}=[a,b]$ and prove that $u$ obeys good Strichartz estimates on the slab
$I_{j_0}\times\R^n$.  In order to do so, we view the solution $u$ as a perturbation to solutions to the
$L_x^2$-critical NLS,
$$
\begin{cases}
iv_t+\Delta v=|v|^{\frac 4n} v\\
v(a)=u(a).
\end{cases}
$$
As this initial value problem is globally wellposed in $H_x^1$, and by Assumption~\ref{assumption} and Lemma~\ref{v implies s1}, the unique solution
enjoys the global spacetime bound
$$
\|v\|_{\sz}\le C(M),
$$
we can subdivide $\R$ into $K=K(M, \eta)$ subintervals $J_k$ such that on each $J_k$,
\begin{align}\label{scat1 y0 bound}
\|v\|_{\dot Y^0(J_k)}\sim \eta,
\end{align}
for a small constant $\eta>0$ to be chosen later.  Of course, we are only interested in those $J_k=[t_k,t_{k+1}]$
which have a nonempty intersection with $I_{j_0}$.  Without loss of generality, we may assume that
$$
[a,b]=\cup_{k=0}^{K'-1}J_k, \ \ t_0=a, \ \ t_{K'}=b.
$$

The nonlinear evolution $v$ being small on $J_k\times\R^n$ implies that the linear evolution $e^{i(t-t_k)\Delta}v(t_k)$
is small as well.  Indeed, by Strichartz and \eqref{scat1 y0 bound}, we get
\begin{align*}
\|e^{i(t-t_k)\Delta}v(t_k)\|_{\dot Y^0(J_k)}
&\le \|v\|_{\dot Y^0(J_k)}+C\||v|^{\frac 4n}v\|_{L_{t,x}^{\frac{2(n+2)}{n+4}}(J_k\times\R^n)}\\
&\le \eta+C\|v\|_{V(J_k)}^{1+\frac 4n}\\
&\le \eta+C\eta^{1+\frac 4n}.
\end{align*}
Choosing $\eta$ sufficiently small, this implies
\begin{align}\label{scat1 free small}
\|e^{i(t-t_k)\Delta}v(t_k)\|_{\dot Y^0(J_k)}\leq 2\eta.
\end{align}

Next, we will compare $u$ to $v$ on the slab $[t_0, t_1]\times\R^n$ via the $L_x^2$-stability lemma
and use the result as an input in the conditions one needs to check in order to compare $u$ to $v$ (again, via
Lemma~\ref{l2 stability}) on the slab $[t_1, t_2]\times\R^n$.  By iteration, we will derive bounds on $u$ from
bounds on $v$ on all slabs $J_k\times\R^n$, $0\leq k< K'$, and hence, we will obtain an estimate on the $\dot S^0$-norm
of $u$ on $I_{j_0}\times\R^n$.

We present the details below.  Recalling that $u(t_0)=v(t_0)$, by Strichartz, \eqref{small up2}, and
\eqref{scat1 free small}, we get
\begin{align*}
\|u\|_{\dot Y^0(J_0)}
&\le \|e^{i(t-t_0)\Delta}u(t_0)\|_{\dot Y^0(J_0)}+C\||u|^{\frac 4n}u\|_{\nnjn}+C\||u|^{p_2}u\|_{\nnjn}\\
&\le 2\eta+C\|u\|_{\dot Y^0(J_0)}^{1+\frac 4n}+C(E,M)\varepsilon^c\|u\|_{\dot Y^0(J_0)},
\end{align*}
which, by a standard continuity argument, yields
\begin{align}\label{scat1 u yo jo}
\|u\|_{\dot Y^0(J_0)}\le 4\eta,
\end{align}
provided $\eta$ and $\eps=\eps(E,M)$ are chosen sufficiently small.  Therefore, in order to apply
Lemma~\ref{l2 stability}, we just need to check that the error term $e=|u|^{p_2}u$ is small in $\nnjn$.
As by \eqref{small up2},
\begin{equation}\label{errorjn}
\|e\|_{\nnjn}\le C(E,M)\eps^c\|u\|_{\dot Y^0(J_0)}\le C(E,M)\eta\eps^c,
\end{equation}
choosing $\eps$ sufficiently small depending only on $E$ and $M$, we obtain
$$
\|u-v\|_{\snjn}\le \eps^{c/2}.
$$
By Strichartz, this implies
\begin{align}
\|u(t_1)-v(t_1)\|_{L_x^2}\le \eps^{c/2},\label{scat1 u-v 1}\\
\|e^{i(t-t_1)\Delta}(u(t_1)-v(t_1))\|_{\dot Y^0(J_1)}\lesssim \eps^{c/2}\label{scat1 u-v 2}.
\end{align}

Before turning to the second interval, $J_1$, let us also remark the following $\dot S^1$-control on $u$ on the slab
$J_0\times\R^n$.  Indeed, by Strichartz, \eqref{small up2}, and \eqref{scat1 u yo jo}, we have
\begin{align*}
\|u\|_{\sfjn}
&\lesssim \|u(a)\|_{\dot H^1_x} +\|u\|_{V(J_0)}^{\frac 4n}\|u\|_{\sfjn}+\||u|^{p_2}u\|_{\nfjn}\\
&\lesssim C(E)+(4\eta)^{\frac 4n}\|u\|_{\sfjn}+C(E,M)\eps^c\|u\|_{\sfjn},
\end{align*}
which for $\eta$ and $\eps=\eps(E,M)$ sufficiently small yields
$$
\|u\|_{\sfjn}\le C(E).
$$

Next, we use \eqref{scat1 u-v 1} and \eqref{scat1 u-v 2} to estimate $u$ on the slab $J_1\times\R^n$.  By Strichartz,
\eqref{small up2}, \eqref{scat1 free small}, and \eqref{scat1 u-v 2}, we estimate
\begin{align*}
\|u\|_{\dot Y^0(J_1)}&\le \|e^{i(t-t_1)\Delta}v(t_1)\|_{\dot
Y^0(J_1)}+\|e^{i(t-t_1)\Delta}(u(t_1)-v(t_1))\|_{\dot Y^0(J_1)}\\
&\qquad +C\|u\|_{\dot Y^0(J_1)}^{1+\frac
4n}+C(E,M)\varepsilon^c\|u\|_{\dot Y^0(J_1)}\\
&\le 2\eta+\varepsilon^{c/2} +C\|u\|_{\dot Y^0(J_1)}^{1+\frac
4n}+C(E,M)\varepsilon^c\|u\|_{\dot Y^0(J_1)}.
\end{align*}
A standard continuity argument yields
$$
\|u\|_{\dot Y^0(J_1)}\le 4\eta,
$$
provided $\eta$ and $\eps=\eps(E,M)$ are chosen sufficiently small.  This implies that
the error, i.e., $|u|^{p_1}u$, obeys \eqref{errorjn} with $J_0$ replaced by $J_1$. Choosing $\eps$ sufficiently small
depending on $E$ and $M$, we can apply Lemma \ref{l2 stability} to derive
$$
\|u-v\|_{\dot S^0(J_1\times \R^n)}\le \eps^{c/4}.
$$
The same arguments as before also yield
$$
\|u\|_{\dot S^1(J_1\times\R^n)}\le C(E).
$$

By induction, taking $\eps$ smaller with each step, for each $0\le k\le K'-1$ we obtain
$$
\|u-v\|_{\dot S^0(J_k\times\R^n)}\le \eps^{c/2^{k+1}}
$$
and
$$
\|u\|_{\dot S^1(J_k\times\R^n)}\le C(E).
$$
Adding these estimates over all the intervals $J_k$ which have a nontrivial intersection with $I_{j_0}$, we obtain
\begin{align*}
\|u\|_{\dot S^0(I_{j_0}\times\R^n)}
&\lesssim \|v\|_{\dot S^0(I_{j_0}\times\R^n)}+\sum_{k=0}^{K'-1}\|u-v\|_{\dot S^0(J_k\times\R^n)}\le C(E,M),\\
\|u\|_{\dot S^1(I_{j_0}\times\R^n)}
&\lesssim \sum_{k=0}^{K'-1}\|u\|_{\dot S^1(J_k\times\R^n)}\le C(E,M).
\end{align*}

As the interval $I_{j_0}$ was arbitrarily chosen, we get
\begin{align*}
\|u\|_{\sz(\R\times\R^n)}&\lesssim \sum_{j=0}^{ J-1}\|u\|_{\dot S^0(I_j\times\R^n)}\lesssim J C(E,M)\le C(E,M),\\
\|u\|_{\so(\R\times\R^n)}&\lesssim \sum_{j=0}^{ J-1}\|u\|_{\dot S^1(I_j\times\R^n)}\lesssim J C(E,M)\le C(E,M),
\end{align*}
and hence
$$
\|u\|_{S^1(\R\times\R^n)}\leq C(E,M).
$$

\subsection{Ode to Morawetz}

In this subsection we upgrade the bound \eqref{Morawetz control} to good Strichartz bounds in the case
$\frac 4n<p_1<p_2<\frac 4{n-2}$ and $\lambda_1$, $\lambda_2$ positive.  For simplicity, we only derive spacetime bounds
for solutions to the initial value problem
\begin{equation}\label{scat2 eq}
\begin{cases}
iu_t+\Delta u=|u|^pu\\
u(0)=u_0\in H^1_x,
\end{cases}
\end{equation}
with $\frac 4n<p<\frac 4{n-2}$.  Treating the NLS with finitely many such nonlinearities introduces only notational
difficulties.

Scattering in $H^1_x$ for solutions to \eqref{scat2 eq} was first proved by J. Ginibre and G. Velo, \cite{gv:scatter}.
Below, we present a new, simpler proof that relies on the interaction Morawetz estimate.

By Theorem~\ref{global wellposedness}, the initial value problem \eqref{scat2 eq} is globally wellposed.
Moreover, by Proposition~\ref{prop Morawetz} and the conservation of energy (E) and mass (M), the unique global
solution satisfies
$$
\|u\|_{L^{n+1}_tL_x^{\frac{2(n+1)}{n-1}}(\rr)}\lesssim \|u\|_{L_t^\infty H^1_x(\rr)}\leq C(E,M).
$$
Let $\eta>0$ be a small constant to be chosen later and divide $\R$ into $J=J(E,M,\eta)$ subintervals
$I_j=[t_j,t_{j+1}]$ such that
$$
\|u\|_{L^{n+1}_tL_x^{\frac{2(n+1)}{n-1}}(\ijr)}\sim \eta.
$$

Then, on $I_j$ $u$ satisfies the integral equation
$$
u(t)=e^{i(t-t_j)\Delta}u(t_j)-i\int_{t_j}^t e^{i(t-s)\Delta}\bigl(|u|^pu\bigr)(s)ds.
$$
By Strichartz,
$$
\|u\|_{S^1(I_j\times \R^n)}
\lesssim \|u(t_j)\|_{H^1_x}+\||u|^pu\|_{L^2_t W_x^{1,\frac{2n}{n+2}}(\ijr)},
$$
which by Lemma~\ref{m interpolation lemma} yields
\begin{align*}
\|u\|_{S^1(I_j\times \R^n)}
\lesssim \|u\|_{L_t^\infty H^1_x(\ijr)} + \eta^{\frac{n+1}{2(2\theta+1)}}\|u\|_{L_t^\infty H^1_x(\ijr)}^{\alpha(\theta)+\beta(\theta)}
  \|u\|_{S^1(I_j\times\R^n)},
\end{align*}
provided $\theta$ is chosen sufficiently large. From the conservation of energy and mass, and choosing $\eta$
sufficiently small (depending on $E$ and $M$), we get
\begin{align*}
\|u\|_{S^1(I_j\times\R^n)}\le C(E,M).
\end{align*}
Summing these bounds over all intervals $I_j$, we obtain
$$
\|u\|_{S^1(\R\times\R^n)}\lesssim J C(E,M) \le C(E,M).
$$

\subsection{Global bounds in the case $\frac 4n<p_1<p_2=\frac 4{n-2}$ and $\lambda_1,\lambda_2>0$}

In this subsection we upgrade the spacetime estimate given by the interaction Morawetz inequality,
\eqref{Morawetz control}, to spacetime bounds that imply scattering.  The approach is similar to that used in
subsection 5.2; this time, we view the first nonlinearity as a perturbation to the energy-critical NLS.
Without loss of generality, we may assume $\lambda_1=\lambda_2=1$.

Let $\eps$ be a small constant to be chosen later.  As by Proposition~\ref{Morawetz control proposition} and
the conservation of energy and mass,
$$
\|u\|_{Z(\R)}=\|u\|_{L_t^{n+1}L_x^{\frac {2(n+1)}{n-1}}(\rr)}\lesssim \|u\|_{L_t^\infty H_x^1(\rr)}\leq C(E,M),
$$
we can split $\R$ into $J=J(E,M,\eps)$ intervals $I_j$, $0\leq j\leq J -1$, such that
\begin{align}\label{ass u mor}
\|u\|_{Z(I_j)} \sim \eps.
\end{align}
We will show that on each slab $I_j\times\R^n$, $u$ obeys good Strichartz bounds.

For a spacetime slab $\ir$, we define the spaces
\begin{align*}
\dot Y^0(I)&:=L^{2+\frac 1{\theta}}_tL_x^{\frac {2n(2\theta+1)}{n(2\theta+1)-4\theta}}(\ir)\cap L_{t,x}^{\frac{2(n+2)}{n}}(\ir)\cap L_t^{\frac{2(n+2)}{n-2}} L_x^{\frac {2n(n+2)}{n^2+4}}(\ir),\\
\dot Y^1(I)&:= \{u:\, \nabla u\in \dot Y^0(I)\}, \quad \text{and}
\quad Y^1(I):=\dot Y^0(I)\cap \dot Y^1(I)
\end{align*}
with the usual topology.  Here, $\theta$ is a sufficiently large constant so that \eqref{minter} holds; that is,
on each slab $I_j\times\R^n$,
\begin{align}\label{nonlin 1}
\|\nabla \bigl(|u|^{p_1}u\bigr)\|_{L^2_tL_x^{\frac {2n}{n+2}}(I_j\times \R^n)}
&\lesssim \|\nabla u\|_{L_t^{2+\frac 1{\theta}}L_x^{\frac {2n(2\theta+1)}{n(2\theta+1)-4\theta}}(I_j\times\R^n)}
     \|u\|_{Z(I_j)}^{\frac{n+1}{2(2\theta+1)}}\|u\|_{L^{\infty}_tH_x^1(I_j\times\R^n)}^{\alpha(\theta)+\beta(\theta)} \notag\\
&\le C(E,M)\eps^{\frac{n+1}{2(2\theta+1)}}\|u\|_{\dot Y^1(I_j)}.
\end{align}

Moreover, in this notation we also have (by Sobolev embedding)
\begin{align}\label{nonlin 2}
\|\nabla\bigl(|u|^{\frac 4{n-2}}u\bigr)\|_{L_{t,x}^{\frac{2(n+2)}{n+4}}(I_j\times\R^n)}
\lesssim \|\nabla u\|_{L_{x,t}^{\frac{2(n+2)}n}(I_j\times\R^n)}\|u\|_{L_{x,t}^{\frac{2(n+2)}{n-2}}(I_j\times\R^n)}^{\frac{4}{n-2}}
\lesssim \|u\|_{\dot Y^1(I_j)}^{\frac {n+2}{n-2}},
\end{align}
for each $0\leq j\leq J-1$.

Fix $I_{j_0}:=[a,b]$.  On the slab $I_{j_0}\times\R^n$, we treat the first nonlinearity as a perturbation to the
energy-critical NLS
\begin{equation}\label{out at a}
\begin{cases}
iw_t+\Delta w= |w|^{\frac 4{n-2}}w\\
w(a)=u(a).
\end{cases}
\end{equation}
By the global well-posedness results in \cite{ckstt:gwp, RV, Monica:thesis}, there exists a unique global solution
$w$ to \eqref{out at a} with initial data $u(a)$ at time $t=a$ and moreover,
\begin{align}
\|w\|_{\sor}\le C(\|u(a)\|_{\hox})\le C(E)\label{v(a) S bdd}.
\end{align}
By Lemma~\ref{w implies s1}, \eqref{v(a) S bdd} implies
\begin{align*}
\|w\|_{\szr}\le C(\|u(a)\|_{\hox})\|u(a)\|_{L_x^2} \le C(E,M).
\end{align*}

Given \eqref{v(a) S bdd}, we can split $\R$ into $K=K(E,\eta)$ subintervals $J_k=[t_k,t_{k+1}]$ such that
\begin{align}\label{ass v 2}
\|w\|_{\dot Y^1(J_k)}\sim \eta,
\end{align}
for some small constant $\eta>0$ to be chosen later.  Using \eqref{ass v 2} and Strichartz, it is
easy to see that the free evolution $e^{i(t-t_k)\Delta}w(t_k)$ is small on $J_k$ for every $0\leq k\leq K-1$. Indeed,
\begin{align}\label{lin evol small}
\|e^{i(t-t_k)\Delta}w(t_k)\|_{\dot Y^1(J_k)}
&\le \|w\|_{\dot Y^1(J_k)}+\Bigl\|\int_{t_k}^te^{i(t-s)\Delta}\bigl(|w|^{\frac 4{n-2}}w\bigr)(s)ds\Bigr\|_{\dot Y^1(J_k)}\notag\\
&\le \|w\|_{\dot Y^1(J_k)}+ C\|w\|_{\dot Y^1(J_k)}^{\frac{n+2}{n-2}}\notag\\
&\le 2\eta,
\end{align}
provided $\eta$ is chosen sufficiently small.

Of course, we are only interested in those $J_k$ that have a nonempty intersection with $I_{j_0}$;
we assume, without loss of generality, that for $0\leq k\leq K'-1$, $J_k\cap [a,b] \neq \varnothing$, and write
$$
[a,b]=\cup_{k=0}^{K'-1}J_k=\cup_{k=0}^{K'-1}[t_k,t_{k+1}], \quad t_0=a, \quad t_{K'}=b.
$$
We are going to estimate $u$ on each $J_k\times\R^n$.

First we estimate $u$ on $J_0\times\R^n$. Noting that $w(t_0)=w(a)=u(a)$, by Strichartz, \eqref{nonlin 1},
\eqref{nonlin 2}, \eqref{ass v 2}, and \eqref{lin evol small}, we get
\begin{align*}
\|u\|_{\dot Y^1(J_0)}
&\le \|e^{i(t-a)\Delta}u(a)\|_{\dot Y^1(J_0)}+C(E,M)\eps^{\frac{n+1}{2(2\theta+1)}}\|u\|_{\dot Y^1(J_0)}+C\|u\|_{\dot Y^1(J_0)}^{\frac{n+2}{n-2}}\\
&\le 2\eta+C(E,M)\eps^{\frac{n+1}{2(2\theta+1)}}\|u\|_{\dot
Y^1(J_0)}+C\|u\|_{\dot Y^1(J_0)}^{\frac{n+2}{n-2}}.
\end{align*}
Choosing  $\eta$ and $\eps=\eps(E,M)$ sufficiently small, a standard continuity argument yields
$$
\|u\|_{\dot Y^1(J_0)}\le 4\eta.
$$
In order to apply Lemma~\ref{h1 stability} on the interval $J_0$, we are left to check that the error term
$e=|u|^{p_1}u$ is small. As by \eqref{nonlin 1},
\begin{align}\label{error 2}
\|\nabla e\|_{\dot N^0(J_0\times \R^n)}
\le C(E,M)\eps^{\frac{n+1}{2(2\theta+1)}}\|u\|_{\dot Y^1(J_0)} \le C(E,M)\eps^{\frac{n+1}{2(2\theta+1)}}\eta,
\end{align}
choosing $\eps$ sufficiently small (depending only on the energy and the mass), we may apply Lemma~\ref{h1 stability}
to obtain
$$
\|u-w\|_{\dot S^1(J_0\times \R^n)}\le C(E)\eps^c,
$$
where $c={c(n, \theta)}$ is a small constant.  By the triangle inequality and \eqref{v(a) S bdd}, this implies
$$
\|u\|_{\dot S^1(J_0\times \R^n)}\le C(E),
$$
while by Strichartz, it implies
\begin{align}
\|u(t_1)-w(t_1)\|_{\ho_x}&\le C(E)\eps^c, \notag \\
\|e^{i(t-t_1)\Delta}(u(t_1)-w(t_1))\|_{\dot Y^1(J_1)}&\le C(E)\eps^c.\label{diff lin evol 2}
\end{align}
Now, we can use these last two estimates to control $u$ on $J_1$. By Strichartz, \eqref{nonlin 1}, \eqref{nonlin 2},
\eqref{lin evol small}, and \eqref{diff lin evol 2}, we get
\begin{align*}
\|u\|_{\dot Y^1(J_1)}
&\le \|\propagateo u(t_1)\|_{\dot Y^1(J_1)}+C(E,M)\eps^{\frac{n+1}{2(2\theta+1)}}
\|u\|_{\dot Y^1(J_1)}+C\|u\|_{\dot Y^1(J_1)}^{\frac{n+2}{n-2}}\\
&\le \|\propagateo w(t_1)\|_{\dot Y^1(J_1)}+\|\propagateo (w(t_1)-u(t_1))\|_{\dot Y^1(J_1)}\\
&\quad + C(E,M)\eps^{\frac{n+1}{2(2\theta+1)}}\|u\|_{\dot Y^1(J_1)}+C\|u\|_{\dot Y^1(J_1)}^{\frac{n+2}{n-2}}\\
&\le 2\eta+C(E)\eps^c+C(E,M)\eps^{\frac{n+1}{2(2\theta+1)}}\|u\|_{\dot Y^1(J_1)}^{\frac {n+2}{n-2}}+C\|u\|_{\dot Y^1(J_1)}^{\frac{n+2}{n-2}}.
\end{align*}
Another continuity argument yields
$$
\|u\|_{\dot Y^1(J_1)}\le 4 \eta,
$$
provided $\eta$ and $\eps=\eps(E,M)$ are chosen sufficiently small. Thus \eqref{error 2} holds with $J_0$ replaced by $J_1$.
Applying again Lemma~\ref{h1 stability} on $I:=J_1$, we obtain
$$
\|u-w\|_{\dot S^1(J_1\times\R^n)}\le C(E)\eps_2^{c^2}
$$
and hence also
$$
\|u\|_{\dot S^1(J_1\times\R^n)}\le C(E).
$$
Choosing $\eps$ sufficiently small depending on $E$ and $M$, for any $0\leq k\leq K'-1$ we obtain, by induction,
$$
\|u\|_{\dot Y^1(J_k)}\le 4\eta
$$
and
$$
\|u\|_{\dot S^1(J_k\times\R^n)}\le C(E).
$$
Adding all these estimates together, we get
\begin{align*}
\|u\|_{\dot Y^1([a,b])}&\lesssim \sum_{k=0}^{K'-1}\|u\|_{\dot Y^1(J_k)}\lesssim 4K'\eta \le C(E)\\
\|u\|_{\dot S^1([a,b])}&\lesssim \sum_{k=0}^{K'-1}\|u\|_{\dot S^1(J_k)}\lesssim K' C(E)\le C(E).
\end{align*}
Thus, as the interval $I_{j_0}=[a,b]$ was arbitrarily chosen, we get
\begin{align}
\|u\|_{\dot Y^1(\R)}&\lesssim \sum_{j=0}^{J-1}\|u\|_{\dot Y^1(I_j)}\le J C(E)\le C(E,M)\notag\\
\|u\|_{\dot S^1(\R\times\R^n)}&\lesssim \sum_{j=0}^{J-1}\|u\|_{\dot S^1(I_j\times\R^n)}\le J C(E)\le C(E,M).\label{u s1 5.4}
\end{align}
Moreover, by Strichartz and Lemma~\ref{m interpolation lemma} we have
\begin{align*}
\|u\|_{\dot S^0(\R\times\R^n)}
&\lesssim \|u_0\|_{L_x^2} + C(E,M)\|u\|_{Z(\R)}^{\frac{n+1}{2(2\theta+1)}}\|u\|_{\dot S^0(\R\times\R^n)}
  +\|u\|_{\dot Y^1(\R)}^{\frac{4}{n-2}}\|u\|_{\dot S^0(\R\times\R^n)}.
\end{align*}
So, subdividing $\R$ into subintervals where both the $Z$-norm and the $\dot Y^1$-norm are sufficiently small depending
on $E$ and $M$, we obtain $\dot S^0$-bounds on $u$ on each of these subintervals that depend only on $E$ and $M$.
Adding these bounds, we get
\begin{align}\label{u s0 5.4}
\|u\|_{\dot S^0(\R\times\R^n)}\leq C(E,M).
\end{align}
Putting together \eqref{u s1 5.4} and \eqref{u s0 5.4}, we obtain
$$
\|u\|_{S^1(\R\times\R^n)}\leq C(E,M).
$$

\subsection{Global bounds for $p_1=\frac 4n$, $p_2=\frac 4{n-2}$, and $\lambda_1, \lambda_2>0$}

In this subsection we prove spacetime bounds that imply scattering for solutions to \eqref{equation} in the case when
both nonlinearities are defocusing and $p_1=\frac 4n$, $p_2=\frac 4{n-2}$.   Without loss of generality, we may assume
$\lambda_1=\lambda_2=1$.  In this case we cannot apply Lemma~\ref{m interpolation lemma}; this already suggests that we cannot
treat the $L_x^2$-critical nonlinearity as a perturbation to the energy-critical NLS, nor can we treat the
energy-critical nonlinearity as a perturbation to the $L_x^2$-critical NLS.  However, we can successfully compare the
low frequencies of the solution to the $L_x^2$-critical problem and the high frequencies to the $\dot H^1_x$-critical
problem.  The result is hence conditional upon a satisfactory theory for the $L_x^2$-critical NLS; specifically,
we obtain a good spacetime bound on $u$ under the Assumption~\ref{assumption}.

We will need a few small parameters for our argument.  More precisely, we will need
$$
0< \eta_3 \ll \eta_2 \ll \eta_1 \ll 1,
$$
where each $\eta_j$ is allowed to depend on the energy and the mass of the initial data, as well as on any of the larger
$\eta$'s. We will choose $\eta_j$ small enough such that, in particular, it will be smaller than any constant depending
on the larger $\eta$'s.

By Theorem~\ref{global wellposedness}, it follows that under our hypotheses, \eqref{equation} admits a unique global
solution $u$.  Moreover, by Proposition~\ref{Morawetz control proposition} and conservation of energy and mass,
we have
$$
\|u\|_{Z(\R)}\leq C(E,M).
$$
We split $\R$ into $K=K(E,M,\eta_3)$ subintervals $J_k$ such that on each slab $J_k\times\R^n$ we have
\begin{align}\label{scat4 u small mor}
\|u\|_{Z(J_k)}\sim \eta_3.
\end{align}

We will show that on every slab $J_k\times\R^n$ the solution $u$ obeys good Strichartz bounds.  Fix, therefore,
$J_{k_0}=[a,b]$. For every $t\in J_{k_0}$, we split $u(t)=u_{lo}(t)+u_{hi}(t)$, where $u_{lo}(t):=P_{<\eta_2^{-1}}u(t)$
and $u_{hi}(t):=P_{\geq\eta_2^{-1}}u(t)$.

On the slab $J_{k_0}\times\R^n$, we compare $u_{lo}(t)$ to the following $L^2_x$-critical problem,
\begin{equation}\label{v equation}
\begin{cases}
iv_t +\Delta v = |v|^{\frac 4n} v\\
v(a) = u_{lo}(a),
\end{cases}
\end{equation}
which is globally wellposed in $H^1_x$ and moreover, by Assumption~\ref{assumption},
$$
\|v\|_{V(\R)}\le C(\|u_{lo}(a)\|_{L_x^2})\leq C(M).
$$
By Lemma~\ref{v implies s1}, this implies
\begin{align}
\|v\|_{\dot S^0 (\R\times\R^n)} &\leq C(M) \label{scat4 v in s0}\\
\|v\|_{\dot S^1 (\R\times\R^n)} &\leq C(E,M) \label{scat4 v in s1}.
\end{align}
We divide $J_{k_0}=[a,b]$ into $J=J(M,\eta_1)$ subintervals $I_j=[t_{j-1},t_j]$ with $t_0=a$ and $t_J=b$, such that
\begin{align}\label{scat4 v small}
\|v\|_{V(I_j)}\sim \eta_1.
\end{align}
By induction, we will establish that for each $j=1,\cdots, J$, we have
\begin{equation}\label{jstep}
P(j):
\begin{cases}
\|u_{lo}-v\|_{\dot S^0(I_1\cup\cdots \cup I_j)}\le \eta_2^{1-2\delta}\\
\|u_{hi}\|_{\dot S^1(I_l)}\le L(E),\ \  \text{for every}\ \  1\leq l\leq j\\
\|u\|_{S^1(I_1\cup\cdots\cup I_j)}\le C(\eta_1,\eta_2),
\end{cases}
\end{equation}
where $\delta>0$ is a small constant to be specified later, and $L(E)$ is a certain large quantity to be chosen later
that depends only on $E$ (but not on the $\eta_j$).  As the method of checking that \eqref{jstep} holds for
$j=1$ is similar to that of the induction step, i.e., showing that $P(j)$ implies $P(j+1)$, we will only prove the latter.

Assume therefore that \eqref{jstep} is true for some $1\leq j<J$.  Then, we will show that
\begin{equation}\label{j+1step}
\begin{cases}
\|u_{lo}-v\|_{\dot S^0(I_1\cup\cdots \cup I_{j+1})}\le \eta_2^{1-2\delta}\\
\|u_{hi}\|_{\dot S^1(I_l)}\le L(E),\ \  \text{for every}\ \  1\leq l\leq j+1\\
\|u\|_{S^1(I_1\cup\cdots\cup I_{j+1})}\le C(\eta_1,\eta_2).
\end{cases}
\end{equation}
In order to prove \eqref{j+1step}, we use a bootstrap argument. Let $\Omega_1$ be the set of all times $T\in I_{j+1}$
such that
\begin{align}
\|u_{lo}-v\|_{\dot S^0(I_1\cup\cdots\cup I_j\cup[t_j,T])}& \le \eta_2^{1-2\delta}\label{omega1a}\\
\|u_{hi}\|_{\dot S^1([t_j,T])}& \le L(E)\label{omega1b}\\
\|u\|_{S^1(I_1\cup\cdots\cup I_j\cup[t_j,T])}& \le C(\eta_1,\eta_2).\label{omega1c}
\end{align}
We need to show that $\Omega_1=I_{j+1}$. First, we notice that $\Omega_1$ is a nonempty\footnote{When proving $P(1)$,
this follows immediately provided $L(E)$ is sufficiently large depending on the energy of the initial data and
$C(\eta_1,\eta_2)$ is sufficiently large depending on the energy and the mass of the initial data.} (as $t_j\in \Omega_1$
by the inductive hypothesis) closed (by Fatou) set.  In order to conclude $\Omega_1=I_{j+1}$, we just need to show that
$\Omega_1$ is an open set as well.  Let therefore $\Omega_2$ be the set of all times $T\in I_{j+1}$ such that
\begin{align}
\|u_{lo}-v\|_{\dot S^0(I_1\cup\cdots\cup I_j\cup[t_j,T])}& \le 2\eta_2^{1-2\delta}\label{omega2a}\\
\|u_{hi}\|_{\dot S^1([t_j,T])} &\le 2L(E)\label{omega2b}\\
\|u\|_{S^1(I_1\cup\cdots\cup I_j\cup[t_j,T])} &\le 2C(\eta_1,\eta_2).\label{omega2c}
\end{align}
We will prove that $\Omega_2\subset\Omega_1$, which will conclude the argument.

\begin{lemma}\label{lemma size tkt}
Let $T\in\Omega_2$.  Then, the low frequencies of the solution satisfy
\begin{align}
\|u_{lo}\|_{V(I)}&\lesssim \eta_1, \ \ \text{where}\ I\in\{I_l, \ 1\leq l\leq j\}\cup\{[t_j,T]\},\label{1}\\
\|u_{lo}\|_{\sntnt}&\leq C(M),\label{2}\\
\|\ulo\|_{W([t_j,T])}&\lesssim \eta_2,\label{3}\\
\|u_{lo}\|_{\dot{S}^1(\ir)}&\lesssim E, \ \ \text{where}\ I\in\{I_l, \ 1\leq l\leq j\}\cup\{[t_j,T]\},\label{4}\\
\|u_{lo}\|_{\sftnt}&\leq C(\eta_1)E, \label{5}
\end{align}
while the high frequencies satisfy
\begin{align}
\|\uhi\|_{\dot{S}^0(\ir)}&\lesssim \eta_2L(E),\ \ \text{where}\ I\in\{I_l, \ 1\leq l\leq j\}\cup\{[t_j,T]\},\label{6}\\
\|\uhi\|_{\sntnt}&\leq \eta_2C(\eta_1)L(E),\label{7}\\
\|\uhi\|_{\sftnt}&\leq C(\eta_1)L(E).\label{8}
\end{align}
\end{lemma}

\begin{proof} Using the triangle inequality, Bernstein, \eqref{scat4 v in s0}, \eqref{scat4 v small}, \eqref{jstep}
as well as \eqref{omega2a} and \eqref{omega2b}, we easily check \eqref{1}, \eqref{2}, and \eqref{6}:
\begin{align*}
\|u_{lo}\|_{V(I)}
   &\le \|v-u_{lo}\|_{V(I)}+\|v\|_{V(I)}\\
   &\lesssim \|v-u_{lo}\|_{\dot{S}^0(\ir)}+\|v\|_{V(I)}\\
   &\lesssim \eta_2^{1-2\delta}+\eta_1\lesssim \eta_1,\\
\|u_{lo}\|_{\sntnt}
   &\le \|v- u_{lo}\|_{\sntnt}+\|v\|_{\sntnt}\\
   &\lesssim \eta_2^{1-2\delta}+C(M)\leq C(M),\\
\|u_{hi}\|_{\dot{S}^0(\ir)}
   &\lesssim \eta_2\|u_{hi}\|_{\dot{S}^1(\ir)}\lesssim \eta_2L(E),
\end{align*}
where $I\in\{I_l, \ 1\leq l\leq j\}\cup\{[t_j,T]\}$.  Moreover, as $J=O(\eta_1^{-C})$, we get \eqref{8}:
\begin{align*}
\|\uhi\|_{\sftnt}
   &\lesssim \sum_{l=1}^j\|\uhi\|_{\dot S^1(I_l\times\R^n)}+\|\uhi\|_{\dot S^1([t_j,T]\times\R^n)} \le C(\eta_1)L(E),
\end{align*}
which, by Bernstein, implies \eqref{7}:
\begin{align*}
\|\uhi\|_{\sntnt}&\lesssim \eta_2\|\uhi\|_{\sftnt}\le \eta_2C(\eta_1)L(E).
\end{align*}

Hence, we are left to prove \eqref{3} through \eqref{5}.  Of course, \eqref{5} follows from \eqref{4} and the fact that
$J=O(\eta_1^{-C})$.  Let therefore $I\in\{I_l, \ 1\leq l\leq j\}\cup\{[t_j,T]\}$.  On the slab $\ir$, $u_{lo}$ satisfies
the equation
\begin{equation*}
u_{lo}(t)=e^{i(t-t_l)\Delta}u_{lo}(t_l)-i\int_{t_l}^t e^{i(t-s)\Delta} P_{lo}(|u|^{\frac 4n}u+|u|^{\nmt}u)(s)ds,
\end{equation*}
where $0\leq l\leq j$.  By Strichartz, this implies
\begin{align}\label{ulo s1}
\|u_{lo}\|_{\dot S^1(\ir)}
&\lesssim \|\ulo(t_l)\|_{\dot H^1_x}+\|\pl(|u|^{\frac 4n}u)\|_{\dot N^1(\ir)}
+\|\pl(|u|^{\nmt}u)\|_{\dot N^1(\ir)}.
\end{align}
By Bernstein, \eqref{scat4 u small mor}, \eqref{omega2c}, and Lemma~\ref{h1 morawetz control}, choosing $\theta>0$
sufficiently small, we get
\begin{align*}
\|P_{lo}(|u|^{\nmt} u)\|_{\dot N^1(\ir)}
&\lesssim \eta_2^{-1}\||u|^{\nmt}u\|_{\dot N^0(\ir)}\\
&\lesssim \eta_2^{-1}\|u\|_{Z(I)}^{\theta}\|u\|^{\frac{n+2}{n-2}-\theta}_{S^1(\ir)}\\
&\lesssim \eta_2^{-1}\eta_3^{\theta}C(\eta_1,\eta_2) \leq \eta_2,
\end{align*}
provided $\eta_3$ is chosen sufficiently small depending on $\eta_1$ and $\eta_2$.

On the other hand, writing $u=u_{lo}+u_{hi}$ and using \eqref{jstep}, \eqref{omega2b}, \eqref{1}, and \eqref{6},
we bound the $L^2_x$-critical term as follows:
\begin{align*}
\|\pl(|u|^{\frac 4n}&u)\|_{\dot N^1(\ir)}\\
&\lesssim \||\nabla \ulo||\ulo|^{\frac4n}\|_{\dot N^0(\ir)}+ \||\nabla \ulo||\uhi|^{\frac4n}\|_{\dot N^0(\ir)}\\
&\quad+ \||\nabla \uhi||\ulo|^{\frac4n}\|_{\dot N^0(\ir)} +\||\nabla \uhi||\uhi|^{\frac4n}\|_{\dot N^0(\ir)}\\
&\lesssim\|\ulo\|_{\dot S^1(\ir)}\|\ulo\|_{V(I)}^{\frac 4n}+\|\ulo\|_{\dot S^1(\ir)}\|\uhi\|_{\dot S^0(\ir)}^{\frac 4n}\\
&\quad +\|\uhi\|_{\dot S^1(\ir)}\|\ulo\|_{V(I)}^{\frac 4n}+\|\uhi\|_{\dot S^1(\ir)}\|\uhi\|_{\dot S^0(\ir)}^{\frac 4n}\\
&\lesssim \|\ulo\|_{\dot S^1(\ir)}(\eta_1+\eta_2L(E))^{\frac 4n}+\eta_1^{\frac 4n}L(E)+\eta_2^{\frac 4n}L(E)^{1+\frac 4n}.
\end{align*}
Therefore, putting everything together, \eqref{ulo s1} becomes
\begin{align*}
\|\ulo\|_{\dot S^1(\ir)}
&\lesssim E+\|\ulo\|_{\dot S^1(\ir)}(\eta_1+2L(E)\eta_2)^{\frac 4n}+\eta_1^{\frac 4n}L(E)+\eta_2^{\frac 4n}L(E)^{1+\frac 4n}.
\end{align*}
Taking $\eta_1$ and $\eta_2$ sufficiently small depending on $E$, this implies
\begin{align*}
\|\ulo\|_{\dot S^1(\ir)}&\lesssim E,
\end{align*}
which settles \eqref{4}.

We turn next to \eqref{3} and write $u_{lo}:=P_{\leq \eta_2}u_{lo} + P_{\eta_2<\cdot<\eta_2^{-1}}u_{lo}$.  We will first verify
\eqref{3} for the medium frequencies of $u$.  As the geometry of the Strichartz trapezoid in dimension $n\geq 5$
is quite different from the geometry of the Strichartz trapezoid in dimensions $n=4$ and $n=3$, we will treat these
cases separately.

In dimension $n\geq 5$, the space $L_{t,x}^{\frac{2(n+2)}{n-2}}$ lies between the two $\dot S^1$-spaces,
$L_t^{n+1}L_x^{\frac{2n(n+1)}{n^2-n-6}}$ and $L_t^2L_x^{\frac{2n}{n-4}}$.  By interpolation and Sobolev embedding,
we get
\begin{align*}
\|P_{\eta_2<\cdot<\eta_2^{-1}}&u_{lo}\|_{W([t_j,T])}\\
& \lesssim \|P_{\eta_2<\cdot<\eta_2^{-1}}u_{lo}\|_{L_t^{n+1}L_x^{\frac{2n(n+1)}{n^2-n-6}}([t_j,T]\times\R^n)}^{c} \|P_{\eta_2<\cdot<\eta_2^{-1}}u_{lo}\|_{L_t^2L_x^{\frac{2n}{n-4}}([t_j,T]\times\R^n)}^{1-c}\\
&\lesssim \||\nabla|^{\frac 3{n+1}}P_{\eta_2<\cdot<\eta_2^{-1}}u_{lo}\|_{Z([t_j,T])}^{c} \|u_{lo}\|_{\dot S^1([t_j,T]\times\R^n)}^{1-c},
\end{align*}
where $c=\frac{4(n+1)}{(n-1)(n+2)}$.  As by Bernstein and \eqref{scat4 u small mor},
\begin{align*}
\bigl\||\nabla|^{\frac 3{n+1}}P_{\eta_2<\cdot<\eta_2^{-1}}u_{lo}\bigr\|_{Z([t_j,T])}
\lesssim \eta_2^{-\frac 3{n+1}}\|u_{lo}\|_{Z([t_j,T])}
\lesssim \eta_2^{-\frac 3{n+1}}\eta_3\le \eta_3^{\frac12},
\end{align*}
\eqref{4} implies
$$
\|P_{\eta_2<\cdot<\eta_2^{-1}}u_{lo}\|_{W([t_j,T])}
\lesssim \eta_3^{c/2}\|u_{lo}\|_{\dot S^1([t_j,T]\times\R^n)}^{1-c}
\lesssim \eta_3^{c/2}E^{1-c}\le \eta_2.
$$

In dimension $n=4$, the space $L_{t,x}^{\frac{2(n+2)}{n-2}}=L_{t,x}^6$ lies between the two $\dot S^1$ spaces
$L_t^{n+1}L_x^{\frac{2n(n+1)}{n^2-n-6}}=L_t^5L_x^{\frac{20}{3}}$ and $L_t^\infty L_x^{\frac{2n}{n-2}}=L_t^\infty L_x^4$.
By interpolation, Sobolev embedding, Bernstein, the conservation of energy, and \eqref{scat4 u small mor}, we obtain
\begin{align*}
\|P_{\eta_2<\cdot<\eta_2^{-1}}u_{lo}\|_{W([t_j,T])}
&\lesssim \|P_{\eta_2<\cdot<\eta_2^{-1}}u_{lo}\|_{L_t^5L_x^{\frac{20}{3}}}^{\frac{5}{6}}\|P_{\eta_2<\cdot<\eta_2^{-1}}u_{lo}\|_{L_t^\infty L_x^4}^{\frac{1}{6}}\\
&\lesssim \||\nabla|^{\frac 3{5}}P_{\eta_2<\cdot<\eta_2^{-1}}u_{lo}\|_{Z([t_j,T])}^{\frac{5}{6}} E^{\frac{1}{6}}\\
&\lesssim \bigl(\eta_2^{-\frac 3{5}}\|P_{\eta_2<\cdot<\eta_2^{-1}}u_{lo}\|_{Z([t_j,T])}\bigr)^{\frac{5}{6}}E^{\frac{1}{6}}\\
&\lesssim (\eta_2^{-\frac 3{5}}\eta_3)^{\frac{5}{6}} E^{\frac{1}{6}}\\
&\leq \eta_2.
\end{align*}

In dimension $n=3$, the space $L_{t,x}^{\frac{2(n+2)}{n-2}}=L_{t,x}^{10}$ lies between the two spaces
$L_t^{n+1}L_x^{\frac{2n(n+1)}{n^2-n-6}}=L_t^4L_x^\infty$ and $L_t^\infty L_x^{\frac{2n}{n-2}}=L_t^\infty L_x^6$.
However, because Sobolev embedding fails at the endpoint, we are forced to take $\eps$ more derivatives in order to bound
the $L_t^4L_x^\infty$-norm in terms of the $L_{t,x}^4$-norm.  More precisely, by interpolation, embedding, Bernstein,
conservation of energy, and \eqref{scat4 u small mor}, we get
\begin{align*}
\|P_{\eta_2<\cdot<\eta_2^{-1}}u_{lo}&\|_{W([t_j,T])}\\
&\lesssim \|P_{\eta_2<\cdot<\eta_2^{-1}}u_{lo}\|_{L_t^4L_x^\infty([t_j,T])\times\R^n}^{\frac{2}{5}} \|P_{\eta_2<\cdot<\eta_2^{-1}}u_{lo}\|_{L_t^{\infty}L_t^6([t_j,T]\times\R^n)}^{\frac{3}{5}}\\
&\lesssim \|(1+|\nabla|)^{\frac 3{4}+\epsilon}P_{\eta_2<\cdot<\eta_2^{-1}}u_{lo}\|_{Z([t_j,T])}^{\frac{2}{5}}E^{\frac{3}{5}}\\
&\lesssim \bigl(\eta_2^{-\frac{3}{4}}\|P_{\eta_2<\cdot<\eta_2^{-1}}u_{lo}\|_{Z([t_j,T])}\bigr)^{\frac{2}{5}}E^{\frac{3}{5}}\\
&\lesssim (\eta_2^{-\frac{3}{4}}\eta_3)^{\frac{2}{5}}E^{\frac{3}{5}}\\
&\leq \eta_2.
\end{align*}

Hence, in all dimensions $n\geq 3$, we have
\begin{align}\label{scat4 med freq}
\|u_{\eta_2<\cdot<\eta_2^{-1}}\|_{W([t_j,T])} \leq \eta_2.
\end{align}

We turn now to estimating the very low frequencies of $u$.  By Sobolev embedding, Bernstein, interpolation,
\eqref{2}, and the conservation of mass, we get
\begin{align}\label{scat4 vlo freq}
\|P_{\leq\eta_2}u_{lo}\|_{W([t_j,T])}
&\lesssim \|\nabla P_{\leq\eta_2}u_{lo}\|_{L_t^{\frac{2(n+2)}{n-2}}L_x^{\frac{2n(n+2)}{n^2+4}}([t_j,T]\times\R^n)}\notag\\
&\lesssim \eta_2 \|u_{lo}\|_{L_t^{\frac{2(n+2)}{n-2}}L_x^{\frac{2n(n+2)}{n^2+4}}([t_j,T]\times\R^n)}\notag\\
&\lesssim \eta_2 \|u_{lo}\|_{V([t_j,T])}^{\frac{n-2}{n}}\|u_{lo}\|_{L_t^\infty L_x^2([t_j,T]\times\R^n)}^{\frac{2}{n}}\notag\\
&\lesssim \eta_2\eta_1^{\frac{n-2}{n}}M^{\frac{2}{n}}
 \leq \eta_2,
\end{align}
provided $\eta_1$ is chosen sufficiently small depending on $M$.

Therefore, by the triangle inequality, \eqref{scat4 med freq} and \eqref{scat4 vlo freq} imply \eqref{3}.
\end{proof}

We are now ready to resume our argument and prove that $\Omega_2\subset\Omega_1$.  Let therefore $T\in\Omega_2$.
We will first show \eqref{omega1a}. The idea is to compare $u_{lo}$ to $v$ via the perturbation result
Lemma~\ref{l2 stability}.

The low frequency-part of the solution, $\ulo$, solves the following initial value problem on the slab $[t_0,T]\times\R^n$
\begin{equation*}
\begin{cases}
(i\partial_t+\Delta)\ulo =|\ulo|^{\frac 4n}\ulo+\pl(|u|^{\frac 4n}u-|\ulo|^{\frac 4n}\ulo)-\ph(|\ulo|^{\frac 4n}\ulo)+\pl(|u|^{\nmt}u)\\
\ulo(t_0)=\ulo(a).
\end{cases}
\end{equation*}
As by \eqref{2},
$$
\|u_{lo}\|_{\sntnt}\leq C(M),
$$
and $ v(t_0)=\ulo(t_0)$, in order to apply Lemma~\ref{l2 stability}, we need only verify that the error term
$$
e_1=\pl(|u|^{\frac 4n}u-|\ulo|^{\frac 4n}\ulo)-\ph(|\ulo|^{\frac 4n}\ulo)+\pl(|u|^{\nmt}u)
$$
is small in $\nntnt$.

By H\"older, \eqref{2} and \eqref{7}, we estimate
\begin{align*}
\|\pl(|u|^{\frac 4n}u-|\ulo|^{\frac 4n}\ulo&)\|_{\nntnt}\\
&\lesssim \||\uhi|^{1+\frac 4n}\|_{\nntnt}+\||\uhi||\ulo|^{\frac 4n}\|_{\nntnt}\\
&\lesssim\|\uhi\|_{\sntnt}^{1+\frac 4n}+\|\uhi\|_{\sntnt}\|\ulo\|^{\frac 4n}_{\sntnt}\\
&\lesssim (\eta_2C(\eta_1)L(E))^{1+\frac 4n}+\eta_2C(\eta_1)L(E)C(M)\\
&\leq \eta_2^{1-\delta},
\end{align*}
provided $\eta_2$ is sufficiently small depending on $E$, $M$, and $\eta_1$; here, $\delta>0$ is a small parameter.

Moreover, by Bernstein, \eqref{2}, and \eqref{5}, we get
\begin{align*}
\|\ph(|\ulo|^{\frac 4n}\ulo)\|_{\nntnt}
&\lesssim \eta_2 \|\nabla \ph(|\ulo|^{\frac 4n}\ulo)\|_{\nntnt}\\
&\lesssim \eta_2\| \ulo\|_{\sftnt}\|\ulo\|_{\sntnt}^{\frac 4n}\\
&\lesssim \eta_2 C(\eta_1)EC(M)\le \eta_2^{1-\delta}.
\end{align*}

By Lemma~\ref{h1 morawetz control} (with $\theta$ sufficiently small), \eqref{scat4 u small mor}, and \eqref{omega2c},
we get
\begin{align*}
\|\pl(|u|^{\nmt}u)\|_{\nntnt}
&\lesssim \|u\|^{\theta}_{Z[t_0,T]}\|u\|^{\frac{n+2}{n-2}-\theta}_{S^1([t_0,T]\times\R^n)}
\lesssim \eta_3^{\theta}(C(\eta_1,\eta_2))^{\frac{n+2}{n-2}-\theta}\\
&\leq \eta_2,
\end{align*}
provided $\eta_3$ is chosen sufficiently small depending on $\eta_1$ and $\eta_2$.

Therefore,
$$
\|e_1\|_{\nntnt}\le 3\eta_2^{1-\delta},
$$
and hence, taking $\eta_2$ sufficiently small depending on $M$, we can apply Lemma~\ref{l2 stability} to obtain
$$
\|\ulo-v\|_{\sntnt}\le C(M)\eta_2^{1-\delta}\leq \eta_2^{1-2\delta}.
$$

Thus, \eqref{omega1a} holds.  We turn now to \eqref{omega1b}; the idea is to compare the high frequency-part
of the solution $u$ to the energy-critical NLS.  Consider therefore the initial value problem
\begin{equation}\label{equation u2}
\begin{cases}
iw_t+\Delta w=|w|^{\nmt} w\\
w(t_j)=\uhi(t_j).
\end{cases}
\end{equation}
Then, by \cite{ckstt:gwp, RV, Monica:thesis}, \eqref{equation u2} is globally wellposed and furthermore,
\begin{align}\label{scat4 w s1}
\|w\|_{\so(\R\times\R^n)}\le C(E).
\end{align}
By Lemma~\ref{w implies s1}, \eqref{6}, and Strichartz, this also implies
\begin{align}\label{scat4 w s0}
\|w\|_{\sz(\R\times\R^n)}\le C(E)\|\uhi(t_j)\|_{L_x^2} \lesssim \eta_2C(E)L(E).
\end{align}
On the other hand, the high frequency-portions of $u$ satisfy
\begin{equation*}
\begin{cases}
(i\partial_t+\Delta)\uhi & =|\uhi|^{\nmt}\uhi+\ph(|u|^{\nmt}u-|\uhi|^{\nmt}\uhi)-\pl(|u|^{\nmt}u)\\
&\quad+\ph(|u|^{\frac 4n}u)\\
\uhi(t_j)&=\uhi(t_j).
\end{cases}
\end{equation*}
By the bootstrap assumption \eqref{omega2b}, we have
$$
 \|\uhi\|_{\dot S^1([t_j,T]\times\R^n)}\le 2L(E).
$$
Therefore, in order to apply Lemma~\ref{h1 stability}, we need only to check that the error term
$$
e_2=\ph(|u|^{\nmt}u-|\uhi|^{\nmt} \uhi)-\pl(|u|^{\nmt} u)+\ph(|u|^{\frac 4n} u)
$$
is small in $\dot N^1([t_j,T]\times\R^n)$.

By the triangle inequality, we easily see that
\begin{align}
\|\ph(|u&|^{\nmt}u-|\uhi|^{\nmt} \uhi)\|_{\dot N^1([t_j,T]\times\R^n)}\notag\\
&\lesssim \bigl\||u|^{\nmt}\nabla\ulo\bigr\|_{\dot N^0([t_j,T]\times\R^n)}
  +\bigl\|\bigl(|u|^{\frac{4}{n-2}}-|u_{hi}|^{\frac{4}{n-2}}\bigr)\nabla\uhi\bigr\|_{\dot N^0([t_j,T]\times\R^n)}.\label{scat4 error 1}
\end{align}

To estimate the first term on the right-hand side of \eqref{scat4 error 1}, we use Remark~\ref{variation}
(for a small constant $\theta>0$), Bernstein, \eqref{scat4 u small mor}, \eqref{omega2c}, \eqref{2}, and \eqref{4} to get
\begin{align*}
\bigl\||u|^{\nmt}\nabla\ulo\bigr\|_{\dot N^0([t_j,T]\times\R^n)}
&\lesssim \|u\|_{Z([t_j,T])}^{\theta}\|u\|_{S^1([t_j,T]\times\R^n)}^{\frac{4}{n-2}-\theta}\|\nabla u_{lo}\|_{S^1([t_j,T]\times\R^n)}\\
&\lesssim \eta_3^{\theta}C(\eta_1,\eta_2)\eta_2^{-1}\|u_{lo}\|_{S^1([t_j,T]\times\R^n)}\\
&\le\eta_2,
\end{align*}
provided $\eta_3$ is chosen sufficiently small depending on $\eta_1$ and $\eta_2$.

Next, we estimate the second term on the right-hand side of \eqref{scat4 error 1}.  When the dimension $3\leq n<6$,
we use H\"older, \eqref{omega2b}, \eqref{3}, and \eqref{4} to obtain
\begin{align*}
\bigl\|\bigl(|u|^{\frac{4}{n-2}}&-|u_{hi}|^{\frac{4}{n-2}}\bigr)\nabla\uhi\bigr\|_{\dot N^0([t_j,T]\times\R^n)}\\
&\lesssim \||u|^{\frac{6-n}{n-2}}\ulo\nabla\uhi\|_{\dot N^0([t_j,T]\times\R^n)}\\
&\lesssim \bigl(\|\uhi\|_{\dot S^1([t_j,T]\times\R^n)}^{\frac{6-n}{n-2}} +\|\ulo\|_{\dot S^1([t_j,T]\times\R^n)}^{\frac{6-n}{n-2}}\bigr)\|\ulo\|_{W[t_j,T]}\|\uhi\|_{\dot S^1([t_j,T]\times\R^n)}\\
&\lesssim (L(E)+E)^{\frac{6-n}{n-2}}\eta_2 L(E) \\
&\le \eta_2^{\frac 34}.
\end{align*}
For $n\ge 6$, we use \eqref{omega2b} and \eqref{3} to estimate
\begin{align*}
\bigl\|\bigl(|u|^{\frac{4}{n-2}}-|u_{hi}|^{\frac{4}{n-2}}\bigr)\nabla\uhi\bigr\|_{\dot N^0([t_j,T]\times\R^n)}
&\lesssim \bigl\||\ulo|^{\nmt}\nabla\uhi\bigr\|_{\dot N^0([t_j,T]\times\R^n)}\\
&\lesssim \|\uhi\|_{\dot S^1([t_j,T]\times\R^n)}\|\ulo\|_{W([t_j,T])}^{\nmt}\\
&\lesssim L(E)\eta_2^{\nmt}\le \eta_2^{\frac 3{n-2}}.
\end{align*}

Collecting these estimates, we get
\begin{align}\label{scat4 error 2}
\|\ph(|u|^{\nmt}u-|\uhi|^{\nmt} \uhi)\|_{\dot N^1([t_j,T]\times\R^n)}
\leq \eta_2 + \eta_2^{\frac{3}{4}}+ \eta_2^{\frac{3}{n-2}}.
\end{align}

To estimate the second term in the expression of the error $e_2$, we use Bernstein, Lemma~\ref{h1 morawetz control}
(for some small constant $\theta>0$), \eqref{scat4 u small mor}, and \eqref{omega2c}:
\begin{align}
\|\pl(|u|^{\nmt}u)\|_{\dot N^1([t_j,T]\times\R^n)}
&\lesssim \eta_2^{-1}\||u|^{\nmt}u\|_{\dot N^0([t_j,T]\times\R^n)}\notag\\
&\lesssim \eta_2^{-1}\|u\|_{Z([t_j,T])}^{\theta}\|u\|_{S^1([t_j,T]\times\R^n)}^{\frac{n+2}{n-2}-\theta}\notag\\
&\lesssim \eta_2^{-1}\eta_3^{\theta} C(\eta_1,\eta_2)\notag\\
&\leq \eta_2.\label{scat4 error 3}
\end{align}

We turn now to the last term in the expression of the error $e_2$. Dropping the projection $\ph$, we estimate
\begin{align*}
\|\nabla\ph (|u|^{\frac 4n}u)&\|_{\dot N^0([t_j,T]\times\R^n)}\\
&\lesssim \bigl\|\nabla\uhi|\ulo|^{\frac 4n}\bigr\|_{\dot N^0([t_j,T]\times\R^n)}+\bigl\|\nabla\uhi|\uhi|^{\frac 4n}\bigr\|_{\dot N^0([t_j,T]\times\R^n)}\\
&\quad+\bigl\|\nabla \ulo|\uhi|^{\frac 4n}\bigr\|_{\dot N^0([t_j,T]\times\R^n)}+\bigl\|\nabla \ulo|\ulo|^{\frac 4n}\bigr\|_{\dot N^0([t_j,T]\times\R^n)}.\\
\end{align*}
Using \eqref{omega2b}, \eqref{1}, \eqref{4}, and \eqref{6}, we estimate the four terms on the right-hand side of the
above inequality as follows:
\begin{align*}
\bigl\|\nabla\uhi|\ulo|^{\frac 4n}\bigr\|_{\dot N^0([t_j,T]\times\R^n)}
&\lesssim \|\nabla \uhi|\ulo|^{\frac 4n}\|_{L^{\frac{2(n+2)}{n+4}}_{t,x}([t_j,T]\times\R^n)}\\
&\lesssim \|\nabla\uhi\|_{V([t_j,T])}\|\ulo\|_{V([t_j,T])}^{\frac 4n}\\
&\lesssim \|\uhi\|_{\dot{S}^1([t_j,T]\times\R^n)}\eta_1^{\frac 4n}\\
&\lesssim L(E)\eta_1^{\frac 4n} \le \eta_1^{\frac 3n},\\
\bigl\|\nabla\uhi|\uhi|^{\frac 4n}\bigr\|_{\dot N^0([t_j,T]\times\R^n)}
&\lesssim \|\uhi\|_{\dot{S}^1([t_j,T]\times\R^n)}\|\uhi\|_{\dot{S}^0([t_j,T]\times\R^n)}^{\frac 4n}\\
&\lesssim \eta_2^{\frac 4n}\|\uhi\|_{\dot{S}^1([t_j,T]\times\R^n)}^{1+\frac 4n}\\
&\lesssim \eta_2^{\frac 4n}L(E)^{1+\frac 4n}\le \eta_1^{\frac 3n},\\
\bigl\|\nabla \ulo|\uhi|^{\frac 4n}\bigr\|_{\dot N^0([t_j,T]\times\R^n)}
&\lesssim \|\ulo\|_{\dot{S}^1([t_j,T]\times\R^n)}\|\uhi\|_{\dot{S}^0([t_j,T]\times\R^n)}^{\frac 4n}\\
&\lesssim E[\eta_2L(E)]^{\frac 4n}\le \eta_1^{\frac 3n},\\
\bigl\|\nabla \ulo|\ulo|^{\frac 4n}\bigr\|_{\dot N^0([t_j,T]\times\R^n)}
&\lesssim \|\ulo\|_{V[t_j,T]}^{\frac 4n}\|\ulo\|_{\dot{S}^1([t_j,T]\times\R^n)}\\
&\lesssim \eta_1^{\frac 4n} E\le \eta_1^{\frac 3n}.
\end{align*}
Adding these estimates, we obtain
\begin{align}\label{scat4 error 4}
\|\nabla\ph (|u|^{\frac 4n}u)&\|_{\dot N^1([t_j,T]\times\R^n)}\lesssim \eta_1^{\frac 3n}.
\end{align}

Collecting \eqref{scat4 error 2} through \eqref{scat4 error 4}, we get
$$
\|e_2\|_{\dot N^1([t_j,T]\times\R^n)}\lesssim \eta_1^{\frac 3n}
$$
which, by Lemma~\ref{h1 stability} implies
$$
\|\uhi-w\|_{\dot S^1([t_j,T]\times\R^n)} \lesssim \eta_1^c
$$
for a small constant $c>0$ depending on the dimension $n$.  By the triangle inequality and \eqref{scat4 w s1},
this implies
$$
\|\uhi\|_{\dot S^1([t_j,T]\times\R^n)}\lesssim C(E)+\eta_1^c\le L(E),
$$
provided $L(E)$ is sufficiently large.

Finally, \eqref{omega1c} follows from \eqref{2}, \eqref{5}, \eqref{7}, and \eqref{8}:
\begin{align*}
\|u\|_{S^1([t_0,T]\times\R^n)}
&\le \|\ulo\|_{S^1([t_0,T]\times\R^n)}+ \|\uhi\|_{S^1([t_0,T]\times\R^n)}\\
&\lesssim C(M) + C(\eta_1)E + \eta_2C(\eta_1)E + C(\eta_1)L(E),
\end{align*}
which is of course bounded by $C(\eta_1,\eta_2)$ provided $\eta_1$ and $\eta_2$ are sufficiently large depending
on $E$ and $M$.

This proves that $\Omega_2\subset \Omega_1$.  Hence, by induction,
$$
\|u\|_{S^1(J_{k_0}\times\R^n)}\leq C(\eta_1,\eta_2).
$$
As $J_{k_0}$ was chosen arbitrarily and the total number of intervals $J_k$ is $K=K(E,M,\eta_3)$, adding these bounds
we obtain
$$
\|u\|_{S^1(\R\times\R^n)}\leq C(\eta_1,\eta_2,\eta_3)=C(E,M).
$$

\subsection{Global bounds for $\frac 4n\le p_1<p_2=\frac 4{n-2}$, $\lambda_2>0$, $\lambda_1<0$, and small mass}

The approach to proving finite global Strichartz bounds for solutions $u$ to \eqref{equation} in this case is similar
to the one used in subsection 5.4.  However, as in this case we do not have a good \emph{a priori} interaction
Morawetz inequality, we will rely instead on the small mass assumption.  As in subsection 5.4, in this case we also
compare \eqref{equation} to the energy-critical problem
$$
\begin{cases}
iw_t+\Delta w=|w|^{\nmt}w\\
w(0)=u_0,
\end{cases}
$$
which by \cite{ckstt:gwp, RV, Monica:thesis} is globally wellposed and moreover,
\begin{align}\label{scat5 w s1}
\|w\|_{\dot S^1(\R^\times\R^n)}\le C(E).
\end{align}
By Lemma~\ref{w implies s1}, \eqref{scat5 w s1} implies
\begin{align}\label{scat5 w s0}
\|w\|_{\dot S^0(\R\times\R^n)}\le C(E)\|u_0\|_{L_x^2}\leq C(E) M^{\frac 12}.
\end{align}

In this subsection we will carry out our analysis in the following spaces:
$$
\dot Y^0(I):=V(I)\cap L_t^{\frac{2(n+2)}{n-2}}L_x^{\frac{2(n+2)}{n^2+4}}(I\times\R^n)
$$
and
$$
\dot Y^1(I):= \{u;\, \nabla u\in \dot{Y}^0\},\ \ \ Y^1(I):=\dot{Y}^0(I) \cap \dot{Y}^1(I).
$$
In this notation, Lemma~\ref{lemma bound by w01} reads
\begin{align}
\bigl\|\nabla^k\bigl(|u|^{p_1}u\bigr)\bigr\|_{\dot N^0(I\times\R^n)}
&\lesssim \|u\|_{\dot Y^0(I)}^{2-\frac{p_1(n-2)}2}\|u\|^{\frac{np_1}2-2}_{\dot Y^1(I)}\|u\|_{\dot Y^k(I)}\label{scat5 est 1}\\
\bigl\|\nabla^k\bigl(|u|^{\nmt}u\bigr)\bigr\|_{\dot N^0(I\times\R^n)}
&\lesssim \|u\|_{\dot Y^1(I)}^{\nmt}\|u\|_{\dot Y^k(I)};\label{scat5 est 2}
\end{align}
here $k=0,1$.

By time reversal symmetry, it suffices to show that $u$ obeys good spacetime bounds on $\R^+\times\R^n$.
Let $\eta>0$ be a small constant to be chosen later and divide $\R^+$ into $J=J(E,\eta)$ subintervals $I_j=[t_j,t_{j+1}]$
such that
\begin{align*}
\|w\|_{\dot Y^1(I_j)}\sim \eta.
\end{align*}
Moreover, choosing $M$ sufficiently small depending on $E$ and $\eta$, by \eqref{scat5 w s0} we may assume
$$
\|w\|_{\dot S^0(\R\times\R^n)}\le \eta.
$$
Therefore,
\begin{align}\label{scat5 w small y1}
\|w\|_{Y^1(I_j)}\sim \eta.
\end{align}

As an immediate consequence of \eqref{scat5 w small y1}, we show that the linear flow of $w$ is small on each slab
$I_j\times\R^n$.  Indeed, by Strichartz and \eqref{scat5 w small y1},
\begin{align}\label{scat5 w linear}
\|e^{i(t-t_j)\Delta}w(t_j)\|_{Y^1(I_j)}
&\le \|w\|_{Y^1(I_j)}+C\|w\|_{Y^1(I_j)}^{\frac{n+2}{n-2}}
 \le \eta+C\eta^{\frac{n+2}{n-2}}\le 2\eta,
\end{align}
provided $\eta$ is sufficiently small.

We first consider the interval $I_0=[t_0,t_1]$. Recalling that $w(t_0)=u(t_0)=u_0$, Strichartz, \eqref{scat5 est 1},
\eqref{scat5 est 2}, \eqref{scat5 w small y1}, and \eqref{scat5 w linear} imply
\begin{align*}
\|u\|_{Y^1(I_0)}\le 2\eta+C\|u\|_{Y^1(I_0)}^{p_1+1}+C\|u\|_{Y^1(I_0)}^{\frac{n+2}{n-2}}.
\end{align*}
By a standard continuity argument, this yields
\begin{align}\label{scat5 u y1}
\|u\|_{Y^1(I_0)}\le 4\eta,
\end{align}
provided $\eta$ is chosen sufficiently small.

Another application of Strichartz together with \eqref{scat5 est 1}, \eqref{scat5 est 2}, and \eqref{scat5 u y1}, yields
\begin{align*}
\|u\|_{\dot Y^0(I_0)}
&\lesssim M^{\frac{1}{2}}+\|u\|_{\dot Y^1(I_0)}^{\frac{np_1}2-2}\|u\|_{\dot Y^0(I_0)}^{3-\frac{p_1(n-2)}2}+\|u\|_{\dot Y^1(I_0)}^{\nmt}\|u\|_{\dot Y^0(I_0)}\\
&\lesssim M^{\frac{1}{2}}+\eta^{\frac{np_1}2-2}\|u\|_{\dot Y^0(I_0)}^{3-\frac{p_1(n-2)}2}+\eta^{\nmt}\|u\|_{\dot Y^0(I_0)}.
\end{align*}
Therefore, choosing $\eta$ sufficiently small and remembering that $M$ is chosen tiny and $3-\frac{p_1(n-2)}2>1$, we get
\begin{align}\label{scat5 u y0}
\|u\|_{\dot Y^0(I_0)}\lesssim M^{\frac{1}{2}}.
\end{align}

In order to apply Lemma~\ref{h1 stability} we need to show that the error, which in this case is the energy-subcritical
perturbation, is small.  By \eqref{scat5 est 1}, \eqref{scat5 u y1}, and \eqref{scat5 u y0}, we find
\begin{align*}
\||u|^{p_1}u\|_{\dot N^1(I_0\times\R^n)}
&\lesssim \|u\|_{\dot Y^0(I_0)}^{2-\frac{p_1(n-2)}2}\|u\|^{\frac{np_1}2-1}_{\dot Y^1(I_0)}
 \lesssim M^{1-\frac{p_1(n-2)}4}\eta^{\frac{np_1}2-1}\le M^{\delta_0},
\end{align*}
for a small constant $\delta_0>0$.  Taking $M$ sufficiently small depending on $E$ and $\eta$, by Lemma~\ref{h1 stability}
we derive
\begin{align*}
\|u-w\|_{\dot S^1(I_0\times\R^n)}\leq M^{c\delta_0},
\end{align*}
for a small constant $c>0$ that depends only on the dimension $n$.  By Strichartz, this implies
\begin{align*}
\|e^{i(t-t_1)\Delta}(u(t_1)-w(t_1))\|_{\dot S^1(I_1\times\R^n)}\lesssim M^{c\delta_0}.
\end{align*}

We turn now to the second interval $I_1=[t_1,t_2]$.  By Strichartz, the triangle inequality, \eqref{scat5 est 1},
\eqref{scat5 est 2}, and \eqref{scat5 w linear}, we get
\begin{align*}
\|u\|_{Y^1(I_1)}
&\le \|e^{i(t-t_1)\Delta}u(t_1)\|_{\dot Y^0(I_1)}+\|e^{i(t-t_1)\Delta}(u(t_1)-w(t_1))\|_{\dot Y^1(I_1)}\\
&\quad +\|e^{i(t-t_1)\Delta}w(t_1)\|_{\dot Y^1(I_1)}+C\|u\|_{Y^1(I_1)}^{p_1+1}+C\|u\|_{Y^1(I_1)}^{\frac{n+2}{n-2}}\\
&\lesssim M^{\frac12}+M^{c\delta_0}+\eta+\|u\|_{Y^1(I_1)}^{p_1+1}+\|u\|_{Y^1(I_1)}^{\frac{n+2}{n-2}}.
\end{align*}
Therefore, choosing $\eta$ sufficiently small, by a standard continuity argument we obtain
$$
\|u\|_{Y^1(I_1)}\le 4\eta.
$$
Moreover, arguing as above, we also get
$$
\|u\|_{\dot Y^0(I_1)}\lesssim M^{\frac{1}{2}},
$$
which allows us to control the energy-subcritical term on the slab $I_1\times\R^n$ in terms of $M$.  For $M$ sufficiently
small, we can apply Lemma~\ref{h1 stability} to obtain
\begin{align*}
\|u-w\|_{\dot S^1(I_0\times\R^n)}\leq M^{c\delta_1},
\end{align*}
for a small constant $0<\delta_1<\delta_0$.

By induction, choosing $M$ smaller at every step (depending only on $E$ and $\eta$), we obtain
$$
\|u\|_{Y^1(I_j)}\le 4\eta.
$$
Summing these estimates over all intervals $I_j$ and recalling that the total number of these intervals is $J=J(E,\eta)$,
we get
$$
\|u\|_{Y^1(R^+)}\lesssim J \eta \le C(E).
$$
By Strichartz, \eqref{scat5 est 1}, and \eqref{scat5 est 2}, this implies
\begin{align*}
\|u\|_{S^1(\R^+\times\R^n)}
&\lesssim \|u_0\|_{H^1_x} + \|u\|_{Y^1(I_0)}^{p_1+1} + \|u\|_{Y^1(I_0)}^{\frac{n+2}{n-2}}\\
&\lesssim M + E + C(E) \leq C(E,M)=C(E).
\end{align*}
By time reversal symmetry, we obtain
$$
\|u\|_{S^1(\R\times\R^n)}\le C(E).
$$

\subsection{Global bounds for $\frac 4n\le p_1<p_2<\nmt$, $\lambda_2>0$, $\lambda_1<0$, and small mass}

In this case, we compare \eqref{equation} to the free Schr\"odinger equation,
$$
i\tilde u_t+\Delta \tilde u=0,\ \ \ \tilde u(0)=u_0.
$$
By Strichartz, the global solution $\tilde u$ obeys the spacetime estimates
\begin{align*}
\|\tilde u\|_{\dot S^1(\ir)}&\lesssim \|u_0\|_{\dot H_x^1}\lesssim E^{\frac12},\\
\|\tilde u\|_{\dot S^0(\ir)}&\lesssim \|u_0\|_{L_x^2} \lesssim M^{\frac12}.
\end{align*}

For a slab $\ir$, let $\dot Y^0(I)$, $\dot Y^1(I)$, and $Y^1(I)$ be the same spaces as in the last subsection.
More precisely,
$$
\dot Y^0(I):=V(I)\cap L_t^{\frac{2(n+2)}{n-2}}L_x^{\frac{2(n+2)}{n^2+4}}(I\times\R^n)
$$
and
$$
\dot Y^1(I):= \{u;\, \nabla u\in \dot{Y}^0\},\ \ \ Y^1(I):=\dot{Y}^0(I) \cap \dot{Y}^1(I).
$$

We divide $\R^+$ into $J=J(E,\eta)$ subintervals $I_j$ such that on each $I_j=[t_j, t_{j+1}]$,
$$
\|\tilde u\|_{\dot Y^1(I_j)}\sim\eta,
$$
where $\eta>0$ is a small parameter.  Choosing $M$ small enough depending on $\eta$, we may assume that on each
slab $I_j\times\R^n$, we have
$$
\|\tilde u\|_{Y^1(I_j)}\sim\eta.
$$

We will prove that $u$ obeys good spacetime bounds on each slab $I_j\times\R^n$.  Consider the first interval,
$I_0=[t_0, t_1]$.  As $u(t_0)=\tilde u(t_0)=u_0$, by Strichartz and Lemma~\ref{lemma bound by w01},
$$
\|u\|_{Y^1(I_0)}
\le \|\tilde u \|_{Y^1(I_0)} + C(\lambda_1,\lambda_2) \sum_{i=1,2} \|u\|_{Y^1(I_0)}^{p_i+1}
\le \eta + C(\lambda_1,\lambda_2)\sum_{i=1,2}\|u\|_{Y^1(I_0)}^{p_i+1}.
$$
A standard continuity argument yields
$$
\|u\|_{Y^1(I_0)}\le 2\eta,
$$
provided $\eta$ is chosen sufficiently small.  Hence, another application of Strichartz and
Lemma~\ref{lemma bound by w01} implies
\begin{align*}
\|u\|_{\dot Y^0(I_0)}
&\lesssim M^{\frac12}+C(\lambda_1,\lambda_2)\sum_{i=1,2}\||u|^{p_i}u\|_{\dot N^0(I_0\times\R^n)}\\
&\lesssim M^{\frac12}+C(\lambda_1,\lambda_2)\sum_{i=1,2}\|u\|_{\dot Y^0(I_0)}^{3-\frac{p_i(n-2)}2}\|u\|^{\frac{np_i}2-2}_{\dot Y^1(I_0)}\\
&\lesssim M^{\frac12}+C(\lambda_1,\lambda_2)\sum_{i=1,2}\|u\|_{\dot Y^0(I_0)}^{3-\frac{p_i(n-2)}2}\eta^{\frac{np_i}2-2}.
\end{align*}
Therefore,
\begin{align}\label{scat6 u y0}
\|u\|_{\dot Y^0(I_0)}\lesssim M^{\frac12},
\end{align}
provided $M$ is chosen sufficiently small.

By Strichartz and \eqref{scat6 u y0}, we can now compute the difference between $u$ and $\tilde u$ on the slab
$I_0\times\R^n$:
\begin{align*}
\|u-\tilde u\|_{\dot S^1(I_0\times\R^n)}
&\lesssim C(\lambda_1,\lambda_2)\sum_{i=1,2}\|u\|_{\dot Y^0(I_0)}^{2-\frac{p_i(n-2)}2}\|u\|_{\dot Y^1(I_0)}^{\frac{np_i}2-1}\\
&\lesssim M^{\delta}\sum_{i=1,2}\eta^{\frac{np_i}2-1} \\
&\le M^{\delta},
\end{align*}
where $\delta>0$ is a small constant, provided $\eta$ is chosen sufficiently small.  By Strichartz, this implies
\begin{align}\label{scat6 diff 1}
\|e^{i(t-t_1)\Delta}(u(t_1)-\tilde u(t_1))\|_{\dot S^1(I_1\times\R^n)}\lesssim M^\delta.
\end{align}

We turn now to the second interval, $I_1=[t_1,t_2]$.  By Strichartz, the triangle inequality,
Lemma~\ref{lemma bound by w01}, and \eqref{scat6 diff 1}, we get
\begin{align*}
\|u\|_{Y^1(I_1)}
&\le \|e^{i(t-t_1)\Delta}u(t_1)\|_{\dot Y^0(I_1)}+\|e^{i(t-t_1)\Delta} \tilde u(t_1)\|_{\dot Y^1(I_1)}\\
&\quad +\|e^{i(t-t_1)\Delta}(u(t_1)-\tilde u(t_1))\|_{\dot Y^1(I_1)}+C(\lambda_1,\lambda_2)\sum_{i=1,2}\|u\|_{Y^1(I_1)}^{p_i+1}\\
&\lesssim M^{\frac12}+\eta+M^\delta+C(\lambda_1,\lambda_2)\sum_{i=1,2}\|u\|_{Y^1(I_1)}^{p_i+1}.
\end{align*}
Another continuity argument yields
\begin{align}\label{scat6 u y1 bis}
\|u\|_{Y^1(I_1)}\le 2\eta,
\end{align}
provided $\eta$ and $M$ are chosen sufficiently small.  By the same arguments as those used to establish
\eqref{scat6 u y0}, this implies
\begin{align}\label{scat6 u y0 bis}
\|u\|_{\dot Y^0(I_1)}\lesssim M^{\frac12}.
\end{align}
Therefore, we can now control the difference between $u$ and $\tilde u$ on $I_1\times\R^n$.  Indeed, by Strichartz,
Lemma~\ref{lemma bound by w01}, \eqref{scat6 diff 1}, \eqref{scat6 u y1 bis}, and \eqref{scat6 u y0 bis},
\begin{align*}
\|u-\tilde u\|_{\dot S^1(I_1\times\R^n)}
&\lesssim \|e^{i(t-t_1)\Delta}(u(t_1)-\tilde u(t_1))\|_{\dot S^1(I_1\times\R^n)}\\
&\quad +\|\ntu\|_{\dot N^1(I_1\times\R^n)}\\
&\lesssim M^\delta + C(\lambda_1,\lambda_2)\sum_{i=1,2}\|u\|_{\dot Y^0(I_1)}^{2-\frac{p_i(n-2)}2}\|u\|_{\dot Y^1(I_1)}^{\frac{np_i}2-1}\\
&\lesssim M^\delta + C(\lambda_1,\lambda_2)\sum_{i=1,2} M^{1-\frac{p_i(n-2)}4}\eta^{\frac{np_i}2-1}\\
&\lesssim M^{\delta},
\end{align*}
provided $\delta>0$ is chosen sufficiently small.  By Strichartz, this immediately implies
$$
\|e^{i(t-t_2)\Delta}(u(t_2)-\tilde u(t_2))\|_{\dot Y^1(I_2)}\lesssim M^\delta
$$
and hence one can repeat the argument on the next slab, that is, $I_2\times\R^n$.  By induction, on every slab
$I_j\times\R^n$, we get
$$
\|u\|_{Y^1(I_j)}\le 2\eta.
$$
Adding these estimates over all intervals $I_j$ and recalling that the total number of these intervals is $J=J(E,\eta)$,
we get
$$
\|u\|_{Y^1(\R^+)}\le C(E).
$$
Strichartz estimates and time reversal symmetry (see the end of subsection 5.6) yield the global spacetime bound
$$
\|u\|_{S^1(\R\times\R^n)}\le C(E).
$$

\subsection{Finite global Strichartz norms imply scattering}

In this subsection we show that finite global Strichartz norms imply scattering.  As the techniques are classical, we
will only present the construction of the scattering states $u_{\pm}$ and show that their linear flow approximates
well the solution as $t\to \pm \infty$ in $H^1_x$.  Standard arguments can also be employed to construct the wave operators;
for details see, for example, \cite{cazenave:book}.

For simplicity, we will only construct the scattering state in the positive time direction. Similar arguments can
be used to construct the scattering state in the negative time direction.

For $0<t<\infty$, we define
$$
u_+(t)=u_0-i\int_0^t e^{-is\Delta}(\ntu)(s) ds.
$$
As $u\in S^1(\R\times\R^n)$, Strichartz estimates and Lemma~\ref{lemma bound by w01} imply that $u_+(t)\in H_x^1$
for all $t\in \R^+$.  We will show that $u_+(t)$ converges in $H^1_x$ as $t\to \infty$.  Indeed, for $0<\tau<t$, we have
\begin{align*}
\|u_+(t)-u_+(\tau)\|_{H_x^1}
&=\Bigl\|\int_{\tau}^{t}e^{-is\Delta}(\ntu)(s)ds\Bigr\|_{H_x^1}\\
&\lesssim \Bigl \|\int_{\tau}^{t}e^{i(t-s)\Delta}(\ntu)(s)ds\Bigr\|_{L_t^{\infty}H_x^1([\tau,t]\times\R^n)},
\end{align*}
which by Strichartz and \eqref{bound by w01} implies
$$
\|u_+(t)-u_+(\tau)\|_{H_x^1}
\lesssim C(\lambda_1,\lambda_2)\sum_{i=1,2}\|u\|_{V([\tau,t])}^{2-\frac{p_i(n-2)}2}\|u\|_{W([\tau,t_])}^{\frac{np_i}2-2} \|(1+|\nabla|)u\|_{V([\tau,t])}.
$$
As
\begin{align}\label{s1 bound 5.8}
\|u\|_{S^1(\R\times\R^n)}<\infty,
\end{align}
we see that for $\eps>0$ there exists $T_\eps>0$ such that
$$
\|u_+(t)-u_+(\tau)\|_{H_x^1}\leq \eps
$$
for any $t, \tau>T_\eps$.  Thus, $u_+(t)$ converges in $H^1_x$ as $t\to\infty$ to some function $u_+$.
Moreover, Strichartz and \eqref{bound by w01} imply that
$$
u_+:=u_0-i\int_0^{\infty} e^{-is\Delta}(\ntu)(s) ds.
$$

Next, we show that the linear evolution $e^{it\Delta}u_+$ approximates $u(t)$ in $H_x^1$ as $t\to \infty$.
Indeed, by Strichartz and Lemma~\ref{lemma bound by w01},
\begin{align*}
\|e^{-it\Delta}u(t)-u_+\|_{H_x^1}
&=\Bigl\|\int_t^{\infty} e^{-is\Delta}(\ntu)(s)ds\Bigr\|_{H_x^1}\\
&=\Bigl\|\int_t^{\infty}e^{i(t-s)\Delta}(\ntu)(s)ds\Bigr\|_{H^1_x}\\
&\lesssim\sum_{i=1,2}\|u\|_{V([t,\infty))}^{2-\frac{p_i(n-2)}2}\|u\|_{W([t,\infty))}^{\frac{np_i}2-2}\|(1+|\nabla|)u\|_{V([t,\infty))},
\end{align*}
which obviously tends to $0$ as $t\to\infty$ (see \eqref{s1 bound 5.8}).

%
%
%
%

\section{Finite time blowup}

To prove blowup under the assumptions on the nonlinearities described in Theorem~\ref{blowup}, we will follow
the convexity method of Glassey, \cite{glassey}.  More precisely, we will consider the variance
$$
V(t)=\int_{\R^n} |x|^2|u(t,x)|^2dx
$$
and show that as a function of $t>0$, it is decreasing and concave, which suggests the existence of a blowup time $T_*$
at which $V(T_*)=0$.

For strong $H^1_x$-solutions $u$ to \eqref{equation} with initial data $u_0\in \Sigma$, it is known that
$V\in C^2(-T_{min},T_{max})$; see, for example, \cite{cazenave:book}.  Formal computations (which are made
rigorous in \cite{cazenave:book}) prove

\begin{lemma} For all $t\in(-T_{min},T_{max})$, we have
\begin{align}\label{Vdot}
V'(t)=-4y(t),
\end{align}
where
$$
y(t)=-\Im\int_{\R^n} r\bar uu_r(t)dx.
$$
Moreover,
\begin{align}\label{Vdotdot}
V''(t)=-4y'(t)=8\|\nabla u(t)\|_2^2+\tfrac{4n\lambda_1p_1}{p_1+2}\|u(t)\|_{p_1+2}^{p_1+2}
  +\tfrac{4n\lambda_2p_2}{p_2+2}\|u(t)\|_{p_2+2}^{p_2+2}.
\end{align}
\end{lemma}

In what follows, we will show that the first and second derivatives of the variance are negative for positive times $t$.
More precisely, in each of the three cases described in Theorem ~\ref{blowup}, we will show that
\begin{equation}\label{ydot}
y'(t)\ge c\|\nabla u(t)\|_2^2>0
\end{equation}
for a small positive constant $c$.  Thus, by \eqref{Vdotdot} it follows that $V''(t)<0$ for all times
$t\in(-T_{min},T_{max})$, which implies that $V(t)$ is concave. Moreover, as by hypothesis
$$
y(0)=y_0>0,
$$
\eqref{ydot} implies that $y(t)>y(0)>0$, for all times $t>0$.  By \eqref{Vdot}, this yields $V'(t)<0$
for positive times $t$ and hence, $V(t)$ is decreasing for $t>0$.

Alternatively, by H\"older, for all times $t\in(-T_{min},T_{max})$,
$$
y(t)\le \|xu(t)\|_2\|\nabla u(t)\|_2
$$
and hence,
\begin{equation}\label{key}
\|\nabla u(t)\|_2\ge \frac{y(t)}{\|xu_0\|_2}.
\end{equation}
By \eqref{ydot} and \eqref{key}, we derive a first order ODE for $y$,
\begin{equation}
\begin{cases}
y'(t)\geq c\frac{y^2(t)}{\|xu_0\|_2^2}\\
y(0)=y_0>0,
\end{cases}
\end{equation}
which implies that there exists $0<T_*\le \frac{\|xu_0\|_2^2}{c y_0}$ such that
$$
\lim_{t\to T_*}\|\nabla u(t)\|_2=\infty.
$$

For the remainder of the section, we will derive \eqref{ydot} in each of the three cases described in Theorem~\ref{blowup}.

\textbf{Case 1):} $\lambda_1>0$ , $0<p_1<p_2$, and $E<0$. \\
By \eqref{Vdotdot}, the conservation of energy, and our assumptions, we get
\begin{align*}
y'(t)
&=-2\|\nabla u(t)\|_2^2 + p_2n\Bigl\{\tfrac 12\|\nabla u(t)\|_2^2 + \tfrac{\lambda_1}{p_1+2}\|u(t)\|_{p_1+2}^{p_1+2}-E\Bigr\}
  -\tfrac{n\lambda_1p_1}{p_1+2}\|u(t)\|_{p_1+2}^{p_1+2}\\
&=\tfrac {p_2n-4}2\|\nabla u(t)\|_2^2 + \Bigl\{\tfrac{p_2n\lambda_1}{p_1+2}-\tfrac{p_1n\lambda_1}{p_1+2}\Bigr\}\|u(t)\|_{p_1+2}^{p_1+2}-p_2nE\\
&\ge \tfrac {p_2n-4}2\|\nabla u(t)\|_2^2 + \tfrac{n\lambda_1(p_2-p_1)}{p_1+2}\|u(t)\|_{p_1+2}^{p_1+2}\\
&\ge \tfrac{p_2n-4}2\|\nabla u(t)\|_2^2,
\end{align*}
and hence \eqref{ydot} holds with $c:=\frac{p_2n-4}2$.

\textbf{Case 2):} $\lambda_1<0$, $\frac 4n<p_1<p_2$, and $E<0$. \\
In this case, by \eqref{Vdotdot} and the conservation of energy, we obtain
\begin{align*}
y'(t)
&=-2\|\nabla u(t)\|_2^2 + p_1n\Bigl\{\tfrac 12\|\nabla u(t)\|_2^2+\tfrac {\lambda_2}{p_2+2}\|u(t)\|_{p_2+2}^{p_2+2}-E\Bigr\}
  -\tfrac{p_2n\lambda_2}{p_2+2}\|u(t)\|_{p_2+2}^{p_2+2}\\
&=\tfrac{p_1n-4}2\|\nabla u(t)\|_2^2+\tfrac{n\lambda_2(p_1-p_2)}{p_2+2}\|u(t)\|_{p_2+2}^{p_2+2}-p_1nE\\
&\ge \tfrac{p_1n-4}2\|\nabla u(t)\|_2^2,
\end{align*}
which implies \eqref{ydot} with $c:=\frac{p_1n-4}2$.

\textbf{Case 3):} $\lambda_1<0$, $ 0<p_1\le \frac 4n$, and $E+C M<0$ for some constant
$C=C(\lambda_1,\lambda_2,p_1,p_2,n)$ to be specified momentarily.

The idea is to use part of the contribution coming from the higher power nonlinearity to obtain a positive multiple
of the kinetic energy and the rest to annihilate the effect of the resulting lower power term.  The details are
as follows.

As $p_2>\frac 4n$, we can find a small constant $\eps$ such that $p_2>\frac{2(2+\varepsilon)}n$.  It is immediate
that $\theta:=\frac {2(2+\eps)}{p_2n}<1$.  Therefore, from the conservation of energy, we get
\begin{align*}
y'(t)
&=-2\|\nabla u(t)\|_2^2- \tfrac{n\lambda_2p_2\theta}{p_2+2}\|u(t)\|_{p_2+2}^{p_2+2}-\tfrac{n\lambda_2p_2(1-\theta)}{p_2+2}\|u(t)\|_{p_2+2}^{p_2+2}
 -\tfrac{n\lambda_1p_1}{p_1+2}\|u\|_{p_1+2}^{p_1+2}\\
&\ge -2\|\nabla u(t)\|_2^2 + p_2n\theta\Bigl\{\tfrac 12\|\nabla u(t)\|_2^2 + \tfrac{\lambda_1}{p_1+2}\|u(t)\|_{p_1+2}^{p_1+2}-E\Bigr\}\\
&\quad-\tfrac{n\lambda_2p_2(1-\theta)}{p_2+2}\|u(t)\|_{p_2+2}^{p_2+2} -\tfrac{n\lambda_1p_1\theta}{p_1+2}\|u(t)\|_{p_1+2}^{p_1+2}\\
&\ge \Bigl\{-2+\tfrac{p_2n\theta}2\Bigr\}\|\nabla u(t)\|_2^2 + \tfrac{n\lambda_1\theta(p_2-p_1)}{p_1+2}\|u(t)\|_{p_1+2}^{p_1+2}-p_2n\theta E\\
&\quad-\tfrac{n\lambda_2p_2(1-\theta)}{p_2+2}\|u(t)\|_{p_2+2}^{p_2+2}\\
&=\eps \|\nabla u(t)\|_2^2+\tfrac{n\lambda_1\theta(p_2-p_1)}{p_1+2}\|u(t)\|_{p_1+2}^{p_1+2}
 - \tfrac{n\lambda_2p_2(1-\theta)}{p_2+2}\|u(t)\|_{p_2+2}^{p_2+2}-p_2n\theta E.
\end{align*}

By Young's inequality, for any positive constants $a$ and $\delta$,
$$
a^{p_1+2}\le C(\delta) a^2+\delta a^{p_2+2}.
$$
Hence,
\begin{align*}
\tfrac{n|\lambda_1|\theta(p_2-p_1)}{p_1+2}&\|u(t)\|_{p_1+2}^{p_1+2}
\le C(\delta)\tfrac{n|\lambda_1|\theta(p_2-p_1)}{p_1+2}\|u(t)\|_2^2+\delta\tfrac{n\theta|\lambda_1|(p_2-p_1)}{p_1+2}\|u(t)\|_{p_2+2}^{p_2+2}.
\end{align*}
Choosing $\delta>0$ sufficiently small such that
$$
\delta\frac{n\theta|\lambda_1|(p_2-p_1)}{p_1+2}<\frac{np_2|\lambda_2|(1-\theta)}{p_2+2},
$$
we obtain
$$
y'(t)\ge \eps\|\nabla u(t)\|_2^2-C(\lambda_1,\lambda_2,p_1,p_2,n)M-p_2n\theta E,
$$
which, as long as
$$
p_2n\theta E+C(\lambda_1,\lambda_2,p_1,p_2,n)M<0,
$$
yields
$$
y'(t)\ge \eps\|\nabla u(t)\|_2^2.
$$
This proves \eqref{ydot} in this case.

Hence, \eqref{ydot} holds in all three cases described in Theorem~\ref{blowup}, which implies for the reasons
given above the existence of a time $T_*\leq \frac{\|xu_0\|_2^2}{y_0}$ such that
$$
\lim_{t\to T_*}\|\nabla u(t)\|_{2}=\infty.
$$

%
%
%
%

\section{Scattering in $\Sigma$}
In this section we prove Theorem~\ref{sigmascattering}.  When $\alpha(n)<p_1<p_2<\frac 4{n-2}$, the result is
well-known (see \cite{tsutsumi:scatter, hayashi:tsutsumi, cazenave:book}) and we will not revisit the proof here.
Instead, we will present the proof of Theorem~\ref{sigmascattering} in the case $\alpha(n)<p_1<p_2=\nmt$, which
is not covered by these earlier results.

As both nonlinearities are assumed defocusing, without loss of generality we may assume $\lambda_1=\lambda_2=1$.
We will also write $p$ instead of $p_1$ to ease notation.  Note that under our assumptions,
Theorem~\ref{global wellposedness} implies the existence of a unique global solution $u$.  Moreover, on every slab $\ir$
$u$ obeys the following estimate:
\begin{align}\label{local stb}
\|u\|_{S^1(I\times\R^n)}\le C(\|u_0\|_{H_x^1},|I|).
\end{align}

Next, we introduce the Galilean operator associated with the linear Schr\"odinger operator $i\partial_t+\Delta$,
$$
H(t):=x+2it\nabla.
$$
Note that we can also write
\begin{align}\label{H(t)}
H(t)=2ite^{\frac{i|x|^2}4t}\nabla(e^{-\frac{i|x|^2}4t}\cdot)=e^{it\Delta}xe^{-it\Delta}.
\end{align}
From the first identity in \eqref{H(t)} we see that $H(t)$ behaves morally like a derivative when applied to the Hamiltonian nonlinearity $F(u)=|u|^pu+|u|^{\nmt} u$, which commutes with phase rotations.  More precisely, we have
\begin{align}\label{dH(t)}
H(t)F(u)=\partial_z F(u)H(t)u-\partial_{\bar z}F(u)\overline{H(t)u}.
\end{align}

The first step to prove scattering in $\Sigma$ is to use the local spacetime bound \eqref{local stb} and the fact that
$xu_0\in L^2_x$ to derive $\dot S^0$ spacetime bounds on $Hu$ on every slab $\ir$.
By time reversal symmetry, we may assume $I=[0,T]$.  By \eqref{local stb}, we can split $[0,T]$ into
$J=J(\|u_0\|_{H^1_x}, T, \eta)$ subintervals $I_j=[t_j,t_{j+1}]$ such that
\begin{align}\label{scatsigma ux1}
\|u\|_{\dot X^1(I_j)}\sim \eta,
\end{align}
where $\eta>0$ is a small constant to be chosen later.  We will derive Strichartz bounds on $Hu$ on every slab
$I_j\times\R^n$.

Fix therefore $I_j$; on this interval $u$ satisfies
$$
u(t)=e^{i(t-t_j)\Delta}u(t_j)-i\int_{t_j}^t e^{i(t-s)\Delta}(|u|^pu+|u|^{\frac 4{n-2}}u)(s)ds.
$$
By \eqref{H(t)}, this yields
$$
H(t)u(t)=e^{i(t-t_j)\Delta}H(t_j)u(t_j)-i\int_{t_j}^t e^{i(t-s)\Delta}H(s)(|u|^pu+|u|^{\frac4{n-2}}u)(s)ds,
$$

By Strichartz, Lemma~\ref{control by x}, \eqref{dH(t)}, and \eqref{scatsigma ux1}, we estimate
\begin{align*}
\|Hu\|_{\snij}
&\lesssim \|H(t_j)u(t_j)\|_{L_x^2}+|I_j|^{1-\frac{p(n-2)}4}\|u\|_{\dot X^1(I_j)}^p\|Hu\|_{\snij}\\
&\quad +\|u\|_{\dot X^1(I_j)}^{\nmt}\|Hu\|_{\snij}\\
&\lesssim \|H(t_j)u(t_j)\|_{L_x^2}+T^{1-\frac{p(n-2)}4}\eta^p\|Hu\|_{\snij}\\
&\quad +\eta^{\nmt}\|Hu\|_{\snij},
\end{align*}
which implies
$$
\|Hu\|_{\snij}\lesssim \|H(t_j)u(t_j)\|_{L_x^2},
$$
provided $\eta$ is chosen sufficiently small depending on $|I|=T$. By induction, for each $j$ we get
$$
\|Hu\|_{\snij}\lesssim \|xu_0\|_{L_x^2}.
$$
Therefore, adding these estimates over all subintervals $I_j$, we obtain
\begin{align}\label{Hu s0 local}
\|Hu\|_{\dot S^0(I\times\R^n)}\le C(\|xu_0\|_{L_x^2}, |I|).
\end{align}

In order to prove scattering, rather than a local spacetime bound, one needs to derive global spacetime bounds.
To accomplish this, we will use the pseudoconformal identity to prove that the global solution $u$ decays as $t\to\infty$.
More precisely, we will show that $t\mapsto\|u(t)\|_{\frac {2n} {n-2}}$ and $t\mapsto\|u(t)\|_{p+2}$ are decreasing.

Introduce the pseudoconformal energy
$$
h(t):=\|H(t)u(t)\|_2^2+8t^2\bigl(\tfrac 1{p+2}\|u(t)\|_{p+2}^{p+2}+\tfrac {n-2}{2n}\|u(t)\|_{\dnnt}^{\dnnt}\bigr).
$$
A standard computation (see e.g. \cite{SulSul}) shows that
\begin{equation}\label{5.3}
h'(t)=t\bigl(\tfrac {4(4-pn)}{p+2}\|u(t)\|_{p+2}^{p+2}-\tfrac {16}n\|u(t)\|_{\dnnt}^{\dnnt}\bigr):=t\theta(t).
\end{equation}
Again, we can justify this computation in the class $\Sigma$ (treating $h'(t)$ as a weak derivative of $h(t)$) by first
mollifying $u_0$ and the nonlinearity and then taking limits at the end using the well-posedness theory; we omit
the standard details.

If $\frac 4n\le p<\nmt$, we integrate \eqref{5.3} with respect to the time variable over the compact interval $[0,t]$
to obtain
\begin{align}\label{5.4}
\|H(t)&u(t)\|_2^2+\tfrac {8t^2}{p+2}\|u(t)\|_{p+2}^{p+2}+\tfrac{4(n-2)t^2}n\|u(t)\|_{\dnnt}^{\dnnt}\\
&=\|xu_0\|_2^2+\tfrac {4(4-pn)}{p+2}\int_0^ts\int_{\R^n}|u(s,x)|^{p+2} dxds- \tfrac{16}n\int_0^t s\int_{\R^n}|u(s,x)|^{\dnnt} dxds.  \nonumber
\end{align}
Thus, in this case,
\begin{align}\label{decay 1}
\|u(t)\|_{p+2}^{p+2}+\|u(t)\|_{\dnnt}^{\dnnt}\lesssim t^{-2}.
\end{align}

If $0<p<\frac 4n$, we integrate \eqref{5.3} over the compact interval $[1,t]$; we get
$$
h(t)= h(1)+\int_1^t s\theta(s)ds, \quad \mbox{for  } t\ge 1.
$$
In particular,
$$
\tfrac {8t^2}{p+2}\|u(t)\|_{p+2}^{p+2} \le h(1)+\tfrac {4(4-pn)}{p+2}\int_1^ts\|u(s)\|_{p+2}^{p+2}ds,
$$
which by Grownwall's inequality implies
$$
\|u(t)\|_{p+2}^{p+2}\lesssim t^{-\frac {pn}2}.
$$
Combining this with \eqref{5.4}, we obtain
$$
h(t)\lesssim 1+t^{2-\frac {pn}2}
$$
and hence,
$$
\|u(t)\|_{\dnnt}^{\dnnt}\lesssim t^{-\frac {pn}2}, \quad \forall t\ge 1.
$$
Thus, for $0<p<\frac 4n$ and $t\geq 1$ we get
\begin{align}\label{decay 2}
\|u(t)\|_{p+2}^{p+2}+\|u(t)\|_{\dnnt}^{\dnnt}\lesssim t^{-\frac {pn}2}.
\end{align}
Collecting \eqref{decay 1} and \eqref{decay 2}, we obtain
\begin{align}\label{decay}
\|u(t)\|_{p+2}^{p+2}+\|u(t)\|_{\dnnt}^{\dnnt}\lesssim t^{-2} + t^{-\frac {pn}2}
\end{align}
for all $0<p<\frac {4}{n-2}$ and $t\geq 1$.

Next, we show that the decay estimate \eqref{decay} implies that the Strichartz norms of the solution
are small when measured over the slab $[T, \infty)\times\R^n$ for a large enough time $T$.
Let $\delta$ be such that $(\delta, p+2)$ is a Schr\"odinger admissible pair. An easy computation shows that
$p>\alpha(n)\Leftrightarrow 2p>\delta-2$. Then, by \eqref{decay}, on the slab $[T, \infty)\times\R^n$ with
$T\geq 1$ we get
\begin{align}
\||u|^p u\|_{L_t^{\delta'}W_x^{1,\frac{p+2}{p+1}}}
&\lesssim \|u\|_{L_t^{\delta}W_x^{1,p+2}}\|u\|^p_{L_t^{\frac{\delta p}{\delta-2}}L_x^{p+2}}
 \lesssim T^{-\nu}\|u\|_{S^1} \label{sigma up}\\
\||u|^{\frac{4}{n-2}}u\|_{L^2_tW_x^{1,\frac {2n}{n+2}}}
&\lesssim \|u\|^{\frac{4}{n-2}}_{L^{\infty}_tL_x^{\dnnt}}\|u\|_{L^2_tW_x^{1,\dnnt}}
 \lesssim T^{-\nu}\|u\|_{S^1} \label{sigma ucrit}
\end{align}
for some $\nu>0$.  As by \eqref{duhamel}, on $[T,\infty)$ $u$ satisfies
$$
u(t)=e^{i(t-T)\Delta}u(T)-i\int_T^{\infty}e^{i(t-s)\Delta}(|u|^pu+|u|^{\nmt})u(s)ds,
$$
by Strichartz, \eqref{sigma up}, and \eqref{sigma ucrit}, we get
$$
\|u\|_{S^1([T,\infty)\times\R^n)}
\lesssim \|u(T)\|_{H_x^1}+T^{-\nu}\|u\|_{S^1([T,\infty)\times\R^n)}.
$$
Hence, taking $T$ sufficiently large, by the conservation of energy and mass we obtain
\begin{align}\label{u s1 far}
\|u\|_{S^1([T,\infty)\times\R^n)}\lesssim \|u(T)\|_{H_x^1} \le C(\|u_0\|_{H_x^1}).
\end{align}

Similarly, by Strichartz, \eqref{H(t)}, \eqref{dH(t)}, \eqref{sigma up}, and \eqref{sigma ucrit}, we estimate
$$
\|Hu\|_{\dot S^0([T,\infty)\times\R^n)}
\lesssim \|H(T)u(T)\|_{L_x^2}+T^{-\nu}\|Hu\|_{\dot S^0([T,\infty)\times\R^n)},
$$
which taking $T$ sufficiently large and using \eqref{Hu s0 local} yields
\begin{align}\label{Hu s0 far}
\|Hu\|_{\dot S^0([T,\infty)\times\R^n)}\lesssim \|H(T)u(T)\|_{L_x^2}\leq C(\|xu_0\|_{L_x^2}).
\end{align}

Combining \eqref{local stb} and \eqref{u s1 far}, \eqref{Hu s0 local} and \eqref{Hu s0 far}, and using the time reversal
symmetry, we get
\begin{align}
\|u\|_{S^1(\R\times\R^n)}&\le C(\|u_0\|_{H^1_x}) \label{scat us1 global}\\
\|Hu\|_{\dot S^0(\R\times\R^n)}&\le C(\|u_0\|_{\Sigma}). \label{scat Hu s0 global}
\end{align}

Next, we construct the scattering state in the positive time direction.   The scattering
state in the negative time direction is constructed similarly.  Let
$$
u_+(t)=u_0-i\int_0^t e^{-is\Delta}(|u|^pu+|u|^{\nmt}u)(s)ds.
$$
We will show that $u_+(t)$ is a Cauchy sequence in $\Sigma$ when $t\to \infty$. Take $t_1<t_2<\infty$.
Then, by Strichartz, \eqref{H(t)}, the fact that $e^{it\Delta}$ is unitary on $L_x^2$, \eqref{sigma up},
and \eqref{sigma ucrit}, we estimate
\begin{align*}
\|u_+(t_1)-u_+(t_2)\|_{\Sigma}
&=\Bigl\|\int_{t_1}^{t_2} e^{-is\Delta} (|u|^pu+|u|^{\nmt}u)(s)ds\Bigr\|_{\Sigma}\\
&\lesssim\Bigl \|\int_{t_1}^{t_2}e^{-is\Delta}(|u|^pu+|u|^{\nmt}u)(s)ds\Bigr\|_{H_x^1}\\
&\quad +\Bigl\|x\int_{t_1}^{t_2}e^{-is\Delta}(|u|^pu+|u|^{\nmt}u)(s)ds\Bigr\|_{L_x^2}\\
&\lesssim \Bigl\|\int_{t_1}^{t_2}e^{i(t-s)\Delta}(|u|^pu+|u|^{\nmt}u)(s)ds\Bigr\|_{H_x^1}\\
&\quad +\Bigl\|\int_{t_1}^{t_2}e^{i(t-s)\Delta}H(s)(|u|^pu+|u|^{\nmt}u)(s)ds\Bigr\|_{L_x^2}\\
&\lesssim t_1^{-\nu}\bigl(\|u\|_{S^1([t_1,t_2]\times\R^n)}+\|Hu\|_{\dot S^0([t_1,t_2]\times\R^n)}\bigr)\\
&\le t_1^{-\nu} C(\|u_0\|_{\Sigma}).
\end{align*}
Therefore, $u_+(t)$ is a Cauchy sequence in $\Sigma$ as $t\to \infty$ and hence, it converges to some function
$u_+\in \Sigma$.  Similar estimates show that
$$
u_+=u_0-i\int_0^{\infty} e^{-is\Delta}(|u|^pu+|u|^{\nmt}u)(s)ds.
$$

Finally, we show that $u_+$ is the asymptotic state of $e^{-it\Delta}u(t)$ as $t\to \infty$.  Indeed,
similar computations as above yield
\begin{align*}
\|e^{-it\Delta}u(t)-u_+\|_{\Sigma}
&\leq \Bigl\|\int_t^{\infty}e^{-is\Delta}(|u|^pu+|u|^{\nmt}u)(s)ds\Bigr\|_{\Sigma}\\
&\lesssim t^{-\nu}\bigl(\|u\|_{S^1(\R\times\R^n)}+\|Hu\|_{\dot S^0(\R\times\R^n)}\bigr)\\
&\lesssim t^{-\nu} C(\|u_0\|_{\Sigma})
\end{align*}
which tends to 0 as $t\to\infty$, as desired.


\begin{thebibliography}{10}
\bibitem{Blowscatter}
J. Bourgain, \emph{Scattering in the energy space and below for 3D NLS}, Journal D'Analyse
Mathematique,\textbf{75} (1998), 267-297.

\bibitem{borg:scatter}
J. Bourgain, \emph{Global well-posedness of defocusing 3D critical NLS in the radial case}, JAMS
\textbf{12} (1999), 145--171.

\bibitem{borg:book}
J. Bourgain, \emph{New global well-posedness results for non-linear Schr\"odinger equations}, AMS
Publications (1999).

\bibitem{cwI}
T. Cazenave, F.B. Weissler, \emph{Critical nonlinear Schr\"odinger Equation}, Non. Anal. TMA \textbf{14}
(1990), 807--836.

\bibitem{cazenave:book}
T. Cazenave, \textit{Semilinear Schr\"odinger equations,} Courant Lecture Notes in Mathematics, 10.
American Mathematical Society, 2003.

\bibitem{ckstt:gwp}
J. Colliander, M. Keel, G. Staffilani, H. Takaoka, T. Tao, \emph{Global well-posedness and scattering
for the energy-critical nonlinear Schr\"odinger equation in $\R^3$}, to appear Annals of Math.

\bibitem{ckstt:french}
J. Colliander, M. Keel, G. Staffilani, H. Takaoka, T. Tao, \emph{Existence globale et diffusion pour
l'\'equation de
  Schr\"odinger nonlin\'eaire r\'epulsive cubique sur $\R^3$ en
  dessous l'espace d'\'energie},
Journ\'ees ``\'Equations aux D\'eriv\'ees Partielles''
              (Forges-les-Eaux, 2002), Exp. No. X, 14, 2002.

\bibitem{glassey}
R. T. Glassey, \emph{On the blowing up of solution to the Cauchy problem for nonlinear Schr\"odinger
operators}, J. Math. Phys. \textbf{8} (1977), 1794--1797.

\bibitem{gv:scatter}
J. Ginibre, G. Velo, \emph{Scattering theory in the energy space for a class of nonlinear Schr\"odinger
equations}, J. Math. Pure. Appl. \textbf{64} (1985), 363--401.

\bibitem{grillakis:scatter}
M. Grillakis, \emph{On nonlinear Schr\"odinger equations.}, Comm. Partial Differential Equations
\textbf{25} no. 9-10 (2000), 1827--1844.

\bibitem{hayashi:tsutsumi}
N. Hayashi, Y. Tsutsumi, \emph{Remarks on the scattering problem for nonlinear Schr\"odinger equations},
Diff. Eq. Math. Phys., (1986), 162--168.

\bibitem{kato}
T. Kato, \emph{On nonlinear Schr\"odinger equations},  Ann. Inst. H. Poincare Phys. Theor.
\textbf{46}  (1987),  113--129.

\bibitem{katounique}
T. Kato, \emph{On nonlinear Schr\"odinger equations, II.  $H^s$-solutions and unconditional well-posedness}, J. d'Analyse. Math. \textbf{67}, (1995), 281--306.

\bibitem{tao:keel}
M. Keel, T. Tao, \emph{Endpoint Strichartz Estimates}, Amer. Math. J. \textbf{120} (1998), 955--980.

\bibitem{linstrauss}
J. Lin, W. Strauss, \emph{Decay and scattering of solutions of a nonlinear Schr\"odinger equation},
Journ. Funct. Anal. \textbf{30} (1978), 245--263.

\bibitem{merle1}
F. Merle, \emph{Determination of blow-up solutions with minimal mass for nonlinear Schr\"odinger equations
with critical power}, Duke Math. J. 69 (1993), 427--454.

\bibitem{merle-raphael1}
F. Merle, P. Raphael, \emph{The blow-up dynamic and upper bound on the blow-up rate for critical
nonlinear Schr\"odinger equation}, Ann. of Math. 161 (2005), 157--222.

\bibitem{merle-raphael2}
F. Merle, P. Raphael, \emph{Profiles and quantization of the blow up mass for critical nonlinear
Schr\"odinger equation}, Comm. Math. Phys. 253 (2005), 675--704.

\bibitem{merle-tsutsumi}
F. Merle, Y. Tsutsumi, \emph{$L^2$-concentration of blow-up
solutions for the nonlinear Schr\"odinger equation with the critical
power nonlinearity}, J. Diffe. Equa. 84 (1990), 205-214.

\bibitem{morawetz}
C. Morawetz, \emph{Time decay for the nonlinear Klein-Gordon equation}, Proc. Roy. Soc. A \textbf{306}
(1968), 291--296.

\bibitem{nak:scatter}
K. Nakanishi, \emph{Energy scattering for non-linear Klein-Gordon and Schr\"odinger equations in spatial
dimensions 1 and 2}, JFA \textbf{169} (1999), 201--225.

\bibitem{nakanishi}
K. Nakanishi, \emph{Scattering theory for nonlinear Klein-Gordon equation with Sobolev critical power},
IMRN \textbf{1} (1999), 31--60.

\bibitem{RV}
E.~Ryckman, M.~Visan \emph{Global well-posedness and scattering for the defocusing energy-critical
nonlinear Schr\"odinger equation in $\R^{1+4}$,} preprint {\tt math.AP/0501462}.

\bibitem{SulSul}
C. Sulem, P.-L. Sulem, \emph{The nonlinear Schr\"odinger equation: Self-focusing and wave collapse}, Springer-Verlag, New York, (1999), p. 25.

\bibitem{tao:focusing}
T. Tao, \emph{On the asymptotic behavior of large radial data for a
focusing non-linear Schr\"odinger equation,} preprint {\tt
math.AP/0309428}.

\bibitem{tao: gwp radial}
T. Tao, \emph{Global well-posedness and scattering for the
higher-dimensional energy-critical non-linear Schr\"odinger equation
for radial data}, New York Journal of Math., \textbf{11} (2005),
57-80.

\bibitem{TV}
T. Tao, M.~Visan \emph{Stability of energy-critical nonlinear Schr\"odinger equations in high
dimensions,} Electron. J. Diff. Eqns., \textbf{118} (2005), 1-28.


\bibitem{tsutsumi:scatter}
Y. Tsutsumi, \emph{Scattering problem for nonlinear Schr\"odinger equations},
Ann. Inst. H. Poincaré Phys. Théor. 43 (1985), 321--347.

\bibitem{thesis:art} M. Visan, \emph{The defocusing energy-critical nonlinear Schr\"odinger equation
in higher dimensions,} preprint \texttt{math.AP/0508298}.

\bibitem{Monica:thesis} M. Visan, \emph{The defocusing energy-critical nonlinear Schr\"odinger equation
in dimensions four and higher,} Ph.D. Thesis, UCLA, in preparation.

\bibitem{xzhang} X. Zhang, \emph{On Cauchy problem of 3-D energy critical Schr\"odinger equation with
subcritical perturbations,} preprint.

\end{thebibliography}
\end{document}